\documentclass[fleqn,12pt]{article}
\usepackage{amssymb,amsmath,amsfonts}
\usepackage{dsfont}
\usepackage{color}
\usepackage{graphicx}
\usepackage{graphics}
\usepackage{latexsym}
\usepackage{amsmath}
\usepackage{amssymb}
\usepackage{graphics}
\usepackage[dvips]{epsfig}
\usepackage{mathrsfs}

\usepackage{cases}

\newcommand{\dt}{\partial_t}
\newcommand{\dv}{\mathrm{div}\,}
\newcommand{\dr}{\partial_r}

\newcommand{\curl}{\mathrm{curl}\,}
\newcommand{\tr}{\mathrm{tr}\,}
\newcommand{\dT}{\nabla_T}
\newcommand{\dN}{\nabla_N}

\newcommand{\iv}{\tilde v}
\newcommand{\iq}{\tilde q}

\newcommand{\subeqref}[2]{$\eqref{#1}_{#2}$}
\newcommand{\abs}[2]{\left| #1 \right|^{#2}}
\newcommand{\norm}[2]{\left\Arrowvert #1 \right\Arrowvert_{#2}}
\newcommand{\Lnorm}{L^2(\Omega)}
\newcommand{\Hnorm}[1]{H^{#1}(\Omega)}
\newcommand{\bLnorm}{L^2(\Gamma)}
\newcommand{\bHnorm}[1]{H^{#1}(\Gamma)}
\newcommand{\supnorm}{L^{\infty}(\Omega)}
\newcommand{\inti}[1]{\int_\Omega #1 \, dx}
\newcommand{\intb}[1]{\int_\Gamma #1 \, dS}

\newcommand{\commut}[2]{\left\lbrack #1, #2 \right\rbrack}

\newtheorem{thm}{Theorem}[section]
\newtheorem{lm}{Lemma}
\newtheorem{cor}{Corollary}
\newtheorem{prop}{Proposition}

\newenvironment{pf}{\paragraph{Proof}}{\hfill$\square$}

\newenvironment{rmk}{\paragraph{Remark}}{\hfill}

\numberwithin{equation}{section}
\allowdisplaybreaks
\title{Compressible Viscous Flows in a Symmetric Domain with Complete Slip Boundary}

\author{Xin Liu}

\begin{document}
\maketitle

\begin{abstract}
This work is devoted to study the global behavior of viscous flows contained in a symmetric domain with complete slip boundary. In such scenario the boundary no longer provides friction and therefore the perturbation of angular velocity lacks decaying structure. In fact, we show the existence of uniformly rotating solutions as steady states for the compressible Navier-Stokes equations. By manipulating the conservation law of angular momentum, we establish a suitable Korn's type inequality to control the perturbation and show the asymptotic stability of the uniformly rotating solutions with small angular velocity. In particular, the initial perturbation which preserves the angular momentum would decay exponentially in time and the solution to the Navier-Stokes equations converges to the steady state as time grows up. 
\end{abstract}

\tableofcontents

\section{Introduction}
\subsection{Description and Related Works}
In this work, we consider the isentropic compressible Navier-Stokes system in the following, which models the motion of viscous gases (or fluids) in a bounded domain $ \Omega \subset \mathbb R^3 $,
\begin{equation}\label{ICNS}
	\begin{cases}
		\dt \rho + \dv (\rho u) = 0 & \text{in} ~ \Omega, \\
		\dt (\rho u ) + \dv ( \rho u \otimes u ) + \nabla P = \dv \mathbb{S}(u) & \text{in} ~ \Omega,
	\end{cases}
\end{equation}
where $ \rho, u, P, \mathbb S(u) $ represent the density, the velocity, the pressure potential, and the viscous tensor respectively. Moreover, the flow is assumed to be Newtonian. For simplicity, the pressure potential and the viscous tensor are taken in the following forms,
\begin{equation*}
	P = P(\rho) = \rho^\gamma, ~ \mathbb S(u) = \mu \left( \nabla u + \nabla u^\top \right) + \lambda \dv u  \mathbb{I}_3, ~ \gamma > 1, ~ \mu, \lambda > 0,
\end{equation*}
where $ \mathbb{I}_3 $ is the $3 \times 3$ identity matrix and $ \mu, \lambda $ are the viscosity coefficients. \eqref{ICNS} can be complemented with various boundary conditions. In this work, the associated boundary condition is taken as the complete slip boundary condition, i.e., 
\begin{equation}\label{SBC}
	u \cdot \vec{n} = 0, ~ \vec{\tau} \cdot \mathbb S(u) \vec{n} = 0, ~~~ \text{on} ~ \Gamma = \partial \Omega,
\end{equation}
where $ \vec{\tau}, \vec{n}$ denote the tangential and normal vectors on the boundary $ \Gamma $. \eqref{ICNS} and \eqref{SBC} are given with the initial data 
\begin{equation}\label{Ini} \left. ( \rho, u)(x,0) \right. = (\rho_0, u_0)(x), ~~~~~ x \in \Omega. \end{equation}

Moreover, we take $ \Omega = B_1 $ being the ball centred at the origin with radius $ 1 $. Indeed, the geometry of $ \Omega $ would be an important factor in our problem. On one hand, as it is illustrated in the study of the stationary problem \cite{ReinhardFarwig1989}, the shape of the boundary $ \Gamma $ is an important factor in determining the steady states of \eqref{ICNS} with the boundary condition \eqref{SBC}. On the other hand, as pointed out in \cite{Desvillettes2002}, the Korn's inequality
\begin{equation}\label{OriginalKorn's}
	C_\Omega \left\Arrowvert \nabla V + \nabla V^\top \right\Arrowvert_{L^2} \geq \left\Arrowvert \nabla V \right\Arrowvert_{L^2},
\end{equation}
for a vector field $ V $ with the tangency boundary condition ($ u \cdot \vec{n} = 0 $ on $ \Gamma = \partial \Omega $) only holds when $ \Omega $ is a non-axisymmetric domain. Therefore the Korn's inequality \eqref{OriginalKorn's}, which gives the dissipation and plays an important role in the global analysis of \eqref{ICNS},  no longer applies directly in our setting. We shall resolve such issue in this work (see \eqref{korn02}). 

As it will be pointed out, \eqref{ICNS} with \eqref{SBC} admits a class of steady states with non-trivial velocity, which satisfy the equations \eqref{SNS} with $ \bar \rho \neq \text{constant} $. In short, instead of external forces, such non-trivial steady states are consequences of the self-rotation and the geometry of occupied domain. The main goal in this article is to investigate the stability of these steady states. 

It should be emphasised that, most of the available stability theories are subject to the no-slip boundary condition ($ \left. u \right|_{\partial \Omega} = 0$). However, the Navier-slip boundary condition is actually more appropriate when studying many phenomena such as hurricanes and tornadoes (see, \cite{Antontsev2007}). The general Navier-slip boundary condition is of the form
\begin{equation}\label{generalNSB}
	u \cdot \vec{n} = 0, ~ \vec{\tau}\cdot \mathbb S(u) \vec{n} = B(u\cdot \vec{\tau}), ~~~ \text{on} ~ \Gamma = \partial \Omega. 
\end{equation}
Such kinds of boundary conditions were first introduced by Navier \cite{Navier1827}. 
The first work on the mathematical rigorous analysis was due to Solonnikov and \v{S}\v{c}adilov \cite{Solonnikov1973} on the linearized stationary equations. A simplified form of \eqref{generalNSB} is by taking $ B(u\cdot \vec{\tau}) = - \kappa^{-1} u \cdot \vec{\tau} $ with $ \kappa $ being a constant. The most studied case in the literatures is when $ \kappa > 0 $, corresponding to slip with friction. In this situation, $\kappa$ is called the slip length. The case $ \kappa < 0 $ corresponds to the case in which the boundary wall accelerates the fluids (see \cite{Ding2016}). Our focus is on the case when the boundary does not provide friction or acceleration to the fluids, i.e. $ B(u\cdot\vec{\tau})=0 $. Such boundary is called the complete slip boundary (see \cite{ChenGQ2010}). 

For a homogeneous incompressible flow ($\dv u = 0, \rho = \text{constant} $), Chen and Qian in \cite{ChenGQ2010} demonstrated that when $ \kappa \rightarrow 0^+ $, the weak solutions to the incompressible Navier-Stokes equations converge to the solution for the problem subject to the no-slip boundary condition ($ \left. u \right|_{\partial \Omega} = 0 $) for almost all time. Also, as $ \kappa \rightarrow +\infty $, the solutions converge to a solution for the problem subject to the complete slip boundary condition. Recently, Ding, Li, Xin \cite{Ding2016} show that some instability may occur when $ \kappa < 0 $ in two spatial dimensional setting. The instability is in the following sense. $ \exists \epsilon > 0 $ such that $ \forall \delta > 0 $, there exists an initial perturbation of the steady state ($u_s=0$) that will grow larger than the size $ \epsilon $ in a suitable function space even thought such perturbation is smaller than the size $ \delta $ initially. Indeed, a critical value of the viscosity coefficient depending on $ \kappa $ serves as the threshold of stability and instability, and the instability would occur if and only if the viscosity coefficient is smaller than the critical value. The problem with complete slip boundary in an axisymmetric domain was studied by Watanabe in \cite{Watanabe2003}. It is shown that the global weak solution would converge to the projection of initial velocity on the rigid body motion \eqref{def:rigidbodymotion} in $ L^2 $ sense. 

For a compressible flow, Hoff \cite{Hoff2012} proved the local-in-time existence of smooth solutions to the Navier-Stokes equations with the boundary condition \eqref{generalNSB}. In \cite{Huang2006}, Huang, Li, Xin established the global dynamic property to the Stokes approximation of Navier-Stokes equations with the no-stick boundary condition in two dimensional setting. It is shown the solution $ (\rho, u) $ would converge the the equilibria $ (\rho_s, 0) $ in $ L^\alpha \times W^{1,\beta} $ space as time grows up. Recently, H. Li and X. Zhang in \cite{Li2016a} establish the nonlinear stability of Couette flows with the moving condition on the top and the Navier-slip boundary condition \eqref{generalNSB} in which $ B(u \cdot \vec{\tau}) = - \kappa^{-1} u\cdot \tau, \kappa > 0 $ on the bottom. The asymptotic stability of the trivial steady state $ (\rho_s = \text{constant}, u_s = 0 )$ to the problem with frictional boundary ($ \kappa > 0 $) was also demonstrated by Zaj\c{a}czkowski in \cite{Zajaczkowski1998}. 

As for the Cauchy problems and initial boundary valued problems with no-slip boundary for the compressible Navier-Stokes equations, there are rich literatures available . We shall only mention a few. In the absence of vacuum ($ \rho \geq \underline \rho > 0 $), the local and global well-posedness of classical solutions have been widely discussed. The uniqueness of both viscous and inviscid compressible flows was studied by Serrin \cite{Serrin1959}. Itaya \cite{Itaya1971} and Tani \cite{Tani1977} showed the local-in-time existence of classical solutions. Matsumura and Nishida \cite{Matsumura1980,Matsumura1983} first established the global well-posedness of classical solutions with a small pertuerbation of a uniform non-vacuum state. In the present of vacuum, as pointed out by Xin in \cite{zpxin1998} (also Xin, Yan in \cite{XinYan2013}), some singular behaviour may occur. With small initial energy, Huang, Li, Xin \cite{HuangLiXin2012} constructed the global smooth solutions for isentropic compressible viscous flows in $ \mathbb R^3 $. For more stability and instability problems, see \cite{Guo2011b,Guo2013a,Wang2014,Wang2016}. 

To study \eqref{ICNS} with \eqref{SBC}, we start with the stability theory as the first step. In this work, we will show the existence of rigid motions as steady solutions which rotate with uniform angular velocity. Unfortunately, such rotating profile may contain vacuum when the angular velocity is large. Indeed, the vacuum would appear around the symmetric axis of the domain. Thus it is supposed to be a vacuum interface problem. Moreover, the density profile admits physical vacuum across the vacuum interface, which contains singularities and is not suited in our functional framework (see \cite{Liu1996}). We will study such vacuum interface problem in the future and focus on the non-vacuum problem here. To study the stability of the rigid motion, we would make use of the conservation of angular momentum (Lemma \ref{lm:ConserAnguMome}) to establish a Korn's type inequality \eqref{korn02}, which would play important roles in global analysis. This is inspired by the study of free boundary problems in \cite{Zajaczkowski1993}. Moreover, it is introduced a spherical frame which matches the geometry of $ \Omega $. With such structure, we would be able to define some differential operators \eqref{def:DiffFrame}, which separate the normal derivatives and tangential derivatives. The benefit of doing so is that it avoids applying the partition of unity of the domain and works without introducing local charts. Therefore, we can calculate in the entire domain at once. The differential operators in this work are natural to the spherical domain and it is convenient to separate $ \dv \mathbb S(u) $ on the right of the momentum equation \subeqref{ICNS}{2} into tangential and normal directions (see \eqref{NDEst002}). 

The rest of this work would be organised as follows. In the next section, we shall construct the steady state which we are going to study. In Section \ref{sec:preliminaries}, we present the main tools that will be used in this work, including the Korn's type inequality, the differential operators and the classical elliptic estimates on the Stokes' problem. The equations and the main theorem in the perturbation variables would be given in Section \ref{sec:perturbedform}. The main energy estimates are listed in Section \ref{sec:enerest}. To illustrate the program, the estimates on the temporal derivatives and lower order spatial derivatives will be recorded in Section \ref{sec:lowerorderest}. Then we move onto the higher order spatial derivatives and interior estimates. With such blocks in hands, we chain them together in Section \ref{sec:globalwellposed} and show the asymptotic stability. In particular, we demonstrate the nonlinearities can be indeed controlled by the energy in Section \ref{sec:nonlinearity}.

\subsection{Steady States and Main Results}

Here we search for non-vacuum steady states of \eqref{ICNS}. That is, to find $ (\bar \rho, \bar u ) $ with $ \bar \rho > 0 $ satisfying 
\begin{equation}\label{SNS} 
	\begin{cases}
		\dv (\bar \rho \bar u) = 0 & \text{in} ~ \Omega,\\
		\dv ( \bar \rho \bar u \otimes \bar u) + \nabla \bar P = \dv \mathbb{S}(\bar u) & \text{in} ~ \Omega,\\
		\bar u \cdot \vec{n}, ~ \vec{\tau} \cdot \mathbb{S}(\bar u) \vec{n} = 0 & \text{on} ~  \Gamma,
	\end{cases} 
\end{equation}
with $ \bar P = \bar \rho^\gamma, \mathbb{S}(\bar u) = \mu ( \nabla \bar u + \nabla \bar{u}^\top ) + \lambda \dv \bar u \mathbb{I}_3 $. After taking inner product of $ \eqref{SNS}_2 $ with $ \bar u $ and integration by parts, record the resulting equation,
\begin{equation}\label{SNS1}
	\begin{aligned}
		& - \dfrac{1}{2} \int_\Omega \dv ( \bar \rho \bar u) |\bar u|^2\,dx + \dfrac{1}{2} \int_\Gamma \rho |\bar u|^2 \bar u \cdot \vec n \,dS - \dfrac{\gamma}{\gamma-1} \int_\Omega \dv(\bar \rho \bar u)   \cdot \bar \rho^{\gamma-1}\,dx \\
		& ~~~ + \dfrac{\gamma}{\gamma-1} \int_\Gamma \bar\rho^\gamma \bar u \cdot \vec{n}\,dS = - \dfrac{\mu}{2} \int_\Omega \left| \nabla \bar u + \nabla {\bar u}^\top \right|^2\, dx 
		- \lambda \int_\Omega \left| \dv \bar u \right|^2\,dx\\
		& ~~~ + \int_\Gamma \bar u \cdot \mathbb{S}(\bar u) \vec{n}\,dS.
	\end{aligned}
\end{equation}
Since $ \bar u \cdot \vec{n} = 0 $ on $ \Gamma $, $ \bar u $ is in the tangential direction of $ \Gamma $. Therefore, $ \eqref{SNS}_3 $ implies $ \bar u \cdot \mathbb{S}(\bar u) \vec{n} = 0 $ on $ \Gamma $. Together with $ \eqref{SNS}_1, \eqref{SNS}_3 $, \eqref{SNS1} can be rewritten as 
\begin{equation*}
	\dfrac{\mu}{2}\int_\Omega \left| \nabla \bar u + \nabla {\bar u}^\top \right|^2\,dx + \lambda \int_\Omega \left|\dv \bar u \right|^2\,dx = 0,
\end{equation*}
and hence
\begin{equation}\label{vshingsymmatrix}
	\nabla \bar u + \nabla {\bar u}^\top = 0, ~ \text{in} ~ \Omega.
\end{equation}
It follows, via similar arguments as in \cite{XinYan2013} (also \cite{Cho2006})
\begin{equation*}
	\bar u = \bar A x + \bar b, ~~ x \in \Omega,
\end{equation*}
for some anti-symmetric matrix $ \bar A $ and constant vector $ \bar b $. Moreover, on $ \Gamma $, $ x $ is paralleled to $ \vec{n} $ as the consequence of the geometry of $\Omega $. Hence $ 0 = \bar u \cdot \vec{n} = \bar A x \cdot \vec{n} + \bar b \cdot \vec{n} = \bar b \cdot \vec{n} $ on $ \Gamma $. It follows $ \bar b = 0 $. Therefore, there is a vector $ \bar a $ such that 
\begin{equation}\label{def:rigidmotion}
	\bar u = x \times \bar a, ~~ \text{in} ~ \Omega.
\end{equation}
Without lost of generality, by assuming $ \bar a = \bar \omega e_3 $ for $ \bar \omega \geq 0 $, \eqref{SNS} is reduced to the following ordinary differential equation(ODE)
\begin{equation}\label{SNS2}
	\begin{aligned}
		& - \bar \omega^2 \bar \rho r + \dr {\bar \rho}^\gamma = 0,
	\end{aligned}
\end{equation}
with $ \bar \rho = \bar \rho (r) $ where $ r = \sqrt{(x_1)^2+(x_2)^2} $. Therefore,
\begin{equation}\label{SNS3}
	\begin{cases}
		& \bar \rho = \left( \bar\rho^{\gamma-1}(0) + \dfrac{\gamma-1}{2\gamma} \bar\omega^2 r^2 \right)^{1/(\gamma-1)},\\
		& \bar u = \bar\omega \left( x_2, - x_1, 0 \right)^\top,
	\end{cases}
\end{equation}
where $ \bar \rho(0) > 0, \bar \omega \geq 0 $. 

Furthermore, we shall denote the rigid body motions
\begin{equation}\label{def:rigidbodymotion}
\begin{aligned}
	& S = \left\lbrace V \in H^1(\Omega,\mathbb R^3); \nabla V + \nabla V^\top = 0, V \cdot\vec{n} = 0  ~ \text{on} ~ \partial\Omega \right\rbrace \\
	& ~~ = \left\lbrace V(x) = x \times \zeta; V \cdot \vec{n} = 0 ~ \text{on} ~ \partial\Omega, \zeta ~ \text{is a constant vector} \right\rbrace.
\end{aligned}
\end{equation}
Also, let $ P_S $ be the orthogonal projection from $ H^1(\Omega,\mathbb R^3) $ onto $ S $. Notice, $ \forall V \in S $, since the domain is symmetric, it holds $ \inti{V} = 0 $. In particular, $ \forall V \in H^1(\Omega,\mathbb R^3) $,
\begin{equation*}
	\inti{P_S V} = 0.
\end{equation*}

\begin{rmk}
	We use the assumption $ \mu, \lambda > 0 $ in the deviation of \eqref{vshingsymmatrix}. However, the Lam\'e relation ($ \mu > 0, 3\lambda + 2\mu > 0 $) is sufficient to show this. 
\end{rmk}
\begin{rmk}
It is assumed $ \inf_{x\in \Omega} \bar \rho (x) > 0 $. From \eqref{SNS2}, the minimum of $\bar \rho $  is achieved on the axis $ \left\lbrace x_1 = x_2 = 0 \right\rbrace $. For a fixed $ \bar\omega > 0 $ and the fixed domain $ \Omega$, it follows that the total mass of the fluid inside $ \Omega $ has a lower bound $ {\underline M}_{\bar\omega} > 0 $ in order to avoid vacuum. In reality, when the total mass drops below $ {\underline M}_{\bar\omega} $, there would be a vacuum area inside $ \Omega $ and a vacuum interface. Indeed, from \eqref{SNS2}, the density profile across the vacuum interface would admit physical vacuum \cite{Liu1996}. Similar phenomena would occur when the total mass is fixed and the angular velocity increases. We leave such vacuum problems as future works.  
\end{rmk}

In this work, it is devoted to study the stability of the steady state \eqref{SNS3}. In fact, we have the following informal statement of our theorem.
\begin{thm}[Informal Statement]
	Provided the angular velocity $ |\bar \omega| $ is small though, the steady state $ ( \bar\rho, \bar u) $ given by \eqref{SNS3} to the compressible Navier-Stokes equations \eqref{ICNS} with the complete slip boundary condition \eqref{SBC} is asymptotically stable in the following sense. For an initial perturbation with the size less than $ \epsilon_1 $ for some $ \epsilon_1>0 $ in some appropriate function space (defined in \eqref{total-Energy} and \eqref{total-energy}), there is a globally defined classical solution to \eqref{ICNS} and the solution converges to the steady state as time grows up exponentially. The perturbation is taken such that it preserves the angular momentum (see \eqref{Ini1}). 
\end{thm}
We shall state the theorem in perturbation variables later in Section \ref{sec:perturbedform}. 


\section{Preliminaries}\label{sec:preliminaries}

\subsection{Notations}

Through out this work, conventionally, for any quantities $ A, B $,
\begin{equation*}
	\text{by} ~~ A \lesssim B, ~~ \text{it is to say} ~~ A \leq C B,
\end{equation*}
where the constant $ C > 0 $ 
may depend on $ \sup_{x\in\Omega} \bar \rho, \inf_{x\in\Omega} \bar\rho , \Omega, \mu, \lambda $ but is independent of $ \bar\omega, \rho, u $.
Similarly, 
\begin{equation*}
	A \simeq B
~~ \text{is equal to say} ~~
	\dfrac{1}{C} A \leq B \leq C A.
\end{equation*}
For a constant $ 0 < \omega < 1 $, the corresponding value $ C_\omega $ would denote a positive value satisfying $ 1 \leq C_\omega \leq 1/\omega $. In the meantime, for a vector $ \vec{w} = (w^1,w^2,w^3)^\top $, the associated differential operator is defined as,  
\begin{equation*}
	\nabla_{\vec{w}} = \vec{w}\cdot \nabla = w^1 \partial_{1} + w^2 \partial_{2} + w^3 \partial_{3},
\end{equation*}
where $ \partial_{i} = \partial_{x_i} $ represents the spatial derivative for $ i = 1,2,3 $. 
Also, the commutator operator is defined by
\begin{equation*}
	\left\lbrack A, B \right\rbrack = A B - B A,
\end{equation*}
where $ A,B $ may stand for functions or differential operators. Notice $ \left\lbrack\cdot, \cdot \right\rbrack $ is bilinear. The Sobolev norm in $\Omega $ and on the boundary $ \Gamma = \partial \Omega $ is denoted as
\begin{align*}
	& \norm{f}{\Lnorm}^2 = \inti{f^2},
	& \norm{f}{\Hnorm{k}} = \sum_{i=0}^{k} \norm{\nabla^i f}{\Lnorm},\\
	& \norm{f}{\bLnorm}^2 = \intb{f^2},
	& \norm{f}{\bHnorm{k}} = \sum_{i=0}^{k} \norm{\bar\nabla^i f}{\bLnorm},
\end{align*}
where $ \nabla, \bar\nabla $ denote the differential operators in $ \Omega $ and on the boundary $ \Gamma $ respectively. 

\subsection{Korn's Inequality}

The following form of Korn's inequality is from \cite[Lemma 5.1]{Zajaczkowski1993}.
\begin{lm}[Korn's Inequality] For $ V = \left( V^1, V^2, V^3 \right)^\top $ defined in $ \Omega $, it holds
	\begin{equation}\label{korn01}
	\int_{\Omega} \left| \nabla V \right|^2\,dx \lesssim \int_{\Omega} \left| \nabla V + \nabla V^\top \right|^2 \,dx + \int_{\Omega} \left| V\right|^2\,dx \end{equation}
	provided the right hand side is finite. 
\end{lm}
For the sake of completeness, we show the proof here.
\begin{pf} Denote the energy of the symmetric part of any vector $ U = (U^1, U^2, U^3)^\top $ as ,
\begin{equation*}
	E_\Omega(U) = \int_\Omega \left| \nabla U + \nabla U^\top \right|^2 \,dx.
\end{equation*}
	Introduce a decomposition of $ V $,
	\begin{equation}\label{korn1}
		V = \sum_{i=1}^{3} b^i \varphi_i(x) + \bar{V}, 
	\end{equation}
	where
	\begin{equation}\label{korn2}
		\varphi_i(x) = (x - \bar x) \times e_i,
	\end{equation}
	with $x = (x_1,x_2,x_3)^\top, \bar x = \frac{1}{|\Omega|}(\int_\Omega x_1\,dx,\int_\Omega x_2\,dx,\int_\Omega x_3,dx)^\top $ and $ e_i = (\delta_{i1}, \delta_{i2}, \delta_{i3} )^\top, i = 1,2,3$. In particular, if $ \Omega $ is a ball, $ \bar x = 0 $. 
	Define $ b=(b^1,b^2,b^3)^\top$ by
	\begin{equation}\label{korn9}
		b = - \dfrac{1}{2|\Omega|} \int_{\Omega} \curl V \,dx.
	\end{equation}
	Then direct calculation yields
	\begin{equation}\label{korn3}
		\int_\Omega V\,dx = \int_\Omega \bar V \,dx,~ \int_\Omega \curl \bar V\,dx = 0,~ E_\Omega(V) = E_\Omega(\bar V).
	\end{equation} 
	It follows from the original Korn's inequality (see \cite{Friedrichs1947,horgan1983}),
	\begin{equation}\label{korn5}
		\int_\Omega \left| \nabla \bar V \right|^2\,dx \lesssim E_\Omega(\bar V) = E_\Omega(V).
	\end{equation}
	From the decomposition \eqref{korn1}, 
	\begin{equation}\label{korn4}
		\int_\Omega \left| \nabla V \right|^2 \,dx \lesssim |b|^2 + \int_\Omega \left| \nabla \bar V \right|^2\,dx \lesssim |b|^2 + E_\Omega(V).
	\end{equation}
	It remains to estimate $ |b|^2 $. From \eqref{korn1}, we consider the equations
	\begin{equation*}
		\sum_{i=1}^{3} b^i \int_{\Omega} \varphi_i(x) \cdot \varphi_j(x) \,dx = \int_{\Omega} (V - \bar V) \cdot \varphi_j(x)\,dx, ~~~~ j = 1,2,3.
	\end{equation*}
	The non-degeneracy of the coefficient matrix $\left\lbrace \int_G \varphi_i(x) \cdot \varphi_j(x) \,dx \right\rbrace $ implies
	\begin{equation}\label{korn6}
		|b|^2 \lesssim \int_\Omega |V|^2 \,dx + \int_\Omega | \bar V |^2\,dx.
	\end{equation}
	Therefore, applying the Poincar\'{e} inequality together with \eqref{korn3} and \eqref{korn5}, it holds,
	\begin{equation}\label{korn7}
	\begin{aligned}
		& \int_\Omega \left|\bar V \right|^2\,dx \lesssim \int_{\Omega} \left| \bar V - \dfrac{1}{|\Omega|} \int_\Omega \bar V \,dx \right|^2 \,dx + \int_\Omega \left| \dfrac{1}{|\Omega|}\int_\Omega \bar V\,dx\right|^2\,dx \\
		& ~~ \lesssim \int_\Omega \left| \nabla \bar V \right|^2\,dx + \left(\int_\Omega V \,dx\right)^2  \lesssim  E_\Omega(V) + \int_\Omega |V|^2\,dx.
	\end{aligned}
	\end{equation}
	\eqref{korn01} follows from \eqref{korn4}, \eqref{korn6} and \eqref{korn7}.
\end{pf}

In addition, the following lemma consists of the $ L^2 $ estimate of the orthogonal complement of $ V $ with respect to the rigid motions \eqref{def:rigidbodymotion}.
\begin{lm}[Poincar\'e-Morrey Inequality]\label{lm:Poincare-morreyIneq}
	For $ V \in H^1(\Omega, \mathbb R^3) $ with $ V \cdot \vec{n} = 0 $ on $ \partial\Omega $, we have
	\begin{equation}\label{korn8}
		\norm{V - P_S V}{\Lnorm}^2 \lesssim \inti{\abs{\nabla V + \nabla V^\top}{2}}.
	\end{equation}
\end{lm}
The proof contains a compactness argument. We refer the proof to \cite[Lemma 4.2]{Watanabe2003}. See also \cite{ReinhardFarwig1989}. 
As a corollary, one can derive the following Poincar\'e inequality.
\begin{cor}
	The same vector $ V $ as in Lemma \ref{lm:Poincare-morreyIneq} would satisfy the following
	\begin{equation}\label{poincareball}
		\norm{V}{\Lnorm} \lesssim \norm{\nabla V}{\Lnorm}.
	\end{equation}
\end{cor}

\begin{pf}
This form of Poincar\'e inequalities can be found in \cite{Bishop1988}. However, a new proof is provided here. 
One can rewrite
\begin{equation*}
	V = V - \dfrac{1}{|\Omega|} \inti{V} + \dfrac{1}{|\Omega|}\inti{(V - P_S V)},
\end{equation*}
where it has been used the fact $ \inti{P_S V} = 0 $. Therefore, the Poincar\'e inequality and \eqref{korn8} then yield
\begin{equation*}
	\begin{aligned}
		& \norm{V}{\Lnorm}^2 \lesssim \norm{\nabla V}{\Lnorm}^2 + \norm{\nabla V + \nabla V^\top}{\Lnorm}^2 \lesssim \norm{\nabla V}{\Lnorm}^2.
	\end{aligned}
\end{equation*} 
\end{pf}


An interesting property of \eqref{ICNS} comes from the complete slip boundary \eqref{SBC}. More precisely, in the symmetric domain $ \Omega = B_1 $, the flow admits conservation of angular momentum. 

\begin{lm}[Conservation of Angular Momentum]\label{lm:ConserAnguMome} With the complete slip boundary condition \eqref{SBC}, for any smooth solution to \eqref{ICNS} with the initial data \eqref{Ini}, it holds
	\begin{equation}
		\int_\Omega \rho u \cdot \varphi_i \,dx = \int_\Omega \rho_0 u_0 \cdot \varphi_i \,dx, ~~ i = 1,2,3,
	\end{equation}
	where $ \varphi_i $ are defined in \eqref{korn2} with $ \bar x = 0 $.
\end{lm}

\begin{pf}
	Take inner product of $ \eqref{ICNS}_2 $ with $ \varphi_i $ ($ i =1,2,3 $), and record the resulting after integration by parts,
	\begin{equation}\label{cnsv1}
		\begin{aligned}
			& \dfrac{d}{dt} \int_\Omega \rho u \cdot \varphi_i \,dx - \int_\Omega  \rho u \cdot (u\cdot \nabla) \varphi_i\,dx + \int_{\Gamma} ( \rho u\cdot \varphi_i) \cdot (u \cdot \vec{n}) \,dS \\
			&  - \int_\Omega P \cdot \dv \varphi_i \,dx + \int_\Gamma P \varphi_i \cdot \vec{n} \,dS = - \int_\Omega  \mathbb S (u) : \nabla \varphi_i \, dx + \int_\Gamma \varphi_i \cdot ( \mathbb S (u) \vec n)\, dS.
		\end{aligned}
	\end{equation}
	Notice $ \varphi_i \cdot \vec n =  ( x \times e_i ) \cdot \frac{x}{|x|} = 0, u \cdot \vec n = 0 $ on $ \Gamma $, which means both $ \varphi_i $ and $ u $ are in the tangential direction on $ \Gamma $. Therefore all the boundary integrations above vanish. In the meantime, $ \nabla \varphi_i $ is an anti-symmetric matrix, while $ \mathbb S (u) $ is symmetric. Hence $ \dv \varphi_i = \tr (\nabla \varphi_i) = 0, u \cdot (u \cdot \nabla)\varphi_i = u \cdot (\nabla \varphi) u = 0, \mathbb S (u) : \nabla \varphi_i = 0 $. Therefore, from \eqref{cnsv1}
	\begin{equation}
		\dfrac{d}{dt} \int_\Omega \rho u \cdot \varphi_i \,dx = 0.
	\end{equation}
\end{pf}


Now we would able to derive the Korn's type inequality which is our first important block in this work. 
\begin{lm}[Korn's Type Inequality]\label{lm:korn'stypeinequality} For any smooth solution to \eqref{ICNS} with
	\begin{equation}
		\int_\Omega \rho u \cdot \varphi_i \,dx = \int_\Omega \bar \rho \bar u \cdot \varphi_i \,dx, ~~ i =1,2,3, ~~ \text{and} ~~ \bar \rho \geq \inf_{x\in \Omega} \bar \rho > 0 ~ \text{in} ~ \Omega,
	\end{equation}
	it holds
	\begin{equation}\label{korn02}
	\begin{aligned}
		& \int_\Omega \left|\nabla ( u - \bar u ) \right|^2 \,dx \lesssim \int_\Omega \left| \nabla u + \nabla u^\top \right|^2 \,dx + \int |\bar u|^2 |\rho - \bar\rho|^2\,dx \\
		& ~~~~~~ + \int_\Omega | \rho  - \bar \rho |^2| u - \bar u|^2 \,dx  + \left( \int_\Omega u \,dx \right)^2.
	\end{aligned}
	\end{equation}
	Moreover, the last term $ \inti{u} $ on the right can be dropped in \eqref{korn02}.
\end{lm}

\begin{pf}
	Decompose $ u - \bar u = \sum_{i=1}^{3} b^i \varphi_i + v $, where $ b = (b^1, b^2, b^3)^\top $ are defined similarly as in \eqref{korn9} with $ V = u - \bar u $. Similarly, it holds 
	\begin{equation}\label{cnsv2}
		\begin{gathered}
			\int_\Omega v\,dx = \int_\Omega u\,dx,~\int_\Omega \left| \nabla v \right|^2 \,dx \lesssim \int_\Omega \left | \nabla u + \nabla u^\top \right|^2 \,dx,\\
			\int_\Omega \left| \nabla (u- \bar u) \right|^2 \,dx \lesssim |b|^2 + \int_\Omega \left | \nabla u + \nabla u^\top \right|^2 \,dx.
		\end{gathered}
	\end{equation}
	Notice
	\begin{equation*}
	\begin{aligned}
		& 0 = \int_\Omega \rho u \cdot \varphi_i \,dx - \int_\Omega \bar \rho \bar u \cdot \varphi_i\,dx = \int_\Omega (\rho - \bar\rho ) ( u - \bar u ) \cdot \varphi_i \,dx \\
		& ~~~~~ + \int_\Omega (\rho - \bar \rho ) \bar u \cdot \varphi_i\,dx + \int_\Omega \bar\rho ( u - \bar u ) \cdot \varphi_i\,dx.
	\end{aligned}
	\end{equation*}
	Therefore, after plugging in the decomposition of $ u - \bar u $, we have the following system of equations ($ i = 1,2,3 $ ),
	\begin{equation*}
		\begin{aligned}
			& \sum_{j=1}^{3} b^j \int_\Omega \bar\rho \varphi_j \cdot \varphi_i \,dx = - \int_\Omega \bar\rho v\cdot \varphi_i\,dx \\
			& ~~~~~ - \int_\Omega (\rho - \bar\rho ) ( u - \bar u ) \cdot \varphi_i \,dx - \int_\Omega (\rho - \bar \rho ) \bar u \cdot \varphi_i\,dx.
		\end{aligned}
	\end{equation*}
	Notice $ \left( \int_\Omega \bar\rho \varphi_i \cdot \varphi_j \,dx \right)$ is a non-degenerate matrix. It follows
	\begin{equation*}
	\begin{aligned}
		& | b |^2 \lesssim \int_\Omega |\bar\rho|^2 |v|^2\,dx + \int_\Omega |\rho - \bar\rho|^2| u - \bar u |^2\,dx + \int_\Omega |\bar u|^2 |\rho - \bar\rho|^2\,dx.
	\end{aligned}
	\end{equation*}
	By using the Poincar\'e inequality and \eqref{cnsv2},
	\begin{equation*}
		\begin{aligned} 
		& \int_\Omega |\bar \rho |^2 |v|^2\,dx \lesssim \int_\Omega \left| v - \frac{1}{|\Omega|}\int_\Omega v \,dx  \right|^2 \,dx + \int_\Omega \left| \frac{1}{|\Omega|} \int_\Omega v\,dx \right|^2\,dx\\
		 & ~~~~ \lesssim \int_\Omega \left| \nabla v \right|^2 \,dx + \left( \int_\Omega v\,dx \right)^2 \lesssim \int_\Omega \left| \nabla u + \nabla u^\top \right|^2\,dx + \left( \int_\Omega u \,dx \right)^2.
		 \end{aligned}
	\end{equation*}
	\eqref{korn02} follows by chaining the above inequalities. In addition, the identity
	\begin{equation*}
		\inti{u} = \inti{(u - P_S u)}
	\end{equation*}
	yields
	\begin{equation*}
		\left( \inti{u} \right)^2 \lesssim \norm{u - P_S u}{\Lnorm}^2 \lesssim \inti{\abs{\nabla u + \nabla u^\top}{2}}
	\end{equation*}
	by applying \eqref{korn8}.
\end{pf}

Similarly, we have the following 
\begin{lm} Under the same assumptions as in Lemma \ref{lm:korn'stypeinequality},
	\begin{equation}\label{korn03}
		\begin{aligned}
			& \int_\Omega \left| \nabla \dt u \right|^2\,dx \lesssim \int_\Omega \left| \nabla \dt u + \nabla \dt u^\top \right|^2\,dx  + \int_\Omega |\dt \rho |^2 |u-\bar u |^2\,dx \\& ~~~~ +  \int_\Omega |\rho - \bar \rho |^2|\dt u|^2\,dx +  \int_\Omega |\bar u|^2 |\dt \rho|^2 \,dx +  \left( \int_\Omega \dt u \, dt \right)^2.
		\end{aligned}
	\end{equation}
	Generally, for any integer $ k \geq 1 $,
	\begin{equation}\label{korn04}
	\begin{aligned}
		& \int_\Omega \left| \nabla \dt^k u \right|^2 \,dx \lesssim \int_\Omega \left| \nabla \dt^k u + \nabla \dt^k u^T \right|^2 \,dx + \int_\Omega |\bar u|^2 |\dt^k \rho |^2 \,dx \\
		& ~~~~ + \sum_{j=0}^{k} \int_\Omega |\dt^j (\rho-\bar \rho)|^2 |\dt^{k-j} (u - \bar u )|^2 \,dx + \left(\int_\Omega \dt^k u \,dx \right)^2.
	\end{aligned}	
	\end{equation}
	In particular,the last terms on the right of \eqref{korn03} and \eqref{korn04} can be dropped.
\end{lm}
The proofs of \eqref{korn03}, \eqref{korn04} are similar to that of \eqref{korn02}.

\subsection{Spherical Differential Frame}\label{sec:sphdiffframes}


Here, we introduce the decomposition of differential operators $ \nabla = ( \partial_{x_1}, \partial_{x_2}, \partial_{x_3} ) $ near the boundary $ \Gamma $ into two kinds of operators, corresponding to tangential and normal derivatives respectively.

To begin with, define the following cut-off function,
\begin{equation*}
		\psi : C_c^\infty(B_{3/4}) \longmapsto [0,1]
\end{equation*}
satisfying $ \left| \nabla \psi(x) \right| \leq 8 $ and
\begin{equation}\label{def:cut-off}
\psi(x) 
	\begin{cases}
		= 1 & \abs{x}{} \leq 1/2, \\
		=  0 & \abs{x}{} \geq 3/4.
	\end{cases}
\end{equation}
Our differential operators are defined as
\begin{equation}\label{def:DiffFrame}
	\begin{aligned}
		& \nabla_{T} = \varphi_i(x)\cdot \nabla, ~ \text{for} ~ i =1,2,3,\\
		& \nabla_{N} = N(x) \cdot \nabla,
	\end{aligned}
\end{equation}
where $\varphi_i $ is defined in \eqref{korn2} with $ \bar x = 0 $, $ N(x) = (x_1, x_2, x_3)^\top $. 
In the following, we study the properties of these differential operators. 
%
First, the commutator of two differential operators is also a differential operator. More precisely,
\begin{lm}[Commutator] For any function $ f: \mathbb{R}^3 \longmapsto \mathbb{R} $ and any vector fields $ \alpha = ( \alpha^1, \alpha^2, \alpha^3)^\top, \beta = (\beta^1,\beta^2,\beta^3)^\top $,
	\begin{equation}\label{commt}
		\commut{\nabla_\alpha}{\nabla_\beta}f(x) = \phi_{\alpha,\beta} \cdot \nabla f(x),
	\end{equation}
	for some $ \phi_{\alpha,\beta} $ satisfying $$ \norm{\phi_{\alpha,\beta}}{\supnorm} \lesssim \norm{\beta}{\supnorm}\norm{\nabla \alpha}{\supnorm} + \norm{\alpha}{\supnorm}\norm{\nabla\beta}{\supnorm} . $$ Moreover, $ \phi_{\alpha,\beta} \in C^\infty $ provided $ \alpha,\beta \in C^\infty $. 
\end{lm}
	
\begin{pf}
	For $ i, j \in \left\lbrace 1,2,3 \right\rbrace $, direct calculation yields
	\begin{equation*}
		\begin{aligned}
			& \commut{\alpha^i \partial_i}{ \beta^j \partial_j } f = \alpha^i \partial_i (\beta^j \partial_j f) - \beta^j \partial_j (\alpha^i \partial_i f) = \alpha^i\partial_i\beta^j \cdot \partial_j f  - \beta^j \partial_j \alpha^i \cdot \partial_i f .
		\end{aligned}
	\end{equation*}
	This finishes the proof.
\end{pf}



In the meantime, the following lemma shows that $ \nabla
 $ can be indeed decomposed into $ \nabla_T $ and $ \nabla_N $ in the boundary subdomain, in the sense that the estimates of $ L^2 $ norms of $ \nabla_T f, \dN f $ would be sufficient to obtain the corresponding estimate of $ \nabla f $. 

\begin{lm} For any smooth function $ f: \mathbb{R}^3 \longmapsto \mathbb{R} $,
	\begin{align}
	& \norm{\nabla f}{L^2(\Omega\backslash B_{1/2})} \lesssim \norm{\nabla_T f}{L^2(\Omega\backslash B_{1/2})} + \norm{\nabla_N f}{L^2(\Omega\backslash B_{1/2})}  , \label{DFrame01} \\
	&	\norm{\nabla f}{L^2(\Omega)} \lesssim \norm{\nabla (\psi f)}{L^2(\Omega)} + \norm{\nabla_{T} f}{L^2(\Omega)} + \norm{\nabla_{N} f}{L^2(\Omega)} \lesssim \norm{f}{\Hnorm{1}}, \label{DFrame02}
	\end{align}
	where $ \norm{\nabla(\psi f)}{L^2(\Omega)} $ can be replaced by $ \norm{\nabla f}{L^2(B_{1/2})} $.
	Here and in the following, $ \norm{ \dT f}{(\cdot)} $ stands for the sum of $$ \norm{\nabla_{\varphi_1} f}{(\cdot)}, \norm{\nabla_{\varphi_2} f}{(\cdot)}, \norm{\nabla_{\varphi_3} f}{(\cdot)}. $$ We adopt this convention through the rest of this work.
\end{lm}

\begin{pf}
	First we separate $ \Omega $ into the interior and boundary subdomains, $ \Omega = B_{1/2} \cup \left( \Omega \backslash B_{1/2}\right) $. Then
	\begin{equation*}
		\norm{\nabla f}{L^2(\Omega)}^2 = \norm{\nabla f}{L^2{(B_{1/2})}}^2 + \norm{\nabla f}{L^2(\Omega\backslash B_{1/2})}^2.
	\end{equation*}
	Then \eqref{DFrame02} follows from \eqref{DFrame01} and the property of $\psi $ in \eqref{def:cut-off}. Thus it remains to show \eqref{DFrame01}. It is sufficient to show that in $ \Omega \backslash B_{1/2} $ , the rank of $ \left\lbrace \varphi_1, \varphi_2, \varphi_3, N \right\rbrace $ is equal to three. Notice, in the boundary subdomain, $ x_1^2 + x_2^2 + x_3^2 \geq 1/4 $ and thus at least one of $ |x_1|, |x_2|, |x_3| $ is no less than 1/4. 
	Without loss of generality, we only consider the case when $ x_1 \geq 1/4 $. Direct calculation yields
	\begin{equation*}
	\det ( \varphi_2, \varphi_3, N ) = 
	\left|\begin{array}{ccc}
			-x_3 & x_2 & x_1 \\
			0 & -x_1 & x_2 \\
			x_1 & 0 & x_3
		\end{array}\right|
		 = x_1 \left( x_1^2 + x_2^2 + x_3^2 \right) \geq (1/4)^2 > 0.
	\end{equation*}
	This finishes the proof.
\end{pf}

For higher order derivatives, it also admits the following lemma. 

\begin{lm}
	\begin{align}
		& \norm{\nabla^2 f}{\Lnorm} \lesssim \norm{\nabla^2 (\psi f)}{\Lnorm} + \norm{\nabla_T^2 f}{\Lnorm} + \norm{\nabla_T \nabla_N f}{\Lnorm}  {\nonumber} \\ & ~~~~~~~~~~~~~~~~~  + \norm{\nabla_N^2 f}{\Lnorm}  + \norm{\nabla f}{\Lnorm},    \label{DFrame03} \\
		& \norm{\nabla^k f}{\Lnorm} \lesssim \norm{\nabla^k (\psi f)}{\Lnorm} + \sum_{i+j = k} \norm{\nabla_T^i\nabla_N^j f}{\Lnorm} {\nonumber} \\ & ~~~~~~~~~~~~~~~~~ + \sum_{i \leq k-1 } \norm{\nabla^i f}{\Lnorm} , \label{DFrame04}
	\end{align}
	where $ \norm{\nabla^2(\psi f)}{L^2(\Omega)}, \norm{\nabla^k(\psi f)}{L^2(\Omega)} $ can be replaced by $ \norm{\nabla^2 f}{L^2(B_{1/2})}$ and $ \norm{\nabla^k f}{L^2(B_{1/2})} $ respectively. 
	In these inequalities, $ \dN $ and $ \dT $ can be interchanged. 
\end{lm}

\begin{pf}
	We show \eqref{DFrame03} only. \eqref{DFrame04} can be proved via a similar argument. 
	Notice,
	\begin{equation*}
		\norm{\nabla^2 f}{L^2} = \norm{\nabla^2 f}{L^2(B_{1/2})} + \norm{\nabla^2 f}{L^2(\Omega\backslash B_{1/2})}.
	\end{equation*}
	Then, apply \eqref{DFrame01} repeatly,
	\begin{equation*}
		\begin{aligned}
			& \norm{\nabla^2 f}{L^2(\Omega\backslash B_{1/2})} \lesssim \norm{\nabla_T \nabla f}{L^2(\Omega\backslash B_{1/2})} + \norm{\nabla_N \nabla f}{L^2(\Omega\backslash B_{1/2})} \\
			& ~~~~~~ \lesssim \norm{\nabla \nabla_T f}{L^2(\Omega\backslash B_{1/2})} + \norm{\nabla \nabla_N f}{L^2(\Omega\backslash B_{1/2})} +  \norm{\nabla f}{L^2(\Omega\backslash B_{1/2})} \\
			& ~~~~~~ \lesssim \norm{\nabla_T \nabla_T f}{L^2(\Omega\backslash B_{1/2})} + \norm{\nabla_N \nabla_T f}{L^2(\Omega\backslash B_{1/2})} \\ & ~~~~~~~~~~ + \norm{\nabla_T\nabla_N f}{L^2(\Omega\backslash B_{1/2})} + \norm{\nabla_N\nabla_N f}{L^2(\Omega\backslash B_{1/2})} + \norm{\nabla f}{L^2(\Omega\backslash B_{1/2})}
		\end{aligned}
	\end{equation*}
	where in the second inequality \eqref{commt} has been applied. Again, after applying \eqref{commt} again, it holds
	\begin{equation*}
		\norm{\nabla_N\nabla_T f}{L^2(\Omega\backslash B_{1/2})} \lesssim \norm{\nabla_T\nabla_N f}{L^2(\Omega\backslash B_{1/2})} + \norm{\nabla f}{L^2(\Omega\backslash B_{1/2})}.
	\end{equation*}
	Thus this finishes the proof of \eqref{DFrame03}.
\end{pf}


The next lemma concerns the calculus on the boundary $ \Gamma $. In fact, $ \left.\dT\right|_{\Gamma} $ is a differential operator on the boundary. 

\begin{lm}\label{FTCB}
	For any smooth function $ g: \Gamma \longmapsto \mathbb{R} $, it holds,
	\begin{equation}\label{FTCB1} \int_\Gamma \dT g \,dS = 0. \end{equation}
	As a corollary, for $ g_1, g_2: \Gamma \longmapsto \mathbb{R} $,
	\begin{equation}\label{FTCB2}
		\int_\Gamma \dT g_1 \cdot g_2 \,dS = - \int_\Gamma g_1 \cdot \dT g_2 \,dS.
	\end{equation}
	Similarly,
	\begin{align}
	 \inti{\dT f} & = 0,  {\label{FTC1}} \\
	 \inti{\dT f_1 \cdot f_2} & = - \inti{f_1 \cdot \dT f_2}. {\label{FTC2}}
	\end{align}
	where $ f, f_1, f_2 : \Omega \longmapsto \mathbb R $.
\end{lm}

\begin{pf}
	Without lost of generality, we show the lemma for $ \dT = \varphi_1 \cdot \nabla = x_3 \partial_{x_2} - x_2 \partial_{x_3} $.  To do so, we introduce the following coordinate representation of $ \Gamma $, 
	\begin{equation*}
		\begin{aligned}
			X = \left(\begin{array}{c}
				x_1 \\x_2\\x_3
			\end{array}\right) = \left( \begin{array}{c}
				r \\ \sqrt{1-r^2} \cos{\theta} \\ \sqrt{1-r^2} \sin{\theta}
			\end{array} \right) \in \Gamma,
		\end{aligned}
	\end{equation*}
	with $ \theta \in \left\lbrack 0, 2\pi \right) $, $-1 \leq r \leq 1 $. Then it holds
	\begin{equation*}
		\dT = - \partial_\theta, ~ \text{and} ~ \det(D_{r,\theta}X^\top D_{r,\theta}X) = 1,
	\end{equation*}
	Therefore, 
	\begin{equation*}
		\intb{\dT g} = \int_{-1}^{1} \int_{0}^{2\pi} \left( - \partial_\theta g\right) \sqrt{\det(D_{r,\theta}X^\top D_{r,\theta}X)} \, d\theta \,dr  = 0.
	\end{equation*}
	This finishes the proof.
\end{pf}


%
%
%
%
%
%
%
%
%
%
%
%
%
%

\subsection{Embedding Theories and the Stokes Problem}


We start with the trace theorem. The following is from \cite{Adams2003}. 

%

\begin{lm}[Trace Theory]
	Let $ \Omega $ be a bounded domain in $ \mathbb{R}^3 $ with smooth boundary $ \Gamma $. Suppose that $ kp < 3 $ and $ p \leq q \leq p^* = 2p/(3-kp) $. The trace operator $ Tr: W^{k,p}(\Omega) \rightarrow L^q(\Gamma) $ is bounded. Moreover, for $ u \in W^{k,p}(\Omega) $,
	\begin{equation*}
		\norm{u}{L^q(\Gamma)} \lesssim \norm{u}{W^{k,p}(\Omega)}.
	\end{equation*}
	If $ kp = 3 $, then the above relation holds for $ p \leq q < \infty $.
	In particular,
	\begin{equation}\label{trace1}
		\norm{u}{\bLnorm} \lesssim \norm{u}{\Hnorm{1/2}} \lesssim \norm{u}{\Lnorm}^{1/2} \norm{u}{\Hnorm{1}}^{1/2}.
	\end{equation}
\end{lm}
\begin{pf}
See \cite[Theorem 5.36]{Adams2003} and \cite{Brezis2001}.
\end{pf}

Meanwhile, we shall record the following form of trace theorem.
\begin{lm}[Trace Theory: Fractional Sobolev space]
	For $ n \geq 1 $, and $ u \in \Hnorm{n} $, it holds,
	\begin{equation}\label{trace}
		\norm{u}{\bHnorm{n-1/2}} \lesssim \sum_{j=0}^{n-1}\norm{\dT^{j} u}{\Hnorm{1}}.
	\end{equation}
\end{lm}

\begin{pf}
	We shall apply the following fact about the trace operator (see \cite[Theorem 1]{Ding1996}): the trace operator $ Tr : H^s(\Omega) \rightarrow H^{s-1/2}(\partial \Omega) $ is bounded for $ 1/2 < s < 3/2 $. Then \eqref{trace} is a consequence of the fact the rank of $ \lbrace \varphi_1, \varphi_2, \varphi_3 \rbrace $ is equal to two and therefore
	\begin{equation*}
		\norm{u}{\bHnorm{n-1/2}} \lesssim \sum_{j=0}^{n-1} \norm{\dT^j u}{\bHnorm{1/2}}.
	\end{equation*}
\end{pf}

In the meantime, we shall record the classical regularity theory for the Stokes problem. 
Consider the Stokes problem, 
\begin{equation}\label{StokesP}
	\begin{cases}
		 - \mu \Delta u  + \nabla P = f & \text{in} ~ \Omega,\\
		 \dv u = g & \text{in} ~ \Omega,\\
		 u = h & \text{on} ~ \Gamma.
	\end{cases}
\end{equation}
$ \Omega $ is a domain with a smooth boundary $ \Gamma $. 
The following is from \cite{Wang2014,Wang2016}, 
\begin{lm}\label{lm:StokesP}
	Let $ n \geq 2 $. If $ f \in H^{n-2}(\Omega), g \in H^{n-1}(\Omega), h \in H^{n-1/2}(\Gamma) $ be given such that 
	\begin{equation*}
		\inti{g} = \intb{h\cdot\vec{n}},
	\end{equation*}
	then there exists unique $ u \in H^{n}(\Omega), P\in H^{n-1}(\Omega)\text{(up to constants)} $ solving \eqref{StokesP}. Moreover,
	\begin{equation}
		\norm{u}{\Hnorm{n}} + \norm{\nabla P}{\Hnorm{n-2}} \lesssim \norm{f}{\Hnorm{n-2}} + \norm{g}{\Hnorm{n-1}} + \norm{h}{\bHnorm{{n-1/2}}}.
	\end{equation}
\end{lm}

\begin{pf}
See \cite{Ladyzhenskaya1969,Temam1984}. 	
\end{pf}

\subsection{Perturbed Formulation and Poincar\'e Inequality}\label{sec:perturbedform}

We aim to study the stability of \eqref{SNS3}. It is assumed the initial data $ (\rho_0, u_0) $ in \eqref{Ini} satisfies, 
\begin{equation}\label{Ini1}
	\begin{aligned}
		& \int_\Omega \rho_0 \,dx = \int_\Omega \bar \rho \,dx, \int_\Omega \rho_0 u_0 \cdot \varphi_i \,dx = \int_\Omega \bar \rho \bar u \cdot \varphi_i\,dx ~~ i = 1,2,3,
	\end{aligned}	
\end{equation}
where $ (\bar \rho, \bar u ) $ is given in \eqref{SNS3} with some $ \bar \omega > 0 $. 
In particular, it is assumed $ \bar \rho(0) > 0 $ so that $ \bar \rho $ admits uniform lower and upper bounds. Moreover, $ \abs{\nabla^k \bar u}{} \lesssim \bar\omega $ for any $ k \geq 0 $.

Now it is time to introduce the perturbed formulation of \eqref{ICNS} around the steady solution \eqref{SNS3}. Define $ q : = \rho - \bar \rho, v : = u - \bar u $. 
Also, write 
\begin{equation}
	\rho^\gamma - \bar \rho^\gamma = \gamma\bar\rho^{\gamma-1} (\rho - \bar\rho) + \gamma(\gamma-1)\int_{\bar\rho}^{\rho} (\rho - y) y^{\gamma-2}\,dy : =\gamma \bar\rho^{\gamma-1} q + \mathcal{R}.\\
\end{equation}
Notice, from the definition of $ \mathcal R $, 
\begin{equation*}
	\begin{aligned}
	& \nabla \mathcal R = \gamma ( \rho^{\gamma-1} - \bar\rho^{\gamma-1} ) \nabla \rho - \gamma(\gamma-1)\bar\rho^{\gamma-2}(\rho-\bar\rho) \nabla\bar\rho \\
	& ~~~~~ = \gamma ( (\gamma-1) \bar\rho^{\gamma-2} q + \mathcal O(2) )(\nabla q + \nabla\bar \rho )  - \gamma(\gamma-1) \bar\rho^{\gamma-2} q \nabla\bar \rho \\
	& ~~~~~ = \gamma(\gamma-1) \bar\rho^{\gamma-2} q \nabla q + (\nabla q + \nabla\bar \rho ) \mathcal O(2),
	\end{aligned}
\end{equation*}
is at least quadratic in $ q, \nabla q $. 
From \eqref{ICNS} and \eqref{SNS}, we easily derive the system of the perturbation variables $ ( q, v ) $,
\begin{equation}\label{PtNS}
	\begin{cases}
		\dt q  + \dv (\bar \rho v) + ( v + \bar u ) \cdot \nabla q = G_1 & \text{in} ~ \Omega,\\
		(\bar \rho + q) \dt v  + \bar \rho v \cdot \nabla v + \gamma \bar\rho\nabla( \bar\rho^{\gamma-2} q) - \dv \mathbb S(v) = F_2  + G_2 & \text{in} ~ \Omega,\\
		v \cdot \vec{n}, ~ \vec{\tau} \cdot \mathbb S(v)\vec{n} = 0 &\text{on} ~ \Gamma, \\
	\end{cases}
\end{equation}
where 
\begin{equation}\label{Nonlinear-1}
	\begin{aligned}
		& G_1 = - q \dv v,\\
		& G_2 = - \nabla \mathcal R - \left( q \bar u \cdot \nabla v + q v \cdot \nabla \bar u + q v \cdot \nabla v \right), \\
		& F_2 = - \dfrac{\gamma}{\gamma-1} q \nabla \bar\rho^{\gamma-1} - ( \bar\rho \bar u\cdot \nabla v + \bar \rho v \cdot \nabla \bar u + q\bar u \cdot \nabla \bar u )\\
		& ~~~~~ = -\bar\rho \bar u \cdot \nabla v - \bar\rho v\cdot \nabla\bar u.\\
	\end{aligned}
\end{equation}
It should be noticed in the definition of $ F_2 $, it has been used the fact
\begin{equation*}
	\bar u \cdot \nabla \bar u + \dfrac{\gamma}{\gamma-1} \nabla \bar\rho^{\gamma-1} = 0.
\end{equation*}
In addition, we have the following Poincar\'e inequality for $ q $ and $ v $. 
\begin{lm}[Poincar\'e Inequality for $ q, v $] For $ i, j \geq 0 $,
	\begin{align}
		& \inti{\abs{\dT^i\dt^j q}{2}} \lesssim \inti{\abs{\nabla\dT^i\dt^j q}{2}}, \label{qPoincare} \\
		& \inti{\abs{\dT^i\dt^j v}{2}} \lesssim \inti{\abs{\nabla\dT^i\dt^j v}{2}}. \label{vPoincare}
	\end{align}
\end{lm}

\begin{pf}
	For $ i = 0 $, it is a direct consequence of \eqref{Ini1} that
	\begin{equation*}
		\inti{\dt^j q} = 0.
	\end{equation*}
	For $ i \geq 1 $, from \eqref{FTC1},
	\begin{equation*}
		\inti{\dT^i \dt^j q} = \inti{\dT^i \dt^j v} = 0.
	\end{equation*}
	Then \eqref{qPoincare}, \eqref{vPoincare} follow from the standard Poincar\'e inequality as well as \eqref{poincareball}.
\end{pf}

In terms of $ (q, v) $, our stability theorem can be stated as follows. 
\begin{thm}[Main Theorem in Perturbed Variables] The steady rigid motion \eqref{SNS3} with $ \abs{\bar\omega}{} < \epsilon_0 $ for some $ \epsilon_0> 0 $ (given in Lemma \ref{lm:total-Energy}) is nonlinearly stable. 

In particular, there is a constants $ \epsilon_1>0 $ (given in Lemma \ref{lm:total-energy}) such that the classical solution $(q, v)$ to \eqref{PtNS} exists globally with the given initial data $ (q_0, v_0) $ satisfying \eqref{Ini1} and $ \mathfrak{\bar E}_L(0), \mathcal{\bar E}_L(0) < \epsilon_1 $ (defined in \eqref{total-Energy} and \eqref{total-energy}) for $ L \geq 3 $. Moreover, the following inequalities hold for the energy functionals,
\begin{equation}
	e^{\sigma t}\mathcal{\bar E}_L(t) \leq \epsilon_0, ~ e^{\sigma t} \mathfrak{\bar E}_L(t) \leq \epsilon_1,
\end{equation}
for some positive constant $ \sigma > 0 $. Consequently, the initial perturbation would decay to zero as time grows up and therefore $ (\rho, u)(x,t) \rightarrow (\bar\rho, \bar u)(x) $ as $ t \rightarrow +\infty $.  
\end{thm}

\begin{rmk}
	Roughly speaking, the energy functional $ \mathcal{\bar E}_L $ consists of the $ \Hnorm{s} $ norms of $ \dt^l q, \dt^l v $ for some $ s, l \geq 0 $. On the other hand, the functional $ \mathfrak{\bar E}_L $ consists of the anisotropic $ \Hnorm{s} $ norms. In particular, $ \mathfrak{\bar E}_L $ mainly includes $ \Lnorm $ norms of the tangential and interior derivatives. However, the relation $ \mathfrak{\bar E}_L \leq \mathcal{\bar E}_L $ does not hold as one can easily check.
\end{rmk}

\section{Energy Estimates}\label{sec:enerest}

To investigate the stability theory, we shall build some blocks in order to describe the propagation of the initial regularities. Thought out this section, the nonnegative integers $ l,m,n $ are not specified, and would be addressed with appropriate values in the next section.

\subsection{On Temporal and Lower Order Spatial Derivatives}\label{sec:lowerorderest}

To begin with, it is shown the estimate of the temporal derivatives of $ (q, v) $ in this section. Moreover, we present the estimates of the first tangential derivatives and the estimates of normal derivatives through a Stokes problem in terms of the differential frames introduced in Section \ref{sec:sphdiffframes}.

Applying the temporal derivative $ \dt^l $ for $ l = 0,1,\cdots $ to \eqref{PtNS}, it results in the following system,
\begin{equation}\label{PtNStl}
	\begin{cases}
		\dt^{l+1} q + \dv (\bar \rho \dt^l v) + (v + \bar u) \cdot \nabla \dt^l q = G_1^{l} & \text{in} ~ \Omega,\\
		(\bar\rho + q)\dt^{l+1} v + \bar \rho v \cdot \nabla \dt^l v + \gamma \bar\rho \nabla (\bar\rho^{\gamma-2}\dt^l q) - \dv \mathbb{S} (\dt^l v) = F_2^l + G_2^{l} & \text{in} ~ \Omega,\\
		\dt^l v\cdot \vec{n}, ~ \vec{\tau} \cdot \mathbb{S}(\dt^l v)\vec{n} = 0 & \text{on} ~ \Gamma,
	\end{cases}
\end{equation}
where by using the Leibniz's rule,
\begin{equation}\label{Nonlinear-2}
	\begin{aligned}
		& G_1^{l} = \dt^l G_1 - \sum_{j=0}^{l-1} C_{j,l} \dt^{l-j} v \cdot \nabla \dt^j q ,\\
		& G_2^{l} = \dt^l G_2 - \sum_{j=0}^{l-1} C_{j+1,l} \dt^{j+1} q \dt^{l-j} v - \bar \rho \sum_{j=0}^{l-1} C_{j+1,l} \dt^{j+1} v \cdot \nabla \dt^{l-j-1} v, \\
		& F_2^l = \dt^l F_2.
	\end{aligned}
\end{equation}



We record the following energy identity.
\begin{lm}\label{lm:BasicEnergy} For any smooth solution $ (\dt^l q,\dt^l v) $ to \eqref{PtNStl}, it holds the following energy identities, for any integer $ l \geq 0 $, 
	\begin{align}
			& \dfrac{d}{dt}\left\lbrace \dfrac{1}{2} \inti{(\bar\rho+q) \abs{\dt^l v}{2}} + \dfrac{\gamma}{2} \inti{\bar\rho^{\gamma-2}\abs{\dt^l q}{2} } \right\rbrace {\nonumber}\\
			& ~~~~ + \inti{ \left( \dfrac{\mu}{2} \abs{\nabla \dt^l v + \nabla \dt^l v^\top}{2} + \lambda \abs{\dv \dt^l v}{2} \right) }  {\nonumber}\\
			& ~~ =  \dfrac{\gamma}{2} \inti{ \abs{\dt^l q}{2} (v + \bar u) \cdot \nabla \bar\rho^{\gamma-2} } + \dfrac{\gamma}{2} \inti{\bar\rho^{\gamma-2}\abs{\dt^l q}{2} \dv( v + \bar u) }   {\label{EneId01}}\\
			& ~~~~ + \gamma \inti{\bar\rho^{\gamma-2} \dt^l q G_1^{l}} + \dfrac{1}{2}\inti{ \left(\dt q + \dv(\bar\rho v) \right) \abs{\dt^l v}{2} }  {\nonumber}\\
			& ~~~~ +  \inti{(F_2^l + G_2^l) \cdot \dt^l v}. {\nonumber}
	\end{align}
\end{lm}

\begin{pf}
	Take inner product of $ \eqref{PtNStl}_2 $ with $ \dt^l v $ and then record the resulting after integration by parts,
	\begin{equation}\label{EneId1}
		\begin{aligned}
			& \dfrac{d}{dt} \dfrac{1}{2} \int_\Omega (\bar\rho + q) \left|\dt^l v\right|^2\,dx  - \gamma \inti{\dv (\bar\rho \dt^l v)\bar\rho^{\gamma-2} \dt^l q} \\
			& ~~~~~~ + \inti{ \left( \dfrac{\mu}{2} \abs{\nabla \dt^l v + \nabla \dt^l v^\top}{2} + \lambda \abs{\dv \dt^l v}{2} \right) } \\
			& ~~~~~~ + \gamma \intb{ \bar\rho^{\gamma-1} \dt^l q \dt^l v \cdot \vec{n} } + \dfrac{1}{2}\intb{ \bar\rho \abs{\dt^l v}{2} v\cdot \vec{n} } - \intb{ \dt^l v\cdot \mathbb S(\dt^l v) \vec{n} }\\
			& ~~~~ = \dfrac{1}{2} \inti{\dt q \abs{\dt^l v}{2} } + \dfrac{1}{2}\inti{\dv (\bar\rho v) \abs{\dt^l v}{2} } +  \inti{(F_2^l + G_2^l) \cdot \dt^l v}.
		\end{aligned}
	\end{equation}
	Similar as before, from \subeqref{PtNStl}{3}, $ \dt^k v $ is in the tangential direction of $ \Gamma $. Therefore, from \subeqref{PtNStl}{3} and \subeqref{PtNS}{3}, the boundary integrals in \eqref{EneId1} vanish. Meanwhile, using \subeqref{PtNStl}{1}, it holds
	\begin{equation*}
		\begin{aligned}
			& - \gamma \inti{\dv(\bar\rho \dt^l v)\bar\rho^{\gamma-2} \dt^l q } = \gamma \inti{\left(\dt^{l+1} q + (v+\bar u) \cdot \nabla \dt^l q - G_1^l \right)\bar\rho^{\gamma-2} \dt^l q }\\
			& = \dfrac{d}{dt}\dfrac{\gamma}{2}\inti{ \bar\rho^{\gamma-2} \abs{\dt^l q}{2} } - \gamma \inti{\bar\rho^{\gamma-2} \dt^l q G_1^{l}}  - \dfrac{\gamma}{2} \inti{ \abs{\dt^l q}{2} (v+\bar u) \cdot \nabla \bar\rho^{\gamma-2} }\\
			& ~~~~~~ - \dfrac{\gamma}{2} \inti{\bar\rho^{\gamma-2}\abs{\dt^l q}{2} \dv (v+\bar u) }, 
		\end{aligned}
	\end{equation*}
	where we have applied integration by parts and the boundary condition \subeqref{PtNS}{3}. Thus \eqref{EneId01} follows after chaining the above identity. 
\end{pf}


\begin{lm} For any smooth solution $ (\dt^l q, \dt^l v) $ to \eqref{PtNStl}, any $ 0 < \omega < 1 $, the following estimate holds,
	\begin{align}\label{Prop:l-estimate}
			& \dfrac{d}{dt}\left\lbrace \dfrac{1}{2} \inti{(\bar\rho+q) \abs{\dt^l v}{2}} + \dfrac{\gamma}{2} \inti{\bar\rho^{\gamma-2}\abs{\dt^l q}{2} } \right\rbrace + \dfrac{\mu}{2} \inti{ \abs{\nabla \dt^l v}{2} } {\nonumber} \\
			& ~~ + \lambda \inti{ \abs{\dv \dt^l v}{2}}  \lesssim (\omega + a_{l,1} ) \inti{\abs{\nabla\dt^l q}{2}} + a_{l,2} \inti{ \abs{\nabla\dt^l v}{2}} {\nonumber} \\ 
			& ~~ + \sum_{j=0}^{[l/2]} b_{l,j}\left( \inti{\abs{\nabla\dt^{l-1-j}q}{2}} + \inti{\abs{\nabla\dt^{l-j}v}{2}} \right) {} \\
			& ~~ + C_\omega \inti{\abs{G_1^{l}}{2}} + \inti{\abs{G_2^{l}}{2}},  {\nonumber}
	\end{align}
	where, for some integer $ a > 0 $, 
	\begin{equation*}
		\begin{aligned}
			& a_{l,1} = \bar\omega^a + \norm{v}{\supnorm}^2 + \norm{\nabla v}{\supnorm},   \\
			& a_{l,2} = \bar\omega^a + \norm{\dt q}{\supnorm} + \norm{\nabla v}{\supnorm} + \norm{v}{\supnorm}^2, \\
			& b_{l,j} = \norm{\dt^{j+1} v}{\supnorm}^2 + \norm{\dt^j q}{\supnorm}^2 ,~~~ 0 \leq j \leq [l/2].
		\end{aligned}
	\end{equation*}

\end{lm}

\begin{pf}
	After chaining \eqref{korn04}, \eqref{EneId01} and applying Cauchy's inequality, it holds for any $ 0 < \delta < 1 $ and $0 < \omega < 1 $, 
	\begin{align*}
			& \dfrac{d}{dt}\left\lbrace \dfrac{1}{2} \inti{(\bar\rho+q) \abs{\dt^l v}{2}} + \dfrac{\gamma}{2} \inti{\bar\rho^{\gamma-2}\abs{\dt^l q}{2} } \right\rbrace \\
			& ~~~~~~ + \dfrac{\mu}{2} \inti{ \abs{\nabla \dt^l v}{2} } + \lambda \inti{ \abs{\dv \dt^l v}{2}} \\
			& ~~ \lesssim \left( \omega + \bar\omega^a +\norm{v}{\supnorm}^2 + \norm{\nabla v}{\supnorm}  \right)\inti{\abs{\dt^l q}{2} } \\
			& ~~~~ + \left( \delta + \norm{\dt q}{\supnorm} + \norm{\nabla v}{\supnorm} + \bar\omega^a \norm{v}{\supnorm} \right) \inti{ \abs{\dt^l v}{2}}\\
			& ~~~~ + \sum_{j=0}^{l-1} \inti{ \abs{\dt^j q}{2} \abs{\dt^{l-j}v}{2}} + C_\delta \inti{ \abs{F_2^l}{2} } \\
			& ~~~~ + C_\omega \inti{\abs{G_1^{l}}{2}} + C_\delta \inti{\abs{G_2^{l}}{2}},
	\end{align*}
	where it has been used the fact $ \norm{\bar u }{\supnorm}, \norm{\nabla \bar u}{L^\infty(\Omega)}, \norm{\nabla \bar \rho}{\supnorm} \lesssim \bar\omega^a $ for some integer $ a > 0 $.
	Meanwhile
	\begin{equation*}
		\begin{aligned}
			& \sum_{j=0}^{l-1} \inti{ \abs{\dt^j q}{2} \abs{\dt^{l-j}v}{2}} \lesssim \sum_{j=0}^{[l/2]} \left( \norm{\dt^{j+1} v}{\supnorm}^2 + \norm{\dt^j q}{\supnorm}^2 \right) \\
			& ~~~~~~ \times \left( \inti{\abs{\dt^{l-1-j}q}{2}} + \inti{\abs{\dt^{l-j}v}{2}} \right). 
		\end{aligned}
	\end{equation*}
%
On the other hand, from the definition of $ F_2^l $, 
	\begin{equation*}
		\inti{\abs{F_2^l}{2}} \lesssim \bar\omega^a \left( \inti{\abs{\nabla \dt^l v}{2}} + \inti{\abs{\dt^l v}{2}} \right).
	\end{equation*}
	Then by chaining these estimates, together with \eqref{vPoincare}, \eqref{qPoincare} and choosing an appropriately small $ \delta > 0 $, \eqref{Prop:l-estimate} holds. 
\end{pf}

Here and after, $ a > 0 $ would denote an integer which might be different from line to line. Next lemma is concerning the estimate of the spatial derivative. 
\begin{lm} Under the same assumptions as in Lemma \ref{lm:BasicEnergy}, 
	\begin{align}
			& \dfrac{d}{dt}\biggl\lbrace \dfrac{\mu}{4} \inti{\abs{\nabla \dt^l v + \nabla \dt^l v^T}{2}} + \dfrac{\lambda}{2}  \inti{\abs{\dv \dt^l v}{2}} {\nonumber}\\
			& ~~ - \gamma \inti{\bar\rho^{\gamma-2}\dv(\bar\rho\dt^lv)\dt^l q}  \biggr\rbrace  + \inti{(1 + q ) \abs{\dt^{l+1}v}{2}} \lesssim (\bar\omega^a {\label{l-est2}} \\
			& ~~ + \norm{v}{\supnorm}^2)  \inti{\abs{\nabla\dt^l q}{2}}  + (1 + \bar\omega^a + \norm{v}{\supnorm}^2  )\inti{\abs{\nabla \dt^l v}{2}}  {\nonumber}\\
			& ~~ + \inti{\abs{G_1^l}{2}} + \inti{\abs{G_2^l}{2}}.{\nonumber}
	\end{align}
\end{lm}

\begin{pf}
	Take inner product of \subeqref{PtNStl}{2} with $ \dt^{l+1} v $ and then record the resulting after integration by parts. Similar arguments as in the previous lemma yield,
	\begin{align*}
			& \dfrac{d}{dt} \left\lbrace \dfrac{\mu}{4}\inti{\abs{\nabla \dt^l v + \nabla \dt^l v^\top}{2}} + \dfrac{\lambda}{2}\inti{\abs{\dv \dt^l v}{2}} \right\rbrace + \inti{(\bar\rho + q ) \abs{\dt^{l+1}v}{2}}\\
			& ~~~~~~ \underbrace{- \gamma \inti{\dv(\bar\rho \dt^{l+1}v) \bar\rho^{\gamma-2}\dt^lq} }_{(i)} = \underbrace{ - \inti{\bar\rho v \cdot \nabla \dt^l v \dt^{l+1}v}}_{(ii)} \\
			& ~~~~~~ + \inti{(F_{2}^l+G_2^l)\cdot\dt^{l+1}v}.
	\end{align*}
Meanwhile, together with \eqref{qPoincare} and \eqref{vPoincare},
\begin{align*}
		& (i) = \dfrac{d}{dt} \left\lbrace -\gamma \inti{\bar\rho^{\gamma-2}\dv(\bar\rho\dt^lv)\dt^l q} \right\rbrace  + \underbrace{ \gamma\inti{\bar\rho^{\gamma-2}\dv(\bar\rho\dt^l v)\dt^{l+1}q} }_{(iii)}, \\
		& (iii) \lesssim (1+\bar\omega^a)\inti{\abs{\nabla\dt^l v}{2}} + \inti{\abs{\dt^{l+1} q}{2}}\\
		& ~~~~ \lesssim (1+\bar\omega^a ) \inti{\abs{\nabla\dt^l v}{2}} + (\norm{v}{\supnorm}^2 + \bar\omega^a) \inti{\abs{\nabla\dt^l q}{2}} + \inti{\abs{G_1^l}{2}},\\
		& (ii) \lesssim \delta \inti{\abs{\dt^{l+1} v}{2}} + C_\delta \norm{v}{\supnorm}^2 \inti{\abs{\nabla\dt^l v}{2}},
\end{align*}
where it has been making use of the fact, from \subeqref{PtNStl}{1} and \eqref{vPoincare}
\begin{equation*}
	\begin{aligned}
		& \inti{\abs{\dt^{l+1}q}{2}} \lesssim (1+\bar\omega^a) \inti{\abs{\nabla\dt^l v}{2}} \\
		& ~~~~~~~~ + (\norm{v}{\supnorm}^2 + \norm{\bar u}{\supnorm}^2) \inti{\abs{\nabla\dt^l q}{2}} + \inti{\abs{G_1^l}{2}}.
	\end{aligned}
\end{equation*}
Similarly, 
\begin{equation*}
	\begin{aligned}
		& \inti{F_{2}^l\cdot\dt^{l+1}v} \lesssim \delta \inti{\abs{\dt^{l+1}v}{2}} + C_\delta \bar\omega^a \inti{\abs{\nabla \dt^l v}{2}},\\
		& \inti{G_2^l\cdot\dt^{l+1}v}\lesssim \delta \inti{\abs{\dt^{l+1} v}{2}} + C_\delta \inti{\abs{G_2^l}{2}}.
	\end{aligned}
\end{equation*}
Thus chaining these estimates with an appropriately small $ \delta > 0 $ leads to \eqref{l-est2}.
\end{pf}

As the consequence of \eqref{Prop:l-estimate} and \eqref{l-est2}, we have the following estimates on the temporal derivatives. 
\begin{prop}[Temporal Derivatives]\label{Prop:temporalderivative}
	Denote $ \Lambda_l = \Lambda_l(\cdot) $ as a polynomial of the following quantities
	$$ \bar\omega, \norm{v}{\supnorm},\norm{\nabla v}{\supnorm}, \norm{\dt q}{\supnorm}, \sum_{j=0}^{[l/2]} ( ||\dt^{j+1}v||_{\supnorm} + ||\dt^j q||_{\supnorm} ) , $$ 
	with the property $ \Lambda_l(0) = 0 $. Then it shall hold the following estimate. For any $ 0 < \omega < 1 $,
	\begin{align}\label{Proposition-time}
			& \dfrac{d}{dt} \biggl\lbrace \inti{(\bar\rho + q)\abs{\dt^l v}{2}} +  c\inti{\abs{\nabla\dt^l v + \nabla\dt^l v^T}{2}} + c\inti{\abs{\dv \dt^l v}{2}} {\nonumber}\\ 
			& + \inti{\bar\rho^{\gamma-2}\abs{\dt^l q}{2}} - c\inti{\bar\rho^{\gamma-2} \dv(\bar\rho \dt^l v) \dt^l q} \biggr\rbrace {\nonumber}\\
			& + (1-c)\inti{\abs{\nabla\dt^l v}{2}} + \inti{\abs{\dv\dt^l v}{2}} + c \inti{(1+q) \abs{\dt^{l+1} v}{2}}\\
			& \lesssim ( \omega + \Lambda_l ) \inti{\abs{\nabla \dt^l q}{2}} + \Lambda_l \inti{\abs{\nabla \dt^l v}{2}}  + \Lambda_l \sum_{j=0}^{[l/2]} \biggl( \inti{\abs{\nabla\dt^{l-1-j}q}{2}} {\nonumber}\\
			& ~~  + \inti{\abs{\nabla\dt^{l-j}v}{2}} \biggr) + (1 +  C_\omega ) \inti{\abs{G_{1}^l}{2}} + \inti{\abs{G_{2}^l}{2}}, {\nonumber}
	\end{align}
	where $ 0 < c < 1 $ is a positive constant such that 
	\begin{equation}\label{little-c}
		\begin{aligned}
		& \inti{\abs{\dt^l v}{2}} +  \inti{\abs{\nabla\dt^l v}{2}} + \inti{\abs{\dt^l q}{2}}\\
		& \lesssim \inti{(\bar\rho + q)\abs{\dt^l v}{2}} +  c\inti{\abs{\nabla\dt^l v + \nabla\dt^l v^T}{2}} + c\inti{\abs{\dv v}{2}}\\ 
		& ~~~~ + \inti{\bar\rho^{\gamma-2}\abs{\dt^l q}{2}} - c\inti{\bar\rho^{\gamma-2} \dv(\bar\rho \dt^l v) \dt^l q}.
		\end{aligned}
	\end{equation}
\end{prop}

\begin{pf} This is a direct consequence of the linear combination $ c \times \eqref{l-est2} + \eqref{Prop:l-estimate} $. 
	We only show the the choice of $ c $ can be justified. By applying the Cauchy's inequality and \eqref{korn01}, \eqref{vPoincare},
	\begin{equation*}
		\begin{aligned}
		& \inti{\abs{\nabla\dt^l v}{2}} - \inti{\abs{\dt^l v}{2}} - \inti{\abs{\dt^l q}{2}} \\
		& \lesssim \inti{\abs{\nabla\dt^l v + \nabla\dt^l v^T}{2}}-\inti{\bar\rho^{\gamma-2} \dv(\bar\rho \dt^l v) \dt^l q}.
		\end{aligned}
	\end{equation*}
	Together with the fact $\inf_{x\in\Omega} \bar\rho > 0 $, \eqref{little-c} holds after choosing $ c > 0 $ sufficiently small. 
\end{pf}

Notice, in Proposition \ref{Prop:temporalderivative}, the estimate \eqref{Proposition-time} contains the term $ \inti{\abs{\nabla\dt^l q}{2}} $ on the right hand side. To derive an estimate in a consistent form, it is desirable to perform the estimate on the spatial derivatives in the rest of this section.


Here we establish the estimates on the tangential derivatives. Starting with the first order tangential derivative, $ \nabla_T \eqref{PtNStl} $ can be written as,
%
\begin{equation}\label{TPt}
	\begin{cases}
		\dt^{l+1} \dT q + \dv \left(\bar\rho \dt^l \dT v \right) + (v+\bar u)\cdot \nabla\dt^l \dT q = F_{1}^{l,1} + G_1^{l,1} & \text{in} ~ \Omega,\\
		(\bar\rho + q)\dt^{l+1} \dT v + \bar\rho v \cdot \nabla \dt^l \dT v + \gamma \bar\rho \nabla\left( \bar\rho^{\gamma-2} \dt^l \dT q \right)\\
		~~~~~~ - \dv \mathbb{S}(\dt^l \dT v) = F_2^{l,1} + \dT G_2^l + G_2^{l,1} & \text{in} ~ \Omega,\\
		\dt^l \dT v \cdot \vec{n} = - \dt^l v \cdot\dT \vec{n}, ~ \dT \left( \vec{\tau} \cdot \mathbb{S}(\dt^l v) \vec{n} \right) = 0 & \text{on} ~ \Gamma,
	\end{cases}
\end{equation}
with
\begin{equation}\label{Nonlinear-3}
	\begin{aligned}
		& G_1^{l,1} = \dT G_1^{l} - \dT v \cdot \nabla\dt^l q - v \cdot \left\lbrack \dT, \nabla\right\rbrack \dt^l q ,\\
		& G_2^{l,1} = - \dT q \dt^{l+1} v - \bar\rho v \cdot \commut{\dT}{\nabla} \dt^l v - \dT(\bar\rho v) \cdot \nabla \dt^l v,\\
		& F_{1}^{l,1} = - \left\lbrack \dT, \dv \right\rbrack (\bar\rho \dt^l v) - \dv(\dT\bar\rho\dt^l v) - \dT \bar u\cdot\nabla\dt^l q - \bar u \cdot  \commut{\dT}{\nabla}\dt^l q ,\\
		& F_{2}^{l,1} = \dT F_2^l - \dT \bar\rho \dt^{l+1} v -\gamma \dT \bar\rho \nabla(\bar\rho^{\gamma-2}\dt^lq) - \gamma\bar\rho \commut{\dT}{\nabla} (\bar\rho^{\gamma-2}\dt^l q) \\
		& ~~~~~~ - \gamma \bar\rho \nabla( \dT \bar\rho^{\gamma-2} \dt^l q)  + \commut{\dT}{\dv} \mathbb{S}(\dt^l v) + \dv ( \commut{\dT}{\mathbb{S}}(\dt^l v) ) .\\
	\end{aligned}
\end{equation}

Take inner product of \subeqref{TPt}{2} with $ \dt^l \dT v $ and then record the resulting after integration by parts, 

\begin{equation}\label{l1-est01}
	\begin{aligned}
		& \dfrac{d}{dt} \dfrac{1}{2} \inti{(\bar\rho + q) \abs{\dt^l \dT v}{2}} \underbrace{- \gamma\inti{\dv{(\bar\rho \dt^l \dT v)}\bar\rho^{\gamma-2} \dt^l\dT q} }_{(i)}\\
		& ~~~~~~ + \inti{\left(\dfrac{\mu}{2} \abs{ \nabla \dt^l \dT v + \nabla \dt^l \dT v^T }{2} + \lambda \abs{\dv( \dt^l \dT v)}{2}  \right)}  \\
		& = \underbrace{ \dfrac{1}{2} \inti{ \dt q \abs{\dt^l \dT v}{2}} + \dfrac{1}{2}\inti{\dv(\bar\rho v) \abs{\dt^l \dT v}{2}} }_{(ii)} \\
		& ~~~~~~ 
		 - \gamma \underbrace{ \intb{\bar\rho^{\gamma-1}\dt^l \dT q \dt^l\dT v\cdot \vec{n} } }_{(iii)}  + \underbrace{\intb{\dt^l\dT v \cdot\mathbb{S}(\dt^l\dT v) \vec{n}}}_{(iv)} \\
		& ~~~~~~ + \underbrace{\inti{(F_2^{l,1} + G_2^{l,1}) \cdot \dt^l \dT v}.}_{(v)} + \underbrace{\inti{\dT G_2^l \cdot \dt^l \dT v }}_{(vi)},
	\end{aligned}
\end{equation}
where we have used the boundary condition \subeqref{PtNS}{3}. 
Meanwhile, from \subeqref{TPt}{1}, \subeqref{PtNS}{3}, 
\begin{align*}
		& (i) = \gamma \inti{\bar\rho^{\gamma-2} \dt^l \dT q \cdot ( \dt^{l+1}\dT q + (v+\bar u)\cdot \nabla\dt^l \dT q - F_{1}^{l,1} - G_{1}^{l,1} ) } \\
		& ~~~~ = \dfrac{d}{dt}\dfrac{\gamma}{2} \inti{\bar\rho^{\gamma-2}\abs{\dt^l\dT q}{2}} - \underbrace{\dfrac{\gamma}{2}\inti{ \dv(\bar\rho^{\gamma-2} (v+\bar u) )\abs{\dt^l\dT q}{2} }}_{(vii)}\\
		& ~~~~~~ - \underbrace{ \gamma\inti{(F_1^{l,1} + G_{1}^{l,1} )\cdot \bar\rho^{\gamma-2}\dt^l \dT q}.}_{(viii)} 
\end{align*}
Applying Cauchy\'s inequality and Poincar\'{e} inequality as follows,
\begin{align*}
		& (vii) \lesssim \left( \bar\omega^a + \norm{v}{\supnorm}^2 + \norm{\nabla v}{\supnorm} \right) \inti{\abs{\dT \dt^l q}{2}}, \\
		& (viii) \lesssim \omega \inti{\abs{\dT \dt^l q}{2}} + C_{\omega} \left(  \inti{\abs{F_1^{l,1}}{2}} + \inti{\abs{G_{1}^{l,1}}{2}} \right), \\
		& (ii) \lesssim \left( \bar\omega^a + \norm{\dt q}{\supnorm} +\norm{v}{\supnorm}^2 + \norm{\nabla v}{\supnorm} \right) \inti{\abs{\dT \dt^l v}{2}}, \\
		& (v) \lesssim  \omega \inti{\abs{F_{2}^{l,1}}{2}}  + C_\omega \left(\inti{\abs{\dT\dt^l v}{2}}+ \inti{\abs{G_{2}^{l,1}}{2}}\right).
\end{align*}
To estimate the boundary terms, from \subeqref{TPt}{3}, as consequences of the trace theorem \eqref{trace1}, H\"older inequality and \eqref{FTCB2}, \eqref{qPoincare}, \eqref{vPoincare}, it holds
\begin{align*}
		& (iii) = - \intb{\bar\rho^{\gamma-1} \dt^l \dT q \dt^l v \cdot \dT \vec{n} } = \intb{\dt^l q \dT(\bar\rho^{\gamma-1} \dt^l v\cdot \dT \vec{n} )}\\
		& \lesssim \norm{\dt^l q}{\bLnorm}\left((1+\bar\omega) \norm{\dt^l v}{\bLnorm} + \norm{\dT \dt^l v}{\bLnorm} \right)\\
		& \lesssim \norm{\dt^l q}{\Lnorm}^{1/2}\norm{\dt^l q}{\Hnorm{1}}^{1/2} \left( (1+\bar\omega) \norm{\dt^l v}{\Lnorm}^{1/2}\norm{\dt^l v}{\Hnorm{1}}^{1/2} \right. \\
		& \left. ~~~~~~ + \norm{\dT\dt^l v}{\Lnorm}^{1/2}\norm{\dT\dt^l v}{\Hnorm{1}}^{1/2}  \right) \\
		& \lesssim \omega \left( \inti{\abs{\nabla \dt^l q}{2}} + \inti{\abs{\nabla\dT \dt^l v}{2}} \right) \\
		& ~~~~ + C_{\omega}(1+\bar\omega^a)\inti{\abs{ \nabla \dt^l v}{2}},
\end{align*}
On the other hand, on $ \Gamma $, we have the following identities from \subeqref{TPt}{3}
\begin{equation}\label{idi:boundarycon01}
	\begin{aligned}
		& 0 = \dT(\vec{\tau} \cdot \mathbb{S}(\dt^l v) \vec{n} ) = \vec{\tau} \cdot \mathbb{S}(\dt^l \dT v) \vec{n} + \dT \vec{\tau} \cdot \mathbb{S}(\dt^l v) \vec{n} \\
		& ~~~~~~ + \vec{\tau} \cdot \commut{\dT}{\mathbb{S}} (\dt^l v) \vec{n} + \vec{\tau}\cdot \mathbb{S}(\dt^l v) \dT \vec{n},\\
		& \dt^l \dT v = ( \dt^l \dT v \cdot \vec{n} ) \vec{n} + v_{l,1} \vec{\tau} = - (\dt^l v \cdot \dT \vec{n}) \vec{n} + v_{l,1} \vec{\tau},
	\end{aligned}
\end{equation}
with $ \abs{v_{l,1}}{} \lesssim \abs{\dt^l \dT v}{} $. 
Therefore, the calculus on the boundary \eqref{FTCB2} then yields
\begin{align}\label{l1-est02}
		& (iv) = - \intb{(\dt^l v \cdot \dT \vec{n}) \vec{n} \cdot \mathbb{S}(\dt^l \dT v) \vec{n} } + \intb{v_{l,1}\vec\tau \cdot \mathbb{S}(\dt^l \dT v) \vec{n} } {\nonumber} \\
		& = - \intb{(\dt^l v \cdot \dT \vec{n}) \vec{n} \cdot \commut{\mathbb{S}}{\dT}(\dt^l v) \vec{n}} + \intb{ \dT((\dt^l v \cdot \dT \vec{n})\vec{n}) \cdot \mathbb{S}(\dt^l v) \vec{n}  } {\nonumber}\\
		& ~~ + \intb{(\dt^l v \cdot \dT \vec{n}) \vec{n} \cdot \mathbb{S}(\dt^l v) \dT \vec{n}} - \intb{v_{l,1}\dT \vec{\tau} \cdot \mathbb{S}(\dt^l v) \vec{n}} {\nonumber}\\
		& ~~ - \intb{v_{l,1}\vec{\tau} \cdot \commut{\dT}{\mathbb{S}} (\dt^l v) \vec{n} } - \intb{v_{l,1} \vec{\tau}\cdot \mathbb{S}(\dt^l v) \dT \vec{n}}\\
		& \lesssim \intb{\left(\abs{\dt^l v}{} + \abs{\dt^l \dT v}{}\right) \cdot\abs{\nabla\dt^l v}{}} \lesssim \norm{\dt^l v}{\bLnorm}^2 + \norm{\nabla\dt^l v}{\bLnorm}^2 {\nonumber}\\
		& ~~ + \norm{\dT \dt^l v}{\bLnorm}^2 \lesssim \norm{\dt^l v}{\Lnorm} \norm{\dt^l v}{\Hnorm{1}} {\nonumber} \\
		& ~~ + \norm{\nabla \dt^l v}{\Lnorm} \norm{\nabla \dt^l v}{\Hnorm{1}} {\nonumber} \lesssim \omega \inti{\abs{\nabla^2 \dt^l v}{2}} + C_{\omega} \inti{\abs{\nabla\dt^l v}{2}}. {\nonumber}
\end{align}

Indeed, the above calculation indicates the next lemma.
\begin{lm}[Tangential Direction Estimate]For any smooth solution $(\nabla_T \dt^l q, \nabla_T \dt^l v)$ to \eqref{TPt}, the following inequality holds, 
	\begin{align}\label{vql1-est} 
			& \dfrac{d}{dt} \left\lbrace \dfrac{1}{2}\inti{(\bar\rho + q) \abs{\dT \dt^l v}{2}} + \dfrac{\gamma}{2} \inti{\bar\rho^{\gamma-2} \abs{\dT \dt^l q}{2}}  \right\rbrace {\nonumber}\\
			& ~~~~ + \inti{\abs{\nabla \dT\dt^l v}{2}} + \inti{\abs{\dv (\dT\dt^l v)}{2}} {\nonumber}\\
			& \lesssim (\omega + \Lambda_l )\inti{\abs{ \nabla^2 \dt^l v}{2}} + (\omega + C_\omega\bar\omega^a + \Lambda_l ) \inti{\abs{ \nabla \dt^l q}{2}} \\
			& ~~~~ + (1 + \omega + C_\omega  + \Lambda_l ) \inti{\abs{\nabla \dt^l v}{2}} + \omega \bar\omega^2 \inti{\abs{\dt^{l+1} v}{2}} {\nonumber}\\
			& ~~~~ + C_{\omega} \biggl( \inti{\abs{G_{1}^{l,1}}{2}} + \inti{\abs{G_{2}^{l,1}}{2}}  \biggr) + \inti{\abs{G_2^l}{2}}, {\nonumber}
	\end{align}
	for any $ 0 < \omega < 1 $. 
\end{lm}

\begin{pf}
	By summing up the estimates between \eqref{l1-est01} and \eqref{l1-est02}
	, and \eqref{korn01}, we only need to estimate $ F_1^{l,1}, F_2^{l,1} $ and $ (vi) $. Direct calculations show
	\begin{align*}
			& \inti{\abs{F_{1}^{l,1}}{2}} 
			\lesssim (1+\bar\omega^a)\inti{\abs{\nabla\dt^l v}{2}} + \bar\omega^a \inti{\abs{\nabla\dt^l q}{2}}, \\
			& \inti{\abs{F_{2}^{l,1}}{2}} \lesssim \inti{\abs{\dT F_2^l}{2}} + \bar\omega^a \inti{\abs{\dt^{l+1} v}{2}} \\
			& ~~~~~~ + (1+\bar\omega^a) \inti{\abs{\nabla\dt^l q}{2}} + \inti{\abs{\nabla^2 \dt^l v}{2}} \\
			& ~~~~ \lesssim \bar\omega^a \inti{\abs{\nabla\dt^l v}{2}} + \bar\omega^a \inti{\abs{\dt^{l+1} v}{2}} \\
			& ~~~~~~ + (1+\bar\omega^a) \inti{\abs{\nabla\dt^l q}{2}} + (1+\bar\omega^a)\inti{\abs{\nabla^2 \dt^l v}{2}},
	\end{align*}
	where $ a > 0 $ and we have applied \eqref{qPoincare}, \eqref{vPoincare}. Also, as the consequence of \eqref{FTC2},
	\begin{equation*}
			(vi) = - \inti{G_2^l \cdot \dT \dt^l \dT v} \lesssim \delta \inti{\abs{\nabla \dT \dt^l v}{2}} + C_\delta \inti{\abs{G_2^l}{2}}.
	\end{equation*}
	Then by choosing an appropriately small $ \delta > 0 $, \eqref{vql1-est} is proved.
\end{pf}

\begin{rmk}\label{rm:boundarytangentialfields}
	The identities \eqref{idi:boundarycon01} on the boundary and the boundary condition \subeqref{TPt}{3} should be understood as follows.
	$ \tau $ stands for one of the tangential vector fields $ \varphi_1, \varphi_2, \varphi_3  $ (defined in \eqref{korn2}), which are smooth and defined globally on $ \Gamma $. Therefore, $ \dT \tau $ is non-singular and smooth. Moreover, the rank of $ \left\lbrace \varphi_1, \varphi_2, \varphi_3 \right\rbrace $ is equal to two and hence any tangental vector on $ \Gamma $ can be represented by them. To show that $ v_{l,1}\tau $ makes sense and $ \abs{v_{l,1}}{} \lesssim \abs{\dt^l \dT v}{} $, we adopt the following representation of tangential vector fields. Indeed, we claim that any tangential vector fields $ V $ on $ \Gamma $ can be denoted as 
	\begin{equation}\label{rp:tangentialfields}
		V = V_1 \varphi_1 + V_2 \varphi_2 + V_3 \varphi_3,
	\end{equation}
	with $ \abs{V_1}{}, \abs{V_2}{}, \abs{V_3}{} $ bounded by $ \abs{V}{} $. To show this is possible, consider a point $ p = (x_1,x_2,x_3) $ on $ \Gamma $. Without loss of generality, we assume $ x_1 > 1/4 $. Then for a neighbourhood $ W_p $ of $ p $, $ \left\lbrace \varphi_2, \varphi_3 \right\rbrace $ forms a non-degenerate basis of the tangential space and $ \abs{\varphi_2}{}, \abs{\varphi_3}{} > 1/4 $. Then $ V $ can be written as
	$$ V = V_{1,p} \varphi_1 + V_{2,p} \varphi_2 + V_{3,p} \varphi_3 $$ with $ V_{1,p}=0 $ and $ \abs{V_{2,p}}{}, \abs{V_{3,p}}{} \lesssim \abs{V}{} $ inside $ W_p $. Since $ \Gamma $ is compact, one can construct a finite cover $ \left\lbrace W_p \right\rbrace $ of $ \Gamma $ and corresponding $ \left\lbrace p; V_{1,p}, V_{2,p}, V_{3,p}\right\rbrace $. Then a partition of unity argument would yield \eqref{rp:tangentialfields}. Notice, $ V_1, V_2, V_3 $ are not necessarily continuous. 
\end{rmk}

Next, we show how to develop the estimates of the normal derivatives in the spherical differential frame. Just as it is classically done (see, for example, \cite{Wang2016}), we shall derive the ordinary differential equation(ODE) satisfied by $ \dt^l \dN q $. In order to do so, 
taking inner product of \subeqref{PtNStl}{2} with $ N $ yields, 
\begin{equation}\label{NDEst001}
	\begin{aligned}
		& (\bar\rho + q) \dt^{l+1} v \cdot N + \bar\rho v\cdot \nabla\dt^l v\cdot N - ( F_2^l + G_2^l )\cdot N \\
		& ~~~~~~ = -  \gamma\bar\rho\nabla\left( \bar\rho^{\gamma-2} \dt^l q \right) \cdot N + \dv \mathbb{S}(\dt^l v)\cdot N.
	\end{aligned}
\end{equation}
Notice, $ N = (x_1,x_2,x_3)^\top $, and 
\begin{align}\label{NDEst002}
		& \dv \mathbb{S}(\dt^l v)\cdot N = \sum_{1\leq i,j\leq 3} x_i \partial_j \left( \mu \left( \partial_i \dt^l v^j + \partial_{j} \dt^l v^i \right)+ \lambda \delta_{ij} \dv \dt^l v \right) {\nonumber}\\
		& ~~~~ = (\mu + \lambda) \dN \dv \dt^l v + \mu \sum_{1\leq i, j \leq 3} x_i \partial_j \partial_j \dt^l v^i = (2\mu +\lambda)\dN \dv \dt^l v \\
		& ~~~~~~~~ + \mu \sum_{1\leq s,m,n\leq 3} \varphi_s \cdot \nabla ( \epsilon_{smn} \partial_m \dt^l v^n), {\nonumber}
\end{align}
where in the last equality we have substituted the following identity,
\begin{align*}
		& \sum_{1\leq i, j \leq 3} x_i \partial_j \partial_j \dt^l v^i = \sum_{1\leq i, j, m, n \leq 3} x_i(\delta_{jm}\delta_{in}) \partial_m \partial_j \dt^l v^n \\
		& = \sum_{1\leq i,j,m,n\leq 3} x_i (\delta_{jm}\delta_{in} - \delta_{jn}\delta_{im}) \partial_m \partial_j \dt^l v^n + \sum_{1\leq i, j, m, n \leq 3} x_i(\delta_{jn}\delta_{im}) \partial_m \partial_j \dt^l v^n \\
		& = \sum_{1\leq i,j,m,n,s\leq 3} x_i \epsilon_{jis} \epsilon_{smn} \partial_m\partial_j\dt^l v^n + \sum_{1\leq i,j\leq 3} x_i \partial_i \partial_j \dt^l v^j \\
		& = \sum_{1\leq s,m,n\leq 3} \varphi_s \cdot \nabla ( \epsilon_{smn} \partial_m \dt^l v^n) + \dN \dv \dt^l v. 
\end{align*}
In the meantime, from \subeqref{PtNStl}{1}
\begin{equation}\label{NDEst003}
	\begin{aligned}
		& \bar\rho \dN \dv \dt^l v = \dN G_1^{l} - \left( \dt^{l+1}\dN q + (v + \bar u) \cdot \nabla \dt^l \dN q  + \dN (v +\bar u) \cdot \nabla \dt^l q \right. \\& ~~~~~~ \left. + (v+\bar u) \cdot \commut{\dN}{\nabla}\dt^l q + \dN \dt^l v \cdot \nabla \bar\rho + \dt^l v \cdot \dN \nabla \bar\rho + \dN \bar\rho \dv \dt^l v \right).
	\end{aligned}
\end{equation}
Combining \eqref{NDEst001}, \eqref{NDEst002}, \eqref{NDEst003}, we can derive the following equations
\begin{equation}\label{PTN1}
	\begin{cases}
		(2\mu+\lambda) \dt^{l+1} \dN q + (2\mu + \lambda) (v+\bar u)\cdot \nabla\dt^l \dN q + \gamma \bar\rho^2\dN \left(\bar\rho^{\gamma-2}\dt^l q \right) \\~~~~~~~~~~ = F_{1,N}^{l} + G_{1,N}^l,\\
		(2\mu + \lambda) \dN \dv \dt^l v = F_{2,N}^{l} + G_{2,N}^l,
	\end{cases}
\end{equation}
with
\begin{align}\label{Nonlinear-4}
		& F_{1,N}^l = - (2\mu+\lambda) \left(\bar u \cdot\commut{\dN}{\nabla}\dt^l q + \dN \bar u \cdot\nabla\dt^l q \right. {\nonumber}\\
		& ~~~~~~\left. + \dN \dt^l v \cdot \nabla \bar\rho + \dt^l v \cdot \dN \nabla \bar\rho + \dN \bar\rho \dv \dt^l v \right) {\nonumber}\\
		& ~~~~~~ + \mu \bar\rho \sum_{1\leq s,m,n\leq 3} \varphi_s \cdot \nabla ( \epsilon_{smn} \partial_m \dt^l v^n) - \bar\rho \left(\bar\rho \dt^{l+1} v \cdot N - F_{2}^l\cdot N    \right), {\nonumber}\\
		& F_{2,N}^l = - \mu \sum_{1\leq s,m,n\leq 3} \varphi_s \cdot \nabla ( \epsilon_{smn} \partial_m \dt^l v^n) + \left(\bar\rho  \dt^{l+1} v \cdot N \right.{\nonumber}\\
		& ~~~~~~ \left. - F_2^l\cdot N + \gamma\bar\rho\nabla\left( \bar\rho^{\gamma-2} \dt^l q \right) \cdot N  \right), \\
		& G_{1,N}^l = (2\mu+\lambda) \left( \dN G_1^l - (\dN v \cdot \nabla \dt^l q  + v \cdot \commut{\dN}{\nabla}\dt^l q) \right) {\nonumber}\\
		& ~~~~~~ - \bar\rho \left( q \dt^{l+1} v \cdot N + \bar\rho v\cdot \nabla\dt^l v\cdot N  - G_2^l\cdot N \right),   {\nonumber}\\
		& G_{2,N}^l = q \dt^{l+1} v \cdot N + \bar\rho v\cdot \nabla\dt^l v\cdot N - G_2^l\cdot N. {\nonumber}
\end{align}

\begin{lm}[Normal Direction Estimate] The following estimate on $ \dN \dt^l q $ holds,
	\begin{equation}\label{ql01-est}
		\begin{aligned}
			& \dfrac{d}{dt} \left\lbrace \dfrac{2\mu+\lambda}{2} \inti{\abs{\dN \dt^l q}{2} } \right\rbrace + \inti{\abs{ \dN \dt^l q}{2}} \lesssim \Lambda_l \inti{\abs{\nabla\dt^l v}{2}} \\
			& ~ + \Lambda_l \inti{\abs{\nabla \dt^l q}{2}} + \inti{\abs{\dt^{l+1}v}{2}} + \inti{\abs{\nabla \dT \dt^l v}{2}} + \inti{\abs{G_{1,N}^l}{2}} .
		\end{aligned}
	\end{equation}
	Also, $ \dv \dt^l v $ admits the estimate,
	\begin{equation}\label{dvl01-est}
		\begin{aligned}
			& \inti{\abs{\nabla \dv \dt^l v}{2}} + \inti{\abs{\dN \dv \dt^l v}{2}} \lesssim \Lambda_l \inti{\abs{\nabla\dt^l q}{2}} \\
			& ~~ + (1+\Lambda_l) \inti{\abs{\nabla \dt^l v}{2}}  + \inti{\abs{\dt^{l+1}v}{2}} + \inti{\abs{\nabla\dT \dt^l v}{2}}\\
			& ~~ + \inti{\abs{\dN \dt^l q}{2}}  + \inti{\abs{G_{2,N}^{l}}{2}} + \int_{B_{1/2}} \abs{\nabla^2 \dt^l v}{2} \,dx.
		\end{aligned}
	\end{equation}
\end{lm}

\begin{pf}
	After taking inner product of \subeqref{PTN1}{1} with $ \dt^l \dN q $	and recording the resulting after integration by parts, it holds,
	\begin{align*}
			& \dfrac{d}{dt} \left\lbrace \dfrac{2\mu+\lambda}{2} \inti{\abs{\dt^l \dN q}{2} } \right\rbrace + \gamma \inti{\bar\rho^\gamma \abs{ \dt^l \dN q}{2}}   \\ 
			& ~~~~~~ = - (\gamma-2)\inti{\dt^l q \dt^l \dN q \dN \bar\rho^{\gamma}} + \dfrac{2\mu+\lambda}{2} \inti{\dv( v+\bar u) \abs{\dt^l \dN q}{2}}\\
			& ~~~~~~ +\inti{F_{1,N}^l\dt^l \dN q } + \inti{G_{1,N}^l \dt^l \dN q}.
	\end{align*}	
	Meanwhile,
	\begin{align*}
			& \inti{\abs{F_{1,N}^l}{2}} \lesssim \bar\omega^a \left(\inti{\abs{\nabla\dt^l v}{2}} + \inti{\abs{\nabla\dt^l q}{2}} \right) + \inti{\abs{\dt^{l+1}v}{2}}\\
			& ~~~~ + \inti{\abs{\nabla \dT \dt^l v}{2}} + \inti{\abs{F_2^l}{2}} \lesssim \bar\omega^a \left( \inti{\abs{\nabla\dt^l v}{2}} + \inti{\abs{\nabla\dt^l q}{2}} \right) \\
			& ~~~~ + \inti{\abs{\dt^{l+1}v}{2}} + \inti{\abs{\nabla \dT \dt^l v}{2}}, 
	\end{align*}
	for some $ a > 0 $. 
	\eqref{ql01-est} then follows from H\"{o}lder's inequality, \eqref{qPoincare}, \eqref{vPoincare} and the fact $ \inf_{x\in\Omega} \bar\rho > 0 $.
	On the other hand,
	\begin{equation*}
		\inti{\abs{\dN \dv \dt^l v}{2}} \lesssim \inti{\abs{F_{2,N}^{l}}{2}} + \inti{\abs{G_{2,N}^{l}}{2}}.
	\end{equation*}
	Similarly,
	\begin{align*}
			& \inti{\abs{F_{2,N}^{l}}{2}} \lesssim \bar\omega^2  \inti{\abs{\nabla \dt^l v}{2}} + \inti{\abs{\dt^{l+1}v}{2}} + \inti{\abs{\nabla\dT \dt^l v}{2}}\\
			& ~~~~~~ + \inti{\abs{\dN \dt^l q}{2}} + \bar\omega^a \inti{\abs{\nabla \dt^l q}{2}}.
	\end{align*}
	Then together with the fact
	\begin{align*}
		& \inti{\abs{\nabla \dv \dt^l v}{2}} \lesssim \int_{B_{1/2}} \abs{\nabla \dv \dt^l v}{2} \,dx + \inti{\abs{\dT \dv \dt^l v}{2}} \\
		& ~~~~~~ + \inti{\abs{\dN \dv \dt^l v}{2}} \lesssim \int_{B_{1/2}} \abs{\nabla^2 \dt^l v}{2} \,dx + \inti{\abs{ \nabla \dT \dt^l v}{2}} \\
		& ~~~~~~ + \inti{\abs{\nabla \dt^l v
		}{2}} + \inti{\abs{\dN \dv \dt^l v}{2}},
	\end{align*}
	as the consequence of \eqref{DFrame02} and \eqref{commt}, \eqref{dvl01-est} follows after chaining the these inequalities.
\end{pf}

%

Now it is time to introduce the associated Stokes problem. 
From \subeqref{PtNStl}{2}, $ (\dt^l v, \gamma \bar\rho^{\gamma-1}\dt^l q)  $ satisfies the following Stokes system. 

\begin{equation}\label{PtS}
	\begin{cases}
		- \dv \mathbb{S}(\dt^l v) + \nabla(\gamma \bar\rho^{\gamma-1} \dt^l q) = F_s^l + G_s^l & \text{in} ~ \Omega,  \\
		\dv \dt^l v = \dv \dt^l v  & \text{in} ~ \Omega, \\
		\dt^l v = \dt^l v & \text{on} ~ \Gamma,  \\
	\end{cases}
\end{equation}
with
\begin{equation}\label{Nonlinear-5}
	\begin{aligned}
		& F_s^l = F_2^l - \bar\rho\dt^{l+1}v + \gamma \bar\rho^{\gamma-2} \dt^l q \nabla \bar\rho, \\
		& G_s^l = G_2^l - q \dt^{l+1}v - \bar\rho v \cdot \nabla\dt^l v.
	\end{aligned}
\end{equation}

\begin{lm} By applying the Stokes estimate in Lemma \ref{lm:StokesP} to \eqref{PtS}, we shall obtain, 
	\begin{equation}\label{qvl02-est}
		\begin{aligned}
			& \norm{\dt^l v}{\Hnorm{2}}^2 + \norm{\nabla\dt^l q}{\Lnorm}^2 \lesssim \norm{\nabla\dv \dt^l v}{\Lnorm}^2 + \norm{\nabla\dT \dt^l v}{\Lnorm}^2 \\
			& ~~~~~~ + \Lambda_l \inti{\abs{\nabla\dt^l q}{2}} + (1+\Lambda_l) \inti{\abs{\nabla\dt^l v}{2}} \\
			& ~~~~~~ + \inti{\abs{\dt^{l+1}v}{2}} + \norm{G_s^l}{\Lnorm}^2.
		\end{aligned}
	\end{equation}
\end{lm}

\begin{pf}
	By noticing
	\begin{equation}\label{id:viscoustensor}
		\dv \mathbb{S} = \mu \Delta + (\mu + \lambda) \nabla \dv,
	\end{equation}
	the Stokes estimate in Lemma \ref{lm:StokesP} then yields
	\begin{align*}
			& \norm{\dt^l v}{\Hnorm{2}}^2 + \norm{\nabla( \gamma \bar\rho^{\gamma-1} \dt^l q)  }{\Lnorm}^2 \leq \norm{F_s^l}{\Lnorm}^2 + \norm{G_s^l}{\Lnorm}^2 \\
			& ~~~~~~ + \norm{\dv\dt^l v}{\Hnorm{1}}^2 + \norm{\dt^l v}{\bHnorm{3/2}}^2.
	\end{align*}
	Meanwhile, the trace embedding inequality \eqref{trace} implies
	\begin{equation*}
		\norm{\dt^l v}{\bHnorm{3/2}}^2 \lesssim \sum_{j=0}^{1} \norm{\dT^j \dt^l v}{\Hnorm{1}}^2 \lesssim \norm{\nabla \dt^l v}{\Lnorm}^2 + \norm{\nabla \dT \dt^l v}{\Lnorm}^2.
	\end{equation*}
	On the other hand, direct calculation gives the following, 
	\begin{align*}
			& \norm{F_s^l}{L^2}^2 \lesssim \bar\omega^a \inti{\abs{\nabla\dt^l q}{2}} + \bar\omega^a \inti{\abs{\nabla\dt^l v}{2}} + \inti{\abs{\dt^{l+1}v}{2}},\\
			& \norm{\nabla\dt^l q}{\Lnorm}^2 \lesssim \norm{\nabla( \gamma \bar\rho^{\gamma-1} \dt^l q)  }{\Lnorm}^2 + \bar\omega^a \norm{\nabla\dt^l q}{\Lnorm}^2.
	\end{align*}
	Hence \eqref{qvl02-est} holds. 
\end{pf}

The following lemma is a consequence of \eqref{vql1-est}, \eqref{ql01-est}, \eqref{dvl01-est} and \eqref{qvl02-est}.

\begin{prop} We have obtained the estimate concerning the spatial derivatives as follows,
	\begin{align}\label{Proposition-tangential}
			& \dfrac{d}{dt}\biggl\lbrace \inti{(\bar\rho + q)\abs{\dT \dt^l v}{2}} + \inti{\bar\rho^{\gamma-2} \abs{\dT\dt^l q}{2}} + \inti{\abs{\dN\dt^l q}{2}} \biggr\rbrace {\nonumber}\\
			& ~~ + \inti{\abs{\nabla^2 \dt^l v}{2}} + \inti{\abs{\nabla \dv \dt^l v}{2}} + \inti{\abs{\nabla \dt^l q}{2}}{\nonumber}\\
			& \lesssim \Lambda_l \biggl( \inti{\abs{\nabla^2 \dt^l v}{2}} +   \inti{\abs{\nabla\dt^l q}{2}} \biggr) +  (1+\Lambda_l) \biggl( \inti{\abs{\nabla \dt^l v}{2}}  \\ &  + \inti{\abs{\dt^{l+1}v}{2}} \biggr) + \inti{\abs{G_1^{l,1}}{2}} + \inti{\abs{G_2^{l,1}}{2}} + \inti{\abs{G_2^l}{2}} {\nonumber}\\
			&  + \inti{\abs{G_{1,N}^l}{2}} + \inti{\abs{G_{2,N}^l}{2}} + \inti{\abs{G_s^l}{2}} + \int_{B_{1/2}} \abs{\nabla^2 \dt^l v}{2}\,dx.{\nonumber}
	\end{align}

\end{prop}
\begin{pf}
	From \eqref{ql01-est}, \eqref{dvl01-est}, \eqref{qvl02-est}, it holds
	\begin{align*}
			& \dfrac{d}{dt}\inti{\abs{\dN \dt^l q}{2}} + \inti{\abs{\dN \dt^l q}{2}} + \inti{\abs{\nabla\dv\dt^l v}{2}} \\
			& + \inti{\abs{\nabla^2 \dt^l v}{2}} + \inti{\abs{\nabla\dt^l q}{2}} \lesssim \Lambda_l  \inti{\abs{\nabla\dt^l q}{2}}  \\
			& +  (1+\Lambda_l) \left( \inti{\abs{\nabla \dt^l v}{2}} + \inti{\abs{\dt^{l+1}v}{2}} \right) + \inti{\abs{\nabla\dT \dt^l v}{2}}\\
			& + \inti{\abs{G_{1,N}^l}{2}} + \inti{\abs{G_{2,N}^l}{2}} + \inti{\abs{G_s^l}{2}} + \int_{B_{1/2}} \abs{\nabla^2 \dt^l v}{2}\,dx.
	\end{align*}
	Together with \eqref{vql1-est} and an appropriate choice of $ \omega > 0 $, \eqref{Proposition-tangential} follows easily.
\end{pf}

\subsection{On Higher Order Spatial Derivatives}
Through the following arguments similar to those in the last section, the estimates involving higher order spatial derivatives would be shown. In particular, the estimates of tangential derivatives are obtained through a high-order version of \eqref{TPt}. Also, by taking mixed derivatives to \eqref{PTN1}, we shall track the regularities of $ q $ and more importantly, the regularities of $ \dv v $. Then a high-order version of the Stokes problem \eqref{PtS} would eventually yield the estimates of $ v $ and $ q $. These estimates would play important roles in the global analysis. 


Apply $ \dT^m $ to \eqref{PtNStl} and record the resulting system as follows,
\begin{equation}\label{Ptlm}
	\begin{cases}
		\dt^{l+1} \dT^m q + \dv ( \bar\rho \dt^l \dT^m v ) + (v + \bar u)\cdot \nabla\dt^l \dT^m q = F_{1}^{l,m} + G_{1}^{l,m} & \text{in} ~ \Omega,\\
   		(\bar\rho + q ) \dt^{l+1} \dT^m v + \bar\rho v\cdot \nabla\dt^l \dT^m v + \gamma \bar\rho \nabla (\bar\rho^{\gamma-2} \dt^l \dT^m q) \\
		~~~~~~ - \dv \dT^m \mathbb{S}(\dt^l v) = F_{2}^{l,m} + \dT^m G_2^l + G_{2}^{l,m} & \text{in} ~ \Omega,\\
		\dt^l\dT^m v\cdot \vec{n} = b_1^m & \text{on} ~ \Gamma,\\
		\vec{\tau} \cdot \dT^m\mathbb{S}(\dt^l v) \vec{n} = b_2^m & \text{on} ~ \Gamma,
	\end{cases}
\end{equation}
with
\begin{align*}
		& F_{1}^{l,m} = - \biggl\lbrace \bar u \cdot\commut{\dT^m}{\nabla}\dt^l q + \sum_{j=0}^{m-1} C_{j,m}\dT^{m-j} \bar u\cdot \dT^j\nabla\dt^l q \\ & ~~~~~~ + \commut{\dT^m}{\dv} (\bar\rho \dt^l v) + \sum_{j=0}^{m-1} C_{j,m} \dv (\dT^{m-j}\bar\rho \dt^l \dT^j v) \biggr\rbrace, \\
		& F_{2}^{l,m} = \dT^m F_{2}^l - \biggl\lbrace \sum_{j=0}^{m-1} C_{j,m} \dT^{m-j} \bar\rho \dt^{l+1}\dT^j v + \gamma \bar\rho \commut{\dT^m}{\nabla} (\bar\rho^{\gamma-2} \dt^l q) \\& ~~~~~~  + \sum_{j=0}^{m-1} C_{j,m}\lbrace \gamma\bar\rho \nabla(\dT^{m-j}\bar\rho^{\gamma-2}\dT^j\dt^l q) + \gamma \dT^{m-j}\bar\rho \dT^j \nabla(\bar\rho^{\gamma-2} \dt^l q) \rbrace \biggr\rbrace\\ & ~~~~~~  + \commut{\dT^m}{\dv} \mathbb{S} (\dt^l v),   \\
		& G_{1}^{l,m} = \dT^m G_1^l - \biggl\lbrace v\cdot \commut{\dT^m}{\nabla}\dt^l q + \sum_{j=0}^{m-1} C_{j,m} \dT^{m-j} v \cdot \dT^j\nabla\dt^l q \biggr\rbrace , \\
		& G_{2}^{l,m} = - \biggl\lbrace \sum_{j=0}^{m-1} C_{j,m} \dT^{m-j} q \dt^{l+1}\dT^j v + \bar\rho v \cdot \commut{\dT^m}{\nabla} \dt^l v \\& ~~~~~~ + \sum_{j=0}^{m-1} C_{j,m} \dT^{m-j}(\bar\rho v) \cdot \dT^j \nabla\dt^l v   \biggr\rbrace ,\\
\end{align*}
and
\begin{align*}
		& b^m_1 = - \sum_{j=0}^{m-1} C_{j,m}\dt^l \dT^j v\cdot \dT^{m-j}\vec{n} ,\\
		& b^m_2 = - \sum_{j=0}^{m-1} \sum_{a+b = m-j} C_{j,m}C_{a,m-j}\dT^a\vec{\tau} \cdot \dT^j\mathbb{S}(\dt^l v) \dT^b \vec{n}.\\
\end{align*}
Applying standard energy estimate to \eqref{Ptlm} then will yield the following lemma.
\begin{lm} For $ m \geq 1 $, the following higher version of \eqref{vql1-est} holds, 
\begin{equation}\label{vqlm-est}
	\begin{aligned}
		& \dfrac{d}{dt} \left\lbrace  \dfrac{1}{2} \inti{(\bar\rho+q) \abs{\dt^l \dT^m v}{2}} + \dfrac{\gamma}{2} \inti{\bar\rho^{\gamma-2}\abs{\dt^l \dT^m q}{2}} \right\rbrace \\
		& ~~~~ + \dfrac{\mu}{4} \inti{\abs{\nabla\dT^m \dt^l v}{2}} + \lambda\inti{\abs{\dv(\dT^m \dt^l v)}{2}} \\
		& \lesssim (\omega + \Lambda_l) \inti{\abs{\nabla^{m+1}\dt^l v}{2}} + (\omega  + C_{\omega} \bar\omega^a + \Lambda_l  ) \inti{\abs{\nabla^m \dt^l q}{2}} \\
		& + (\omega + C_\omega + \Lambda_l ) \biggl( \sum_{j=1}^{m} \inti{\abs{\nabla^j\dt^l v}{2}} + \sum_{j=1}^{m-1} \inti{\abs{\nabla^j\dt^l q}{2}} \biggr)\\
		& + \omega\bar\omega^a \sum_{j=0}^{m-1} \inti{\abs{\nabla^j\dt^{l+1} v}{2}} + C_\omega\inti{\abs{G_1^{l,m}}{2}} + C_\omega \inti{\abs{G_{2}^{l,m}}{2}} \\
		& + \inti{\abs{\dT^{m-1}G_2^l}{2}}.
	\end{aligned}
\end{equation}
with any $ 0 < \omega < 1 $.  
\end{lm}

\begin{pf}
Taking inner product of \subeqref{Ptlm}{2} with $ \dt^l \dT^m v $ and recording the resulting after integration by parts,
\begin{align*}
	& \dfrac{d}{dt} \left\lbrace \dfrac{1}{2} \inti{(\bar\rho+q) \abs{\dt^l \dT^m v}{2}} \right\rbrace \underbrace{ - \gamma \inti{ \bar\rho^{\gamma-2} \dt^l \dT^m q \dv(\bar\rho \dt^l \dT^m v) }  }_{(i)} \\
	&  + \dfrac{\mu}{2} \inti{\abs{\nabla\dt^l \dT^m v + \nabla \dt^l \dT^m v^T}{2}} + \lambda \inti{\dv(\dt^l \dT^m v)} \\
	&  + \underbrace{ \inti{\commut{\dT^m}{\mathbb{S}}(\dt^l v) : \nabla \dt^l \dT^m v } }_{(ii)} = \underbrace{\inti{(F_2^{l,m}+G_2^{l,m}) \cdot \dt^l \dT^m v} }_{(iii)}\\
	& + \underbrace{\inti{\dT^mG_2^l \cdot \dt^l \dT^m v}}_{(iv)} +\underbrace{ \dfrac{1}{2} \inti{\dt q \abs{\dt^l\dT^m v}{2}} + \dfrac{1}{2} \inti{\dv(\bar\rho v) \abs{\dt^l \dT^m v}{2}} }_{(v)}\\
	& - \gamma\underbrace{\intb{\bar\rho^{\gamma-1}\dt^l \dT^m q \dt^l \dT^m v \cdot \vec{n}} }_{(vi)} + \underbrace{ \intb{ \dt^l \dT^m v \cdot \dT^m\mathbb{S}(\dt^l v) \vec{n}} }_{(vii)}.
\end{align*}
From \subeqref{Ptlm}{1}, the following identity holds,
\begin{equation*}
	\begin{aligned}
		& (i) = \gamma \inti{\bar\rho^{\gamma-2} \dt^l \dT^m q \cdot (\dt^{l+1} \dT^m q + (v + \bar u) \cdot \nabla\dt^l\dT^m q - F_{1}^{l,m} -G_1^{l,m} )}\\
			& ~~~~ = \dfrac{d}{dt} \dfrac{\gamma}{2} \inti{\bar\rho^{\gamma-2} \abs{\dt^l\dT^m q}{2}} - \dfrac{\gamma}{2} \underbrace{\inti{ \dv(\bar\rho^{\gamma-2}(v+\bar u))\abs{\dt^l \dT^m q}{2} }}_{(viii)} \\
		& ~~~~~~ - \gamma\underbrace{\inti{(F_1^{l,m}+G_1^{l,m}) \cdot \bar\rho^{\gamma-2} \dt^l \dT^m q}}_{(ix)} .
	\end{aligned}
\end{equation*}
By applying Cauchy's inequality, Poincar\'{e} inequality and \eqref{FTC2}, 
\begin{align*}
		& (viii) \lesssim ( \bar\omega^a + \bar\omega\norm{v}{\supnorm} + \norm{\nabla v}{\supnorm}) \inti{\abs{\dt^l \dT^m q}{2}}, \\
		& (ix) \lesssim \omega \inti{\abs{\dT^m\dt^l q}{2}} + C_{\omega}\biggl( \inti{\abs{F_{1}^{l,m}}{2}} + \inti{\abs{G_{1}^{l,m}}{2}} \biggr), \\
		& (ii) \lesssim \delta \inti{\abs{\nabla \dT^m\dt^l v}{2}} + C_\delta \sum_{j=0}^{m}\inti{\abs{\nabla^j \dt^l v }{2}}, \\
		& (iii) \lesssim \omega \inti{\abs{F_2^{l,m}}{2}} + C_{\omega} \biggl( \inti{\abs{G_2^{l,m}}{2}} + \inti{\abs{\dt^l \dT^m v}{2}} \biggr),\\
		& (iv) = - \inti{\dT^{m-1}G_2^l\cdot\dT\dt^l\dT^m v} \lesssim \delta \inti{\abs{\nabla\dT^m \dt^l v}{2}}  \\ & ~~~~~~ + C_\delta \inti{\abs{\dT^{m-1}G_2^l}{2}},\\
		& (v) \lesssim (\norm{\dt q}{\supnorm} + \bar\omega\norm{v}{\supnorm} + \norm{\nabla v}{\supnorm}) \inti{\abs{\dt^l \dT^m v}{2}}.
\end{align*}
Also, followed from the definition of $F_1^{l,m}, F_2^{l,m}$, it holds
\begin{align*}
		& \inti{\abs{F_1^{l,m}}{2}} \lesssim (1+\bar\omega^a) \sum_{j=0}^{m}\inti{\abs{\nabla^j \dt^l v}{2}} + \bar\omega^a \sum_{j=0}^{m} \inti{\abs{\nabla^j\dt^l q}{2}}, \\
		& \inti{\abs{F_2^{l,m}}{2}} \lesssim (1+\bar\omega^a) \inti{\abs{\nabla^{m+1}\dt^l v}{2}} + (1+\bar\omega^a)\inti{\abs{\nabla^m\dt^l q}{2}}\\
		& ~~~~~~ + \bar\omega^a \sum_{j=0}^{m-1} \biggl(\inti{\abs{\nabla^j\dt^l q}{2}} +\inti{\abs{\nabla^{j}\dt^{l+1} v}{2}}\biggr) \\
		& ~~~~~~ + (1+\bar\omega^a) \sum_{j=0}^{m}\inti{\abs{\nabla^j\dt^l v}{2}}.
\end{align*}
Meanwhile, to handle the boundary integration $ (vi) $, from \subeqref{Ptlm}{3} and the calculus on the boundary \eqref{FTCB2},
\begin{equation*}
	\begin{aligned}
		& (vi) = \intb{ \bar\rho^{\gamma-1}\dt^l\dT^m q b_1^m } = - \intb{\dt^l \dT^{m-1} q \dT ( \bar\rho^{\gamma-1} b_1^m )}\\
		& ~~~~ \lesssim \norm{\dt^l \dT^{m-1} q }{\bLnorm} \left(\bar\omega^a \norm{b_1^m}{\bLnorm} + \norm{\dT b_1^m}{\bLnorm} \right) \\
		& ~~~~ \lesssim (1+\bar\omega^a) \norm{\dt^l \dT^{m-1} q }{\bLnorm}  \times\sum_{j=0}^{m} \norm{\dt^l \dT^j v}{\bLnorm},
	\end{aligned}
\end{equation*}
where it has been used from the definition of $ b_1^m $,
\begin{equation*}
	\begin{aligned}
		& \abs{b_1^m}{} \lesssim \sum_{j=0}^{m-1} \abs{\dt^l \dT^j v}{}, ~~ \abs{\dT b_1^m}{} \lesssim \sum_{j=0}^{m} \abs{\dt^l \dT^j v}{}.
	\end{aligned}
\end{equation*}
Then the trace embedding inequality \eqref{trace1} together with \eqref{qPoincare} and \eqref{vPoincare} implies
\begin{equation*}
	\begin{aligned}
		& (vi) \lesssim \omega \left( \inti{\abs{\nabla^m\dt^l q}{2}} +  \inti{\abs{\nabla^{m+1}\dt^l v}{2}} \right) \\
		& ~~~~~~ + C_\omega (1+\bar\omega^a) \biggl( \inti{\abs{ \nabla^{m-1} \dt^l q}{2}} +  \sum_{j=0}^{m} \inti{\abs{\nabla^j \dt^l v}{2}} \biggr).
	\end{aligned}
\end{equation*}
On the other hand, \subeqref{Ptlm}{3} implies the following decomposition on $ \Gamma $,
\begin{equation*}
		\dt^l \dT^m v = b_1^m \vec{n} + v_{l,m} \vec{\tau}
\end{equation*}
with $ \abs{v_{l,m}}{} \leq \abs{\dt^l\dT^m v}{} $. Then together with \subeqref{Ptlm}{4} 
\begin{equation*}
	\begin{aligned}
		& (vii) = \intb{b_1^m \vec{n} \cdot \dT^m\mathbb{S}(\dt^l v) \vec{n} }+ \intb{v_{l,m}\vec{\tau} \cdot \dT^m\mathbb{S}(\dt^l v) \vec{n} }\\
		& ~~~~ = - \intb{\dT ( b_1^m \vec{n}) \cdot \dT^{m-1}\mathbb{S}(\dt^l v) \vec{n} } - \intb{b_1^m \vec{n} \cdot \dT^{m-1} \mathbb{S}(\dt^l v) \dT\vec{n}}\\
		& ~~~~~~ + \intb{v_{l,m} \cdot b_2^m} \lesssim \norm{\dT^{m-1}\mathbb S(\dt^l v)}{\bLnorm}\left( \norm{b_1^m}{\bLnorm} + \norm{\dT b_1^m}{\bLnorm} \right)\\
		& ~~~~~~ + \norm{\dt^l \dT^m v}{\bLnorm} \norm{b_2^m}{\bLnorm}.
	\end{aligned}
\end{equation*}
From the definition,
\begin{equation*}
	\abs{b_2^m}{} \lesssim \sum_{j=0}^{m-1}\abs{\dT^j \mathbb S(\dt^l v) }{}.
\end{equation*}
Therefore, by applying trace embedding inequality \eqref{trace1}, Cauchy's inequality and Poincar\'{e} inequality,
\begin{equation*}
	(vii) \lesssim \omega \inti{\abs{\nabla^{m+1} \dt^l v}{2}} + C_{\omega} \sum_{j = 0}^{m}\inti{\abs{\nabla^{j}  \dt^l v}{2}}.
\end{equation*}
Summing up these estimates with an appropriately small $\delta>0 $ and \eqref{korn01}, it holds,
\begin{equation*}
	\begin{aligned}
		& \dfrac{d}{dt} \left\lbrace  \dfrac{1}{2} \inti{(\bar\rho+q) \abs{\dt^l \dT^m v}{2}} + \dfrac{\gamma}{2} \inti{\bar\rho^{\gamma-2}\abs{\dt^l \dT^m q}{2}} \right\rbrace \\
		& ~~~~ + \dfrac{\mu}{4} \inti{\abs{\nabla\dT^m \dt^l v}{2}} + \lambda\inti{\abs{\dv(\dT^m \dt^l v)}{2}} \\
		& \lesssim (\omega + \Lambda_l) \inti{\abs{\nabla^{m+1}\dt^l v}{2}} + (\omega + C_\omega \bar\omega^a+ \Lambda_l  ) \inti{\abs{\nabla^m \dt^l q}{2}} \\
		& + (\omega + C_\omega + \Lambda_l ) \left( \sum_{j=0}^{m} \inti{\abs{\nabla^j\dt^l v}{2}} + \sum_{j=0}^{m-1} \inti{\abs{\nabla^j\dt^l q}{2}} \right)\\
		& + \omega\bar\omega^a \sum_{j=0}^{m-1} \inti{\abs{\nabla^j\dt^{l+1} v}{2}} + C_\omega\inti{\abs{G_1^{l,m}}{2}} + C_{\omega} \inti{\abs{G_{2}^{l,m}}{2}} \\
		& + \inti{\abs{\dT^{m-1}G_2^l}{2}}.
	\end{aligned}
\end{equation*}

Therefore, \eqref{vqlm-est} follows after applying \eqref{qPoincare} and \eqref{vPoincare}.

\end{pf}

Next it is to derive the estimates of mixed derivatives of $ q $. 
Applying $ \dT^m \dN^{n-1} $ to \eqref{PTN1}, it holds
\begin{equation}\label{Ptlmn}
	\begin{cases}
		(2\mu+\lambda) \dt^{l+1}\dT^m \dN^n q + (2\mu+\lambda) (v+\bar u)\cdot \nabla \dt^l \dT^m\dN^nq + \gamma \bar\rho^{\gamma} \dt^l \dT^m \dN^n q\\
		~~~~~~~~~~~~~~~~~~~~~~~~~~~~  = F_{1,N}^{l,m,n} + G_{1,N}^{l,m,n}, \\
		(2\mu+\lambda)\dT^m\dN^n \dv \dt^l v = F_{2,N}^{l,m,n} + G_{2,N}^{l,m,n},
	\end{cases}
\end{equation}
with
\begin{align*}
		& F_{1,N}^{l,m,n} = \dT^m\dN^{n-1} F_{1,N}^l - (2\mu+\lambda)  \bar u \cdot\commut{\dT^m\dN^{n-1}}{\nabla}\dt^l\dN q  \\
		& ~~~~~~ - (2\mu+\lambda)\sum_{j=0}^{m+n-2}\sum_{a+b=j}C_{a,m}C_{b,n-1}\dT^{m-a}\dN^{n-1-b} \bar u\cdot \dT^{a}\dN^b\nabla\dt^l \dN q  \\
		& ~~~~~~ - (\gamma-2)\sum_{j=0}^{m+n-1} \sum_{a+b = j} C_{a,m}C_{b,n-1} \dT^{m-a}\dN^{n-b} \bar\rho^\gamma \dt^l \dT^a \dN^b q \\ & ~~~~~~ - \gamma \sum_{j=0}^{m+n-2} \sum_{a+b=j} C_{a,m}C_{b,n-1} \dT^{m-a} \dN^{n-1-b}\bar\rho^\gamma \dt^l \dT^{a}\dN^{b+1} q, \\
		& F_{2,N}^{l,m,n} = \dT^m\dN^{n-1} F_{2,N}^l, \\
		& G_{1,N}^{l,m,n} = \dT^m\dN^{n-1} G_{1,N}^l - (2\mu+\lambda) v \cdot \commut{\dT^m\dN^{n-1}}{\nabla} \dt^l \dN q  \\&  ~~~~~~ - (2\mu+\lambda) \sum_{j=0}^{m+n-2} \sum_{a+b = j} C_{a,m} C_{b,n-1} \dT^{m-a}\dN^{n-1-b} v \cdot \dT^a \dN^b \nabla \dt^l \dN q, \\
		& G_{2,N}^{l,m,n} = \dT^m\dN^{n-1} G_{2,N}^l. \\
\end{align*}
Then the following estimate holds,
\begin{lm} For $ n \geq 1 , m+n\geq 2 $, 
\begin{equation}\label{qlmn-est}
	\begin{aligned}
		& \dfrac{d}{dt} \biggl\lbrace \dfrac{2\mu+\lambda}{2} \inti{\abs{\dT^m \dN^n \dt^l q}{2}}  \biggr\rbrace + \gamma \inti{\bar\rho^\gamma \abs{\dT^m\dN^n\dt^l q}{2}} \\
		& \lesssim \inti{\abs{\nabla \dT^{m+1} \dN^{n-1} \dt^l v}{2}}  + \inti{\abs{\dT^m\dN^{n-1}\dt^{l+1} v}{2}} \\
		& ~~~~  + ( 1 +\Lambda_l)\sum_{j=1}^{m+n} \inti{\abs{\nabla^j \dt^l v}{2}} + \Lambda_l \sum_{j=1}^{m+n} \inti{\abs{\nabla^{j}\dt^l q}{2}}\\
		& ~~~~  + (1+\Lambda_l) \sum_{j=0}^{m+n-2}\inti{\abs{\nabla^{j}\dt^{l+1} v}{2}}  + \inti{\abs{G_{1,N}^{l,m,n}}{2}}.
	\end{aligned}
\end{equation}
Also,
\begin{align}
		& \inti{\abs{\dT^m\dN^n\dv \dt^l v}{2}} \lesssim \inti{\abs{\dT^m\dN^n\dt^l q}{2}} + \inti{\abs{\nabla \dT^{m+1} \dN^{n-1} \dt^l v}{2}} {\nonumber} \\ & ~~~~~~ + \inti{\abs{\dT^m\dN^{n-1}\dt^{l+1} v}{2}} + ( 1 + \Lambda_l)\sum_{j=1}^{m+n} \inti{\abs{\nabla^j \dt^l v}{2}} { \label{dvlmn-est}} \\
		& ~~~~~~ + \Lambda_l \sum_{j=1}^{m+n-1} \inti{\abs{\nabla^{j}\dt^l q}{2}} + (1+\Lambda_l) \sum_{j=0}^{m+n-2}\inti{\abs{\nabla^{j}\dt^{l+1} v}{2}}{\nonumber}\\
		& ~~~~~~ + \inti{\abs{G_{2,N}^{l,m,n}}{2}}. {\nonumber} 
\end{align}

\end{lm}

\begin{pf}
	After taking inner product of \subeqref{Ptlmn}{1}	with $ \dt^l \dT^m\dN^n q $ and recording the resulting after integration, it admits, 
	\begin{equation}
		\begin{aligned}
			& \dfrac{d}{dt} \biggl\lbrace \dfrac{2\mu+\lambda}{2} \inti{\abs{\dT^m \dN^n \dt^l q}{2}}  \biggr\rbrace + \gamma \inti{\bar\rho^\gamma \abs{\dT^m\dN^n\dt^l q}{2}} \\
			& = \dfrac{2\mu+\lambda}{2}\inti{\dv (v+\bar u) \abs{\dT^m \dN^n \dt^l q}{2}} + \inti{(F_{1,N}^{l,m,n} + G_{1,N}^{l,m,n}) \dT^m\dN^n\dt^l q}\\
			& \lesssim (\delta + \bar\omega + \norm{\nabla v}{\supnorm} ) \inti{\abs{\dT^m\dN^n\dt^l q}{2}} + C_{\delta} \inti{\abs{F_{1,N}^{l,m,n}}{2} + \abs{G_{1,N}^{l,m,n}}{2}}.
		\end{aligned}
	\end{equation}
	Moreover repeatedly applying the formula \eqref{commt} would imply
	\begin{equation}
		\begin{aligned}
			& \inti{\abs{F_{1,N}^{l,m,n}}{2}} \lesssim \inti{\abs{\dT^m\dN^{n-1}F_{1,N}^l}{2}} + \bar\omega^a \sum_{j=0}^{m+n} \inti{\abs{\nabla^{j}\dt^l q}{2}} \\
			& ~~~~ \lesssim \inti{\abs{\nabla \dT^{m+1} \dN^{n-1} \dt^l v}{2}} + \inti{\abs{\dT^m\dN^{n-1}\dt^{l+1} v}{2}} \\
			& ~~~~~~ + ( 1 + \bar\omega^a)\sum_{j=0}^{m+n} \inti{\abs{\nabla^j \dt^l v}{2}} + \bar\omega^a \sum_{j=0}^{m+n} \inti{\abs{\nabla^{j}\dt^l q}{2}} \\
			& ~~~~~~ + (1+\bar\omega^a) \sum_{j=0}^{m+n-2}\inti{\abs{\nabla^{j}\dt^{l+1} v}{2}},
		\end{aligned}
	\end{equation}
	from which \eqref{qlmn-est} follows together with \eqref{qPoincare}, \eqref{vPoincare} and an appropriate $ \delta > 0 $. On the other hand, \eqref{dvlmn-est} is the consequence of \subeqref{Ptlmn}{2}, \eqref{qPoincare}, \eqref{vPoincare} and
	\begin{equation}
		\begin{aligned}
			& \inti{\abs{F_{2,N}^{l,m,n}}{2}} \lesssim \inti{\abs{\dT^m\dN^n\dt^l q}{2}} + \inti{\abs{\nabla \dT^{m+1} \dN^{n-1} \dt^l v}{2}} \\ & ~~~~~~ + \inti{\abs{\dT^m\dN^{n-1}\dt^{l+1} v}{2}} + ( 1 + \bar\omega^a)\sum_{j=0}^{m+n} \inti{\abs{\nabla^j \dt^l v}{2}} \\
			& ~~~~~~ + \bar\omega^a \sum_{j=0}^{m+n-1} \inti{\abs{\nabla^{j}\dt^l q}{2}} + (1+\bar\omega^a) \sum_{j=0}^{m+n-2}\inti{\abs{\nabla^{j}\dt^{l+1} v}{2}}.
		\end{aligned}
	\end{equation}
\end{pf}

The last estimate we shall obtain in this section is from the following Stokes problem. 
After applying $ \dT^m $ to \eqref{PtS}, record the resulting system as follows.
\begin{equation}\label{PtSm}
	\begin{cases}
		- \dv \mathbb{S} (\dt^l \dT^m v) + \nabla (\gamma\bar\rho^{\gamma-1} \dt^l \dT^m q) = F_{s}^{l,m} + G_{s}^{l,m} & \text{in} ~ \Omega,\\
		\dv \dt^l\dT^m v = \dv \dt^l \dT^m v & \text{in} ~ \Omega, \\
		\dt^l \dT^m v = \dt^l \dT^m v & \text{on} ~ \Gamma,
	\end{cases}
\end{equation}
with
\begin{equation}
	\begin{aligned}
		& F_{s}^{l,m} = \dT^m F_{s}^{l} + \commut{\dT^m}{\dv\mathbb{S}}(\dt^l v) - \gamma\commut{\dT^m}{\nabla}(\bar\rho^{\gamma-1} \dt^l q) \\
		& ~~~~~~ - \gamma\sum_{j=0}^{m-1}C_{j,m}\nabla (\dT^{m-j} \bar\rho^{\gamma-1}\dt^l \dT^j q), \\
		& G_{s}^{l,m} = \dT^m G_{s}^{l}.
	\end{aligned}
\end{equation}
As consequences of the Stokes estimate in Lemma \ref{lm:StokesP}, the following lemma indicates the estimates of normal derivatives.
\begin{lm} For $ n \geq 2, m+n > 2 $,
	\begin{align}
		& \norm{\nabla^n \dT^m \dt^l v}{\Lnorm}^2 + \norm{\nabla^{n-1}\dT^m \dt^l q}{\Lnorm}^2 \lesssim \norm{\dT^m\nabla^{n-1} \dv \dt^l v}{\Lnorm}^2 {\nonumber}\\
		& ~~~~ + \norm{\nabla\dT^{n+m-1}\dt^l v}{\Lnorm}^2 + \norm{\nabla^{m+n-2}\dt^{l+1}v}{\Lnorm}^2 {\nonumber}\\
		& ~~~~ + (1+\Lambda_l) \biggl\lbrace \sum_{j=0}^{n+m-3} \norm{\nabla^{j}\dt^{l+1}v}{\Lnorm}^2 + \sum_{j=1}^{n+m-1} \norm{\nabla^j \dt^l v}{\Lnorm}^2  {\label{qvlmn-est}}\\
		& ~~~~ +   \sum_{j=1}^{n+m-2}\norm{\nabla^j \dt^l q}{\Lnorm}^2 \biggr\rbrace + \norm{G_s^{l,m}}{\Hnorm{n-2}}^2. {\nonumber}
	\end{align}
\end{lm}

\begin{pf}
Similarly as before, by noticing \eqref{id:viscoustensor}, the Stokes estimate in Lemma \ref{lm:StokesP} then yields,
\begin{align*}
	& \norm{\dt^l \dT^m v}{\Hnorm{n}}^2 + \norm{\nabla(\bar\rho^{\gamma-1} \dt^l \dT^m q)}{\Hnorm{n-2}}^2 \leq \norm{F_{s}^{l,m}}{\Hnorm{n-2}}^2 \\
	& ~~~~~~ + \norm{G_{s}^{l,m}}{\Hnorm{n-2}}^2 + \norm{\dv \dt^l \dT^m v}{\Hnorm{n-1}}^2 + \norm{\dt^l \dT^m v}{\bHnorm{n-1/2}}^2.
\end{align*}
On one hand, triangle inequality implies
\begin{align*}
		& \norm{\nabla^{n-1}\dT^m \dt^l q}{\Lnorm}^2 \lesssim \norm{\nabla(\bar\rho^{\gamma-1} \dt^l \dT^m q)}{\Hnorm{n-2}}^2\\
		& ~~~~~~~~ + \bar\omega^a \sum_{j=0}^{n-2} \norm{\nabla^j\dT^m\dt^l q}{\Lnorm}^2, \\
		& \norm{\nabla^n \dT^m \dt^l v}{\Lnorm}^2 \lesssim \norm{\dt^l \dT^m v}{\Hnorm{n}}^2 + \sum_{j=0}^{n-1}\norm{\nabla^j\dT^m\dt^l v}{\Lnorm}^2.
\end{align*}
Again, the trace embedding inequality \eqref{trace}, together with \eqref{vPoincare}, yields
\begin{align*}
		& \norm{\dt^l \dT^m v}{\bHnorm{n-1/2}}^2 \lesssim \sum_{j=0}^{n-1}\norm{\dt^l \dT^{m+j} v}{\Hnorm{1}}^2 \lesssim \sum_{j=0}^{n-1} \norm{\nabla \dT^{m+j} \dt^l v}{\Lnorm}^2.
\end{align*}
Meanwhile, from the definition of $ F_s^{l,m} $, it admits
\begin{align*}
		& \norm{F_s^{l,m}}{\Hnorm{n-2}}^2 \lesssim \norm{\dT^m F_s^l}{\Hnorm{n-2}}^2 + \sum_{j=0}^{m+1} \norm{\nabla^j \dt^l v}{\Hnorm{n-2}}^2 \\
		& ~~~~ + (1+\bar\omega^a) \norm{\nabla^{m}\dt^l q}{\Hnorm{n-2}}^2 + \bar\omega^a \sum_{j=0}^{m-1}\norm{\nabla^j \dt^l q}{\Hnorm{n-2}}^2 \\
		& \lesssim  \norm{\nabla^{m + n-2}\dt^{l+1}v}{\Lnorm}^2 + (1+\bar\omega^a) \biggl\lbrace\sum_{j=0}^{n+m-3} \norm{\nabla^{j}\dt^{l+1}v}{\Lnorm}^2 \\
		& ~~~~  + \sum_{j=0}^{n+m-1}\norm{\nabla^j\dt^l v}{\Lnorm}^2 +  \sum_{j=0}^{n+m-2}\norm{\nabla^j\dt^l q }{\Lnorm}^2 \biggr\rbrace,
\end{align*}
where we have employed the following estimate of $ F_s^l $, 
\begin{align*}
		& \norm{\dT^mF_s^l}{\Hnorm{n-2}}^2 \lesssim \norm{\nabla^m \dt^{l+1}v}{\Hnorm{n-2}}^2  + \bar\omega^a \biggl( \sum_{j=0}^{m-1} \norm{\nabla^j\dt^{l+1} v}{\Hnorm{n-2}}^2 \\
		& ~~~~  + \sum_{j=0}^{m+1} \norm{\nabla^j\dt^l v}{\Hnorm{n-2}}^2 + \sum_{j=0}^m \norm{\nabla^j\dt^lq}{\Hnorm{n-2}}^2 \biggr).
\end{align*}
Moreover, the commutator formula \eqref{commt} implies
\begin{equation}
	\begin{aligned}
		& \norm{\dv \dt^l \dT^m v}{\Hnorm{n-1}}^2 \lesssim \norm{\dT^m\nabla^{n-1} \dv \dt^l v}{\Lnorm}^2 + \sum_{j=0}^{n+m-1} \norm{\nabla^j \dt^l v}{\Lnorm}^2.
	\end{aligned}
\end{equation}

Therefore, after summing up these estimates, together with \eqref{qPoincare}, \eqref{vPoincare}
\eqref{qvlmn-est} follows.
\end{pf}

\subsection{Interior Estimates}

We will need one last block to establish the stability theory. To track the propagation of regularity away from the boundary, we shall write down the corresponding system in the interior subdomain. Let $ \tilde v = \psi \cdot v, ~ \tilde q = \psi \cdot q $. Then from \eqref{PtNS}, the system satisfied by $ \tilde v, \tilde q $ admits the form
\begin{equation}\label{iPtNS}
	\begin{cases}
		\dt \iq + \dv(\bar\rho \iv) + (v+\bar u)\cdot \nabla \iq = \tilde F_1 + \tilde G_1 & \text{in} ~ \Omega,\\
		(\bar\rho + \rho) \dt \iv + \bar\rho v \cdot \nabla \iv + \gamma \bar\rho \nabla(\bar\rho^{\gamma-2}\iq) - \dv \mathbb S (\iv) = \tilde F_2 + \tilde G_2 & \text{in} ~ \Omega, \\
		(\iq, \iv) = (0,0) & \text{in} ~ \Omega\backslash B_{1/2},
	\end{cases}
\end{equation}
where
\begin{align*}
		& \tilde F_1 = \bar\rho v \cdot \nabla \psi + q \bar u \cdot \nabla \psi, \\
		& \tilde F_2 = \psi F_2 + \gamma \bar\rho^{\gamma-1}q \nabla \psi - (\mu+\lambda)\dv v \nabla \psi - 2\mu (\nabla \psi \cdot \nabla) v \\
		& ~~~~~~ - \mu \Delta \psi v - (\mu+\lambda) \nabla ( v\cdot \nabla \psi), \\
		& \tilde G_1 = \psi G_1 + q v \cdot \nabla \psi , \\
		& \tilde G_2 = \psi G_2 + \bar\rho ( v \cdot \nabla \psi) v.
\end{align*}

After applying the differential operator $ \partial^m \dt^l $ to \eqref{iPtNS}, record the resulting system in the following, where $ \partial $ represents the spatial derivative in any direction (i.e., $ \partial = \partial_i $ for $ i = 1,2,3 $),
\begin{equation}\label{iPtlm}
	\begin{cases}
		\dt \partial^m \dt^l \iq + \dv (\bar\rho \partial^m\dt^l \iv) + (v+\bar u) \cdot \nabla\partial^m\dt^l \iq = \tilde F_1^{l,m} + \tilde G_1^{l,m} & \text{in} ~ \Omega, \\
		(\bar\rho + q ) \dt \partial^m \dt^l \iv + \bar\rho v\cdot \nabla\partial^m\dt^l \iv + \gamma\bar\rho \nabla ( \bar\rho^{\gamma-2}\partial^m\dt^l \iq ) \\
		~~~~~~ - \dv \mathbb S (\partial^m\dt^l \iv) = \tilde F_{2}^{l,m} + \partial^m\dt^l\tilde G_2 + \tilde G_{2}^{l,m} & \text{in} ~ \Omega, \\
		( \partial^m\dt^l\iq,\partial^m\dt^l\iv ) = (0,0) & \text{in} ~ \Omega \backslash B_{3/4},
	\end{cases}
\end{equation}
where
\begin{align*}
		& \tilde F_{1}^{l,m} = \partial^m \dt^l \tilde F_1 - \commut{\partial^m}{\dv}(\bar\rho \dt^l \iv) - \sum_{j=0}^{m-1} C_{j,m} \dv (\partial^{m-j}\bar\rho \partial^{j} \dt^l \iv )  \\
		& ~~~~ - \sum_{j=0}^{m-1} C_{j,m}\partial^{m-j}\bar u \cdot\nabla\partial^j\dt^l \iq ,\\
		& \tilde F_{2}^{l,m} = \partial^m \dt^l \tilde F_2 - \sum_{j=0}^{m-1} C_{j,m} \partial^{m-j} \bar\rho \partial^j\dt^{l+1} \iv\\
		& ~~~~ - \gamma \sum_{j=0}^{m-1} \sum_{a=0}^{m-j} C_{j,m} C_{a,m-j}\partial^a \bar\rho \nabla ( \partial^{m-j-a}\bar\rho^{\gamma-2} \partial^j \dt^l \iq ) , \\
		& \tilde G_{1}^{l,m} = \partial^m \dt^l \tilde G_1 - \sum_{j=0}^{m+l-1} \sum_{a+b=j} C_{a,m}C_{b,l} \partial^{m-a}\dt^{l-b} v \cdot \nabla \partial^{a}\dt^b \iq , \\
		& \tilde G_{2}^{l,m} = - \sum_{j=0}^{m+l-1} \sum_{a+b=j} C_{a,m}C_{b,l} \partial^{m-a}\dt^{l-b} q \partial^a\dt^{b+1} \iv \\
		& ~~~~ - \sum_{j=0}^{m+l-1} \sum_{a+b=j} \sum_{c=0}^{m-a} C_{a,m}C_{b,l} C_{c,m-a} \partial^c\bar\rho \partial^{m-a-c}\dt^{l-b}v\cdot\nabla\partial^a\dt^b\iv.
\end{align*}

Next lemma is concerning the regularity of $ (q,v) $ in the interior sub-domain.

\begin{lm} For $ m \geq 1 $,
	\begin{equation}\label{ivplm-est}
	\begin{aligned}
		& \dfrac{d}{dt}\biggl\lbrace \dfrac{1}{2}\inti{(\bar\rho + q) \abs{\nabla^m \dt^l \iv}{2}} + \dfrac{\gamma}{2} \inti{\bar\rho^{\gamma-2}\abs{\nabla^m\dt^l \iq }{2}} \biggr\rbrace \\
		& ~~~~~~ + \mu \inti{\abs{ \nabla^{m+1} \dt^l \iv}{2}}\\
		& \lesssim (\omega + C_{\omega} \bar\omega^a + \Lambda_l) \inti{\abs{\nabla^m\dt^l q}{2}}  + (\omega + \Lambda_l) \inti{\abs{\nabla^{m+1}\dt^l v}{2}}\\
		& ~~ + (C_\omega  + \Lambda_l) \biggl( \sum_{j=1}^{m}\inti{\abs{\nabla^j \dt^l v}{2}} + \sum_{j=1}^{m-1}\inti{\abs{\nabla^j\dt^l q}{2}} \biggr)\\& ~~  + \omega \Lambda_l \sum_{j=0}^{m-1} \inti{\abs{\nabla^j\dt^{l+1} v}{2}}\\
		& ~~ + C_\omega \inti{\abs{\tilde G_{1}^{l,m}}{2}} +  \inti{\abs{\tilde G_2^{l,m}}{2}} + \inti{\abs{\nabla^{m-1}\dt^l\tilde G_2}{2}},
	\end{aligned}
\end{equation}
where $ 0 < \omega < 1 $ would be addressed later. 
\end{lm}

\begin{pf}
	Take inner product of \subeqref{iPtlm}{2} with $ \partial^m \dt^l \iv $ and then record the resulting after integration by parts,
	\begin{align*}
			& \dfrac{d}{dt} \dfrac{1}{2} \inti{(\bar\rho +q) \abs{\partial^m \dt^l \iv }{2}} \underbrace{- \gamma \inti{ \dv(\bar\rho \partial^m\dt^l \iv)\bar\rho^{\gamma-2}\partial^m\dt^l\iq }}_{(i)} \\
			& ~~~~ + \mu \inti{\abs{\nabla \partial^m\dt^l \iv}{2} } + (\mu+\lambda) \inti{\abs{\dv \partial^m\dt^l \iv}{2}}\\
			& ~~ = \underbrace{ \dfrac{1}{2} \inti{\dt q \abs{\partial^m\dt^l \iv}{2}} + \dfrac{1}{2}\inti{\dv (\bar\rho v) \abs{\partial^m\dt^l \iv}{2}} }_{(ii)}\\
			& ~~~~ + \underbrace{ \inti{\tilde F_{2}^{l,m} \cdot \partial^m\dt^l\iv} }_{(iii)} + \underbrace{\inti{\tilde G_{2}^{l,m} \cdot \partial^m \dt^l \iv}}_{(iv)} + \underbrace{\inti{\partial^m\dt^l \tilde G_2 \cdot\partial^m\dt^l\iv}}_{(v)}.
	\end{align*}
Using \subeqref{iPtlm}{1}
\begin{align*}
		& (i) = \gamma\inti{(\dt \partial^m\dt^l \iq + (v+\bar u)  \cdot\nabla\partial^m\dt^l\iq - \tilde F_1^{l,m} - \tilde G_{1}^{l,m})\bar\rho^{\gamma-2}\partial^m\dt^l\iq} \\
		& = \dfrac{d}{dt} \dfrac{\gamma}{2} \inti{\bar\rho^{\gamma-2} \abs{\partial^m\dt^l \iq}{2}} - \dfrac{\gamma}{2} \underbrace{\inti{\dv(\bar\rho^{\gamma-2}(v+\bar u)) \abs{ \partial^m \dt^l \iq }{2}}}_{(vi)} \\
		& ~~~~ -\gamma \underbrace{ \inti{\tilde F_{1}^{l,m} \bar\rho^{\gamma-2} \partial^m\dt^l\iq} }_{(vii)} - \gamma\underbrace{\inti{\tilde G_1^{l,m}\bar\rho^{\gamma-2}\partial^m\dt^l\iq}}_{(viii)}.
\end{align*}
Now, we shall analyse $ \tilde F_1^{l,m}, \tilde F_2^{l,m} $. From the definition, it holds
\begin{align*}
		& \norm{\tilde F_1^{l,m}}{\Lnorm}^2 \lesssim (1+\bar\omega^a) \sum_{j=0}^{m} \inti{\abs{\nabla^j\dt^l v}{2}} + \bar\omega^a\sum_{j=0}^{m} \inti{\abs{\nabla^j\dt^l q}{2}}, \\
		& \norm{\tilde F_2^{l,m}}{\Lnorm}^2 \lesssim \norm{\partial^m\dt^l \tilde F_2 }{\Lnorm}^2 + \bar\omega^a \sum_{j=0}^{m-1} \inti{\abs{\partial^j\dt^{l+1} v}{2}} \\
		& ~~ + \bar\omega^a \sum_{j=0}^{m}\inti{\abs{\partial^j\dt^l q}{2}} \lesssim (1+\bar\omega^a ) \biggl( \inti{\abs{\nabla^{m+1}\dt^l v}{2}} + \inti{\abs{\nabla^m \dt^l q}{2}} \biggr)\\
		& ~~ +  \bar\omega^a \biggl(  \sum_{j=0}^{m-1} \inti{\abs{\nabla^j\dt^{l+1} v}{2}} + \sum_{j=0}^{m-1}\inti{\abs{\nabla^j\dt^l q}{2}} + \sum_{j=0}^{m}\inti{\abs{\nabla^j \dt^l v}{2}} \biggr).
\end{align*}
Therefore, Cauchy's inequality then yields, for any $ \omega > 0, \delta > 0 $, 
\begin{align*}
		& (ii),(vi) \lesssim (\norm{\dt q}{\supnorm} + \bar\omega \norm{v}{\supnorm} + \norm{\nabla v}{\supnorm} + \bar\omega^a ) \\ & ~~~~~ \times \sum_{j=0}^{m} \biggl( \inti{\abs{\nabla^j\dt^l v}{2}} + \inti{\abs{\nabla^j \dt^l q}{2}} \biggr),\\ 
		& (vii) \lesssim \omega \inti{\abs{\nabla^m \dt^l q}{2}} + C_\omega (1+\bar\omega^a)\sum_{j=0}^{m}\inti{\abs{\nabla^j \dt^l v}{2}}\\
		& ~~~~~~ + C_{\omega} \bar\omega^a \sum_{j=0}^{m} \inti{\abs{\nabla^j\dt^l q}{2}},  \\
		& (iii) \lesssim \omega \norm{\tilde F_2^{l,m}}{\Lnorm}^2 + C_\omega \sum_{j=0}^m\inti{\abs{\nabla^j \dt^l v}{2}},\\
		& (v) = -\inti{\partial^{m-1}\dt^l\tilde G_2 \cdot \partial^{m+1}\dt^l \iv} \lesssim \delta \inti{\abs{\nabla^{m+1}\dt^l\iv}{2}} \\
		& ~~~~~~ + C_\delta \inti{\abs{\nabla^{m-1}\dt^l \tilde G_2}{2}}, \\
		& (iv),(viii) \lesssim \omega \sum_{j=1}^{m}\inti{\abs{\nabla^j\dt^l q}{2}} + \sum_{j=0}^{m}\inti{\abs{\nabla^j\dt^l v}{2}} \\ & ~~~~~~ + C_\omega \inti{\abs{\tilde G_{1}^{l,m}}{2}} + \inti{\abs{\tilde G_{2}^{l,m}}{2}}.
\end{align*}
After summing up these estimates with an appropriately small $ \delta > 0 $, \eqref{ivplm-est} follows.
\end{pf}

\section{Synthesis and Asymptotic Stability}\label{sec:globalwellposed}

In this section, we shall chain the energy estimates obtained in the previous sections. In particular, the full energy dissipation and evolution would be achieved, which would capture the propagation of initial regularities as time grows up. Moreover, the decay of the corresponding energy functionals implies that the solution with small initial perturbations would converge to the stationary solution given in \eqref{SNS3}. 

\subsection{Synthesis}
To chain the energy estimates, we first introduce some notations.  
For $ l \geq 0 $, denote the energy and dissipation functionals
\begin{align*}
		& \mathfrak E_{l} = \sum_{i=0}^{l} \biggl\lbrace \inti{(\bar\rho+q)\abs{\dt^i v}{2}} + c\inti{\abs{\nabla\dt^i v + \nabla\dt^i v^T}{2}}  \\ &  ~~~~ + c\inti{\abs{\dv \dt^i v }{2}} + \inti{\bar\rho^{\gamma-2}\abs{\dt^i q}{2}} - c\inti{\bar\rho^{\gamma-2}\dv(\bar\rho \dt^i v)\dt^i q} \biggr\rbrace, \\
		& \mathfrak D_{l} = \sum_{i=0}^{l} \biggl\lbrace \inti{\abs{\nabla\dt^iv}{2}} + \inti{\abs{\dv \dt^i v}{2}} + \inti{(1+q)\abs{\dt^{i+1} v}{2}} \biggr\rbrace,
\end{align*}
where $ c > 0 $ is determined by \eqref{little-c}.
Also, the tangential and interior energy and dissipation functionals are denoted as, 
\begin{align*}
		& \mathfrak E_{l,m}^\tau = \sum_{i=0}^{l}\sum_{j=0}^{m} \biggl\lbrace \inti{((\bar\rho+q)\abs{\dT^j \dt^i v}{2}} + \inti{\bar\rho^{\gamma-2}\abs{\dT^j\dt^i q}{2}}  \\ &  ~~~~ + \inti{(\bar\rho+q) \abs{\nabla^j\dt^i \iv}{2}} + \inti{\bar\rho^{\gamma-2}\abs{\nabla^j\dt^i \iq}{2}} \biggr\rbrace, \\
		& \mathfrak D_{l,m}^\tau = \sum_{i=0}^{l}\sum_{j=0}^{m} \biggl\lbrace \inti{\abs{\nabla \dT^j\dt^i v}{2}} + \inti{\abs{\dv \dT^j \dt^i v}{2}}  \\ &  ~~~~ + \inti{\abs{\nabla^{j+1}\dt^i \iv}{2}} \biggr\rbrace,
\end{align*}
for $ m \geq 0 $.
The normal direction energy and dissipation functionals are represented by, for $ n \geq 1$, 
\begin{align*}
		& \mathfrak E_{l,m,n} = \sum_{i=0}^{l} \sum_{j=0}^{m} \sum_{k=1}^{n} \biggl\lbrace \inti{\abs{\dT^j\dN^k\dt^i q}{2}} \biggr\rbrace, \\
		& \mathfrak D_{l,m,n} = \sum_{i=0}^{l} \sum_{j=0}^{m} \sum_{k=1}^{n} \biggl\lbrace \inti{\abs{\dT^j\dN^k\dt^i q}{2}} + \inti{\abs{\dT^j\dN^k\dv\dt^i v}{2}} \\ &  ~~~~ + \inti{\abs{ \nabla^{k+1} \dT^j \dt^i v}{2}} + \inti{\abs{\nabla^k \dT^j \dt^i q}{2}} \biggr\rbrace.
\end{align*}
For $ n = 0 $, denote
\begin{equation*}
	\mathfrak D_{l,m,0} = \mathfrak D^\tau_{l,m}.
\end{equation*}
Additionally, let $ \mathfrak E_{l,m}, \mathfrak D_{l,m} $ be defined as
\begin{align*}
		& \mathfrak E_{l,m} = \mathfrak E^\tau_{l,m} + \sum_{j=1}^{m} \mathfrak E_{l,m-j,j}, \\
		& \mathfrak D_{l,m} = \mathfrak D^\tau_{l,m} + \sum_{j=1}^{m} \mathfrak D_{l,m-j,j},
\end{align*} 
with $ m \geq 1 $. Moreover, a comparable quantity of $ \mathfrak D_{l,m} $ can be written as
\begin{equation*}
		\mathfrak{\bar D}_{l,m} = \sum_{i=0}^{l} \sum_{j=1}^{m} \biggl( \inti{\abs{\nabla^{j+1} \dt^i v}{2}} + \inti{\abs{\nabla^j\dt^i q}{2}} \biggr).
\end{equation*}
Then it holds, 
\begin{equation}\label{eqlm}
	\mathfrak{\bar D}_{l,m} \lesssim \mathfrak D_{l,m}, ~ m \geq 1.
\end{equation}
Also, denote 
\begin{equation*}
	\mathfrak{\bar D}_{l,0} = \mathfrak D_{l}, ~ \mathfrak{\bar D}_{l+1, -1} = \mathfrak D_{l}.
\end{equation*}
Notice, the above quantities are monotone increasing in each index $ l,m,n $.
\begin{lm}
	\begin{equation}\label{l=L,m=1}
	\begin{aligned}
		& \dfrac{d}{dt} \left\lbrace \mathfrak E_l + \mathfrak E_{l,1} \right\rbrace + \mathfrak D_l + \mathfrak D_{l,1} \lesssim \Lambda_l \left( \mathfrak D_l + \mathfrak D_{l,1} \right) + \mathfrak G_{l}, 
	\end{aligned}
\end{equation}
where $ \mathfrak G_{l} $ is defined as
\begin{equation}\label{quadratic-l}
	\begin{aligned}
		& \mathfrak G_{l} = \sum_{i=0}^{l} \biggl( \inti{\abs{G_1^i}{2}} + \inti{\abs{G_2^i}{2}} + \inti{\abs{G_1^{i,1}}{2}} + \inti{\abs{G_2^{i,1}}{2}}  \\
		& ~~ + \inti{\abs{\tilde G_1^{i,1}}{2}} + \inti{\abs{\tilde G_2^{i,1}}{2}}+ \inti{\abs{\dt^i\tilde G_2}{2}}\\ & ~~~~  + \inti{\abs{G_{1,N}^{i}}{2}} + \inti{\abs{G_{2,N}^i}{2}} + \inti{\abs{G_{s}^{i}}{2}} \biggr).
	\end{aligned}
\end{equation}
\end{lm}
 
\begin{pf}

From \eqref{Proposition-time},
\begin{equation}\label{L}
	\begin{aligned}
		& \dfrac{d}{dt} \mathfrak E_{l} + \mathfrak D_l \lesssim (\omega + \Lambda_l ) \mathfrak D_{l,1} + \Lambda_l \mathfrak D_l \\
		& ~~~~~~ + ( 1+C_\omega) \sum_{i=0}^{l} \inti{\abs{G_1^i}{2}} + \sum_{i=0}^{l} \inti{\abs{G_2^i}{2}}.
	\end{aligned}
\end{equation}
Meanwhile, as the consequences of \eqref{Proposition-tangential} and \eqref{ivplm-est} with $ m = 1 $, by choosing an appropriately small $ \omega > 0 $,
\begin{align}
		& \dfrac{d}{dt} \mathfrak E_{l,1} + \mathfrak D_{l,1} \lesssim (1+\Lambda_l) \mathfrak D_l + \Lambda_l \mathfrak D_{l,1} + \sum_{i=0}^{l} \biggl(\inti{\abs{G_1^{i,1}}{2}} + \inti{\abs{G_2^{i,1}}{2}}  {\nonumber}\\
		& ~~~~ + \inti{\abs{G_2^i}{2}} + \inti{\abs{\tilde G_1^{i,1}}{2}} + \inti{\abs{\tilde G_2^{i,1}}{2}} + \inti{\abs{\dt^i \tilde G_2}{2}} {\label{L,1}} \\ & ~~~~  + \inti{\abs{G_{1,N}^{i}}{2}} + \inti{\abs{G_{2,N}^i}{2}} + \inti{\abs{G_{s}^{i}}{2}}\biggr). {\nonumber}
\end{align}
Then a sum of \eqref{L,1} and \eqref{L} with an appropriately small $ \omega $ in \eqref{L} yields \eqref{l=L,m=1}. 
\end{pf}

From \eqref{qlmn-est}, \eqref{dvlmn-est} and \eqref{qvlmn-est}, the following estimates holds with $ m + n \geq 2 , n \geq 1 $,
\begin{align}
		& \dfrac{d}{dt} \inti{\abs{\dT^m \dN^n \dt^l q}{2}} + \inti{\bar\rho^\gamma\abs{\dT^m\dN^n\dt^l q}{2}} \lesssim \Lambda_l \mathfrak{\bar D}_{l,m+n} {\nonumber}\\
		& ~~  + \mathfrak D_{l,m+1,n-1} +\mathfrak {\bar D}_{l+1, m+n-2} + (1+\Lambda_l ) \left( \mathfrak{\bar D}_{l,m+n-1} + \mathfrak{\bar D}_{l+1,m+n-3} \right) \\
		& ~~ + \inti{\abs{G_{1,N}^{l,m,n}}{2}}, {\nonumber} \\
		& \inti{\abs{\dT^m \dN^n \dv \dt^l v}{2}} \lesssim \inti{\abs{\dT^m\dN^n\dt^lq}{2}} + \mathfrak D_{l,m+1,n-1} + \mathfrak {\bar D}_{l+1, m+n-2} {\nonumber}\\
		& ~~ + ( 1+\Lambda_l ) \left( \mathfrak{\bar D}_{l,m+n-1} + \mathfrak{\bar D}_{l+1,m+n-3} \right) + \inti{\abs{G_{2,N}^{l,m,n}}{2}}, \\
		& \inti{\abs{\nabla^{n+1}\dT^m\dt^l v}{2}} + \inti{\abs{\nabla^{n}\dT^m\dt^l q}{2}} \lesssim \inti{\abs{\dT^m\nabla^{n}\dv \dt^l v}{2}}  {\nonumber}\\
		& ~~ + \mathfrak D^\tau_{l,m+n}  + \mathfrak {\bar D}_{l+1,m+n-2} + (1+\Lambda_l) \left( \mathfrak{\bar D}_{l+1,n+m-3} + \mathfrak{\bar D}_{l,m+n-1} \right) \\
		& ~~ + \norm{G_s^{l,m}}{\Hnorm{n-1}}^2 {\nonumber}.
\end{align}
Additionally, by using \eqref{DFrame04} and the commutator property \eqref{commt} repeatedly, 
\begin{equation*}
	\begin{aligned}
		& \inti{\abs{\dT^m \nabla^n \dv \dt^l v}{2}} \lesssim \inti{\abs{\dT^m \dN^n \dv \dt^l v}{2}} + \inti{\abs{\nabla^{n}\dT^{m+1} \dt^l v}{2}} \\
		& ~~ + \int_{B_{1/2}} \abs{\nabla^{m+n+1} \dt^l  v}{2} \,dx  + \sum_{j=1}^{m+n} \inti{\abs{\nabla^j \dt^l v}{2}} \lesssim \inti{\abs{\dT^m \dN^n \dv \dt^l v}{2}} \\
		& ~~ + \mathfrak D_{l,m+1,n-1} + \mathfrak D^\tau_{l,m+n} + \mathfrak{\bar D}_{l,m+n-1}.
	\end{aligned}
\end{equation*}
Therefore, chaining these inequalities together with \eqref{ql01-est}, \eqref{dvl01-est} and \eqref{qvl02-est} then yields, for $ m + n \geq 2 , n \geq 1 $, 
\begin{align}
		& \dfrac{d}{dt} \mathfrak E_{l,m,n} + \mathfrak D_{l,m,n} \lesssim \Lambda_l \mathfrak{\bar D}_{l,m+n} + \mathfrak D_{l,m+1,n-1}+ \mathfrak D^\tau_{l,m+n} {\nonumber}\\
		& ~~  + \mathfrak{\bar D}_{l+1, m + n-2}+ (1 + \Lambda_l ) \biggl( \mathfrak{\bar D}_{l,m+n-1} + \mathfrak{\bar D}_{l+1,m+n-3} \biggr) {\label{lmn}} \\
		& ~~ + \sum_{i=0}^{l} \sum_{j=0}^{m}\sum_{k=1}^{n} \biggl( \inti{\abs{G_{1,N}^{i,j,k}}{2}} + \inti{\abs{G_{2,N}^{i,j,k}}{2}} + \norm{G_s^{i,j}}{\Hnorm{k-1}}^2 \biggr).{\nonumber}
\end{align}
On the other hand, as the consequence of \eqref{vqlm-est} and \eqref{ivplm-est} together with \eqref{Prop:l-estimate} for $ m \geq 1 $, 
\begin{align}
		& \dfrac{d}{dt}\mathfrak E^\tau_{l,m} + \mathfrak D^\tau_{l,m} \lesssim ( \omega + C_\omega \bar\omega^a + \Lambda_l ) \mathfrak{\bar D}_{l,m} + ( \omega + C_\omega + \Lambda_l )\mathfrak{\bar D}_{l,m-1} {\nonumber} \\
		& ~~ + \omega\Lambda_l \mathfrak{\bar D}_{l+1,m-2} + C_\omega \sum_{i=0}^{l}\sum_{j=1}^{m} \biggl( \inti{\abs{G_1^{i,j}}{2}} + \inti{\abs{G_2^{i,j}}{2}}   {\label{taulm}}\\ 
		& ~~  + \inti{\abs{\dT^{j-1}G_2^i}{2}} +  \inti{\abs{\tilde G_1^{i,j}}{2}} + \inti{\abs{\tilde G_2^{i,j}}{2}} + \inti{\abs{\nabla^{j-1}\dt^i\tilde G_2}{2}} \biggr). {\nonumber}
\end{align}
Now we would be able to chain the total energy estimates. 
\begin{lm}
	For $ k \geq 2 $, 
	\begin{equation}\label{est-lk}
		\begin{aligned}
			& \dfrac{d}{dt} \mathfrak E_{l,k} + \mathfrak D_{l,k} \lesssim \Lambda_l \mathfrak D_{l,k} + (1+\Lambda_l)\mathfrak{\bar D}_{l+1,k-2} + ( 1 + \Lambda_l ) \mathfrak{\bar D}_{l,k-1} \\
			& ~~~~~ + (1+\Lambda_l) \mathfrak{\bar D}_{l+1,k-3} + \mathfrak G_{l,k},
		\end{aligned}
	\end{equation}
	where $ \mathfrak G_{l,k} $ is defined in \eqref{quadratic-lk}. In particular,
	\begin{align}
			& \dfrac{d}{dt} \mathfrak E_{l,2} + \mathfrak D_{l,2} \lesssim \Lambda_l \mathfrak D_{l,2} +  (1+\Lambda_l)\left( \mathfrak D_{l+1} + \mathfrak D_{l,1} \right) + \mathfrak G_{l,2} \label{est-l,2}, \\
			& \dfrac{d}{dt} \mathfrak E_{l,3} + \mathfrak D_{l,3} \lesssim \Lambda_l \mathfrak D_{l,3} + (1+\Lambda_l) \left( \mathfrak D_{l+1,1} + \mathfrak D_{l+1}+\mathfrak D_{l,2} \right) + \mathfrak G_{l,3}. \label{est-l,3}
	\end{align}
\end{lm}

\begin{pf}
	Consider \eqref{lmn} with $ m = k - j, n = j $ for $ j \geq 1 $,
	\begin{align*}
			& \dfrac{d}{dt} \mathfrak E_{l,k-j, j} + \mathfrak D_{l,k-j, j} \lesssim \Lambda_l \mathfrak{\bar D}_{l,k} + \mathfrak D_{l,k-j+1, j-1} \\
			& ~~ + \mathfrak D^\tau_{l,k} + \mathfrak{\bar D}_{l+1,k-2} + (1+\Lambda_l)\left(\mathfrak{\bar D}_{l,k-1} + \mathfrak{\bar D}_{l+1,k-3} \right) + \sum_{a=0}^{l}\sum_{b=0}^{k-j}\sum_{c=1}^{j} \mathfrak G_{a,b,c},
	\end{align*}
	where
	\begin{equation*}
		\begin{aligned}
			\mathfrak G_{a,b,c} = \inti{\abs{G_{1,N}^{a,b,c}}{2}} + \inti{\abs{G_{2,N}^{a,b,c}}{2}} + \norm{G_s^{a,b}}{\Hnorm{c-1}}^2.
		\end{aligned}
	\end{equation*}
	Thus, after summing from $ j = 1 $ to $ j = k $, together with \eqref{taulm} with $ m = k $ and an appropriately small $ \omega > 0 $, it holds
	\begin{equation*}
		\begin{aligned}
			& \dfrac{d}{dt} \mathfrak E_{l,k} + \mathfrak D_{l,k} \lesssim \Lambda_l \mathfrak D_{l,k} + (1+\Lambda_l)\mathfrak{\bar D}_{l+1,k-2} + ( 1 + \Lambda_l ) \mathfrak{\bar D }_{l,k-1} \\
			& ~~~~~ + (1+\Lambda_l) \mathfrak{\bar D}_{l+1,k-3} + \mathfrak G_{l,k},
		\end{aligned}
	\end{equation*}
	where
	\begin{equation}\label{quadratic-lk}
		\begin{aligned}
			& \mathfrak G_{l,k} = \sum_{j=1}^{k} \sum_{a=0}^{l}\sum_{b=0}^{k-j}\sum_{c=1}^{j} \mathfrak G_{a,b,c} + \sum_{a=0}^{l}\sum_{b=1}^{k} \biggl( \inti{\abs{G_1^{a,b}}{2}} + \inti{\abs{G_2^{a,b}}{2}} \\ 
			&  +\inti{\abs{\dT^{b-1}G_2^a}{2}} +  \inti{\abs{\tilde G_1^{a,b}}{2}} + \inti{\abs{\tilde G_2^{a,b}}{2}} + \inti{\abs{\nabla^{b-1}\dt^a\tilde G_2}{2}} \biggr).
		\end{aligned}
	\end{equation}
\end{pf}

Thus we have shown the following
\begin{prop} For $ L \geq 2 $, 
	\begin{equation}\label{est-total}
		\begin{aligned}
			& \dfrac{d}{dt} \mathfrak{\bar E}_L(t) + (1-\Lambda_L ) \mathfrak{\bar D}_L(t)  \lesssim \mathfrak G_L(t) + \sum_{i=1}^{L} \mathfrak G_{L-i,2i+1}(t),
		\end{aligned}
	\end{equation}
	where
	\begin{align}
		& \mathfrak{\bar E}_L(t) = \mathfrak E_L(t) + \mathfrak E_{L,1}(t) + \sum_{i=1}^{L} \left( \mathfrak E_{L-i,2i}(t) + \mathfrak E_{L-i,2i+1}(t) \right), {\label{total-Energy}}\\
		& \mathfrak{\bar D}_L(t) = \mathfrak D_{L}(t) + \mathfrak D_{L,1}(t) + \sum_{i=1}^{L} \left( \mathfrak D_{L-i,2i}(t) + \mathfrak D_{L-i,2i+1}(t) \right). {\label{total-Dissipation}}
	\end{align}
\end{prop}

\begin{pf}

	For $ i = 1 $, as the consequence of \eqref{est-l,2} and \eqref{est-l,3},
	\begin{equation}\label{est-L-1,2}
		\begin{aligned}
			& \dfrac{d}{dt}\left( \mathfrak E_{L-1,2} + \mathfrak E_{L-1,3} \right) + \mathfrak D_{L-1,2} + \mathfrak D_{L-1,3} \lesssim \Lambda_L \left(\mathfrak D_{L-1,2} + \mathfrak D_{L-1,3} \right) \\
			& ~~~~~~ + (1+\Lambda_L) \left(\mathfrak D_{L} + \mathfrak D_{L,1} \right) + \mathfrak G_{L-1,3}.
		\end{aligned}
	\end{equation}
	Similarly, for $ i \geq 2, 1 \leq j \leq i $, from \eqref{est-lk},
		\begin{align*}
			& \dfrac{d}{dt} \mathfrak E_{L-i,2j} + \mathfrak D_{L-i,2j} \lesssim \Lambda_L \mathfrak D_{L-i,2j} {\nonumber} \\
			& ~~~~ + (1+\Lambda_L) \left(\mathfrak D_{L-i+1,2j-2} + \mathfrak D_{L-i,2j-1} \right) +\mathfrak G_{L-i,2j}, \\
			& \dfrac{d}{dt} \mathfrak E_{L-i,2j+1} + \mathfrak D_{L-i,2j+1} \lesssim \Lambda_L \mathfrak D_{L-i,2j+1} {\nonumber} \\
			& ~~~~ + (1+\Lambda_L) \left( \mathfrak D_{L-i+1,2j-1} + \mathfrak D_{L-i,2j} \right) + \mathfrak G_{L-i,2j+1}.
		\end{align*}
	Hence, after summing over $ j $, it holds 
	\begin{equation}\label{est-L-i,2i}
		\begin{aligned}
			& \dfrac{d}{dt} \left( \mathfrak E_{L-i,2i} + \mathfrak E_{L-i,2i+1} \right) + \mathfrak D_{L-i,2i} + \mathfrak D_{L-i,2i+1} \lesssim \Lambda_L ( \mathfrak D_{L-i,2i} + \mathfrak D_{L-i,2i+1} ) \\
			& ~~ + (1+\Lambda_L) ( \mathfrak{D}_{L-i,1} + \mathfrak{D}_{L-i+1,2i-1}) + \mathfrak G_{L-i,2i+1}\lesssim \Lambda_L( \mathfrak D_{L-i,2i} \\
			& ~~ + \mathfrak D_{L-i,2i+1} )  + (1+\Lambda_L)(\mathfrak{D}_{L,1} + \mathfrak D_{L-i+1,2i-1}) + \mathfrak G_{L-i,2i+1}.
		\end{aligned}
	\end{equation}
	Then \eqref{est-total} is the consequence of \eqref{l=L,m=1} with $ l = L $, \eqref{est-L-1,2} and summing \eqref{est-L-i,2i} with $ i $ from $ 2 $ to $ L $.
\end{pf}

\subsection{Estimates on Nonlinearities}\label{sec:nonlinearity}

In this section, we shall perform the estimates on the nonlinearities. To do so, denote the intermediate energy and dissipation $ \mathcal E_L, \mathcal D_L $ as
	\begin{align}
		& \mathcal{\bar E}_L =\norm{\dt^L v}{\Hnorm{1}}^2 + \sum_{i=0}^{L-1} \norm{\dt^i v}{\Hnorm{2L-2i}}^2 + \sum_{i=0}^L \norm{\dt^i q}{\Hnorm{2L-2i+1}}^2,  {\label{total-energy}} \\
		& \mathcal{\bar D}_L = \sum_{i=0}^{L} \biggl( \norm{\dt^i v}{\Hnorm{2L-2i+2}}^2 + \norm{\dt^i q}{\Hnorm{2L-2i+1}}^2  {\label{total-dissipation}} \biggr) + \norm{\dt^{L+1}v}{\Lnorm}^2.
	\end{align}
Then it holds, 
	\begin{align}
		\sum_{2a+b\leq 2L-2}\norm{\nabla^b\dt^a v}{\supnorm}^2 + \sum_{2a+b\leq 2L-1} \norm{\nabla^b\dt^a q}{\supnorm}^2 & \lesssim \mathcal{\bar E}_L, {\label{sup-est}}\\
		\sum_{2a+b\leq 2L+2}\norm{\nabla^b\dt^a v}{\Lnorm}^2 + \sum_{2a+b\leq 2L+1}\norm{\nabla^b\dt^a q}{\Lnorm}^2 &\lesssim \mathcal{\bar D}_L. {\label{L2-est}} 
	\end{align}


We shall use \eqref{sup-est} and \eqref{L2-est} to manipulate the nonlinear terms on the right hand side of \eqref{est-total}.

\begin{lm} For $ L \geq 3 $, 
	\begin{equation}\label{est-Nonlinear}
		\begin{aligned}
			& \mathfrak G_L + \sum_{i=1}^{L} \mathfrak G_{L-i,2i+1} \lesssim \mathfrak P(\mathcal{\bar E}_L) \mathcal{\bar D}_L.
		\end{aligned}
	\end{equation}
	In the meantime,
	\begin{equation}\label{est-coefficience}
		\Lambda_L \lesssim \mathfrak P( \bar\omega, \mathcal{\bar E}_L).
	\end{equation}
	where $ \mathfrak P(\cdot) $ is a polynomial with the property $ \mathfrak P(0) = 0 $.
\end{lm}

\begin{pf}
	We only show $ \mathfrak G_L $, and others can be handled in a similar way. Recall,
	\begin{align*}
		& \mathfrak G_L =  \sum_{i=0}^{L} \biggl( \underbrace{ \inti{\abs{G_1^i}{2}} + \inti{\abs{G_2^i}{2}} }_{(i)} + \underbrace{\inti{\abs{G_1^{i,1}}{2}} + \inti{\abs{G_2^{i,1}}{2}} }_{(ii)}  \\
		& ~~ + \underbrace{ \inti{\abs{\tilde G_1^{i,1}}{2}} + \inti{\abs{\tilde G_2^{i,2}}{2}} + \inti{\abs{\dt^i\tilde G_2}{2}} }_{(iii)}\\ & ~~~~  + \underbrace{ \inti{\abs{G_{1,N}^{i}}{2}} + \inti{\abs{G_{2,N}^i}{2}} + \inti{\abs{G_{s}^{i}}{2}}}_{(iv)} \biggr).
	\end{align*}
	We shall consider the highest order terms only. In the following, let $ l.o.t $ denote the terms involving relatively lower order derivatives than others, which might be different from line by line. 
	To evaluate $(i)$, from \eqref{Nonlinear-1}, \eqref{Nonlinear-2},
	\begin{align*}
			& \abs{G_1^L}{} \lesssim \abs{\dt^L G_1}{} + \sum_{j=0}^{L-1}\abs{\dt^{L-j}v}{} \abs{\nabla\dt^j q}{} \lesssim \sum_{j=0}^{L} \abs{\dt^j q}{} \abs{\nabla \dt^{L-j} v}{} \\
			& ~~~~ + \sum_{j=0}^{L-1}\abs{\dt^{L-j}v}{} \abs{\nabla\dt^j q}{} \lesssim \sum_{j=0}^{[L/2]} \left( \abs{\dt^j q}{} + \abs{\nabla\dt^j v}{} \right) \left( \abs{\nabla\dt^{L-j} v}{} + \abs{\dt^{L-j} q}{} \right) \\
			& ~~~~ + \sum_{j=0}^{[L/2]} \left( \abs{\dt^{j+1} v}{} + \abs{\nabla\dt^j q}{} \right) \left( \abs{\nabla\dt^{L-j-1} q}{} + \abs{\dt^{L-j} v}{} \right).
	\end{align*}
	Similarly, 
	\begin{align*}
			& \abs{G_2^L}{} \lesssim \abs{\dt^L G_2}{} + \sum_{j=0}^{L-1} \left( \abs{\dt^{j+1}q}{}\abs{\dt^{L-j}v}{} + \abs{\dt^{j+1}v}{} \abs{\nabla\dt^{L-j-1}v}{} \right) \\
			& \lesssim \sum_{j=0}^{[L/2]} \left(\abs{\dt^j q}{} + \abs{\nabla\dt^j q}{} + \abs{\dt^j v}{} + \abs{\nabla\dt^j v}{} \right) \\ & ~~~~~~~~~~~~~~~ \times \left(\abs{\dt^{L-j} q}{} + \abs{\nabla\dt^{L-j} q }{} + \abs{\dt^{L-j} v}{} + \abs{\nabla\dt^{L-j} v}{} \right)\\
			& ~~~~ + \sum_{j=0}^{[L/2]} \left(\abs{\dt^{j+1}q}{} + \abs{\dt^{j+1}v}{}\right) \left( \abs{\dt^{L-j}v}{} + \abs{ \dt^{L-j} q}{} \right) \\
			& ~~~~ + \sum_{j=0}^{[L/2]} \left(\abs{\dt^{j+1}v}{} + \abs{\nabla\dt^{j}v}{} \right) \left( \abs{\nabla\dt^{L-j-1}v}{} + \abs{\dt^{L-j}v}{} \right) + \mathfrak C,
	\end{align*}
	where $ \mathfrak C $ denotes cubic terms. 
	Therefore, 
	\begin{align*}
			& (i) \lesssim \mathcal A \cdot \sum_{j=0}^{[L/2]} \biggl(\inti{\abs{\dt^{L-j}q}{2}} + \inti{\abs{\nabla\dt^{L-j}q}{2}} + \inti{\abs{\nabla\dt^{L-j-1}q}{2}} \\ 
			& ~~  + \inti{\abs{\dt^{L-j}v}{2}} + \inti{\abs{\nabla\dt^{L-j}v}{2}} + \inti{\abs{\nabla\dt^{L-j-1}v}{2}}   \biggr) ,
	\end{align*}
	where $ \mathcal A = \mathcal A(\cdot) $ is a polynomial of
	\begin{align*}
			& \biggl\lbrace \norm{\dt^j q}{\supnorm}, \norm{\dt^{j+1} q}{\supnorm}, \norm{\nabla\dt^j q}{\supnorm}  \\ & ~~~~ \norm{\dt^j v}{\supnorm}, \norm{\dt^{j+1}v}{\supnorm} , \norm{\nabla\dt^j v}{\supnorm}  \biggr\rbrace_{0\leq j\leq [L/2]},
	\end{align*}
	with the property $ \mathcal A(0) = 0 $. Thus from \eqref{sup-est} and \eqref{L2-est}, $ (i) \lesssim \mathfrak P(\mathcal{\bar E}_L) \mathcal{\bar D}_L $ provided
	\begin{equation*}
		\begin{aligned}
			& 2([L/2] + 1) \leq 2L-2, ~ \text{or equivalently} ~ L \geq 3.
		\end{aligned}
	\end{equation*}
	Similarly, from \eqref{Nonlinear-3}
	\begin{align*}
			& \abs{G_1^{L,1}}{} \lesssim \abs{\nabla G_1^L}{} + \abs{\nabla v}{} \abs{\nabla\dt^L q}{} + \abs{v}{} \abs{\nabla\dt^L q}{} \\
			& \lesssim \sum_{j=0}^{L} \left( \abs{\nabla\dt^j q}{} \abs{\nabla\dt^{L-j} v}{} + \abs{\dt^j q}{} \abs{\nabla^2\dt^{L-j}v}{} \right) \\
			& ~~~~ + \sum_{j=0}^{L-1} \left( \abs{\nabla\dt^{L-j}v}{} \abs{\nabla\dt^j q}{} + \abs{\dt^{L-j}v}{} \abs{\nabla^2 \dt^j q}{} \right) + \underbrace{ \abs{\nabla v}{} }_{L^\infty} \underbrace{ \abs{\nabla\dt^L q}{} }_{L^2}\\
			& ~~~~ + \underbrace{ \abs{v}{} }_{L^\infty} \underbrace{ \abs{\nabla\dt^L q}{} }_{L^2} + l.o.t, \\
			& \abs{G_2^{L,1}}{} \lesssim \underbrace{ \abs{\nabla q}{} }_{L^\infty} \underbrace{ \abs{\dt^{ L+1}v}{} }_{L^2} + \underbrace{ \abs{v}{} }_{L^\infty} \underbrace{ \abs{\nabla\dt^L v}{} }_{L^2} + \underbrace{ \left( \abs{\nabla v}{} + \abs{v}{} \right) }_{L^\infty} \underbrace{ \abs{\nabla\dt^L v}{}}_{L^2} + l.o.t,
	\end{align*}
	where the subscripts represent the corresponding norms to bound the terms. The same notations would be adopted in the following. Meanwhile, the rests in $ \abs{G_{1}^{L,1}}{} $ are bounded in a similar way. 
	\begin{align*}
			& \sum_{j=0}^{L} \abs{\nabla\dt^j q}{} \abs{\nabla\dt^{L-j} v}{} \lesssim \sum_{ \begin{array}{c}2j+1\leq 2L-1,\\ 2(L-j)+1\leq 2L+2 \end{array} } \abs{\nabla\dt^j q}{} \abs{\nabla\dt^{L-j} v}{} \\
			& ~~~ + \sum_{ \begin{array}{c} 2(L-j) + 1 \leq 2L-2,\\ 2j+1\leq 2L+1 \end{array}} \abs{\nabla\dt^{L-j} v}{}\abs{\nabla\dt^j q}{} = \sum_{0 \leq j \leq L-1 } \underbrace{\abs{\nabla\dt^j q}{}}_{L^\infty} \underbrace{\abs{\nabla\dt^{L-j} v}{}}_{L^2} \\
			& ~~~ + \sum_{3/2\leq j\leq L} \underbrace{\abs{\nabla\dt^{L-j} v}{}}_{L^\infty} \underbrace{\abs{\nabla\dt^j q}{}}_{L^2},
	\end{align*}
	provided
	\begin{equation*}
		3/2 \leq L-1, ~ \text{or equivalently} ~ L \geq 3. 
	\end{equation*}
	Similar arguments then yield,
	\begin{align*}
			& \sum_{j=0}^L \abs{\dt^j q}{} \abs{\nabla^2\dt^{L-j}v}{} \lesssim \sum_{0\leq j\leq L-1}  \underbrace { \abs{\dt^j q}{}}_{L^\infty} \underbrace{ \abs{\nabla^2\dt^{L-j}v}{} }_{L^2} + \sum_{2\leq j\leq L} \underbrace{ \abs{\nabla^2\dt^{L-j}v}{} }_{L^\infty} \underbrace{ \abs{\dt^j q}{} }_{L^2}, \\
			& \sum_{j=0}^{L-1} \abs{\nabla\dt^{L-j}v}{} \abs{\nabla\dt^j q}{} \lesssim \sum_{2\leq j\leq L} \underbrace{ \abs{\nabla\dt^{L-j}v}{}}_{L^\infty} \underbrace{ \abs{\nabla\dt^j q}{} }_{L^2} + \sum_{0 \leq j\leq L-1} \underbrace{ \abs{\nabla\dt^j q}{} }_{L^\infty} \underbrace{ \abs{\nabla\dt^{L-j}v}{} }_{L^2} ,   \\
			& \sum_{j=0}^{L-1} \abs{\dt^{L-j}v}{} \abs{\nabla^2 \dt^j q}{} \lesssim \sum_{1\leq j\leq L-1} \underbrace{ \abs{\dt^{L-j}v}{} }_{L^\infty} \underbrace { \abs{\nabla^2 \dt^j q}{} }_{L^2} + \sum_{0 \leq j\leq L-2} \underbrace{ \abs{\nabla^2 \dt^j q}{} }_{L^\infty} \underbrace { \abs{\dt^{L-j}v}{} }_{L^2},
	\end{align*}
	provided
	\begin{equation*}
		L \geq 3.
	\end{equation*}
	Thus, $ (ii) \lesssim \mathfrak P(\mathcal{\bar E}_L) \mathcal{\bar D}_L $.	Similarly $ (iii) \lesssim \mathfrak P(\mathcal{\bar E}_L) \mathcal{\bar D}_L $. To handle $ (iv) $, from \eqref{Nonlinear-4}, \eqref{Nonlinear-5}
	\begin{equation*}
		\begin{aligned}
			& \abs{G_{1,N}^L}{}, \abs{G_{2,N}^L}{}, \abs{G_s^L}{} \lesssim \abs{\nabla G_1^L}{} + \abs{\nabla v}{} \abs{\nabla\dt^L q}{} + \abs{v}{} \abs{\nabla\dt^L q}{}\\
			& ~~~~ + \abs{q}{} \abs{\dt^{L+1}v}{} + \abs{v}{} \abs{\nabla\dt^L v}{} + \abs{G_2^L}{},
		\end{aligned}
	\end{equation*}	
	in which the terms on the right hand side have already appeared before. Hence, $ (iv) \lesssim \mathfrak P(\mathcal{\bar E}_L) \mathcal{\bar D}_L $. We have shown $ \mathfrak G_L  \lesssim \mathfrak P(\mathcal{\bar E}_L) \mathcal{\bar D}_L$. The estimate on $ \Lambda_L $ is the direct consequence of \eqref{sup-est}.	
\end{pf}



Next lemma concerns the quantity analysis of $ \mathcal{\bar E}_L $ and $ \mathcal{\bar D}_L $. 
\begin{lm}It holds the equivalent relation 
	\begin{equation}\label{energycmpare-d}
		\begin{aligned}
			& \mathcal{\bar D}_L \simeq \mathfrak{\bar D}_L. 
		\end{aligned}
	\end{equation}
	In the meantime,
		\begin{gather}
			\norm{\dt^L v}{\Hnorm{1}}^2 + \sum_{i=0}^L \norm{\dt^i q}{\Hnorm{2L-2i+1}}^2 \lesssim \mathfrak{\bar E}_L(t) \lesssim \mathcal{\bar D}_L(t)  ,\label{energycmpare-e1} \\
			\mathcal{\bar E}_{L} \lesssim \mathcal{\bar D}_L. \label{enerycmpare-e3}
		\end{gather}
	Also, a direct calculation yields,
	\begin{equation}\label{energycmpare-e2}
		\begin{aligned}
			& \dfrac{d}{dt} \sum_{i=0}^{L-1}\norm{\dt^i v}{\Hnorm{2L-2i}}^2 \lesssim \sum_{i=0}^L \norm{\dt^i v}{\Hnorm{2L-2i+2}}^2 \lesssim \mathcal{\bar D}_{L}.
		\end{aligned}
	\end{equation}
	
\end{lm}

\begin{pf}
	As a consequence of the definitions, it is easy to show, 
	\begin{equation*}
		\begin{aligned}
			& \mathfrak D_L + \mathfrak D_{L,1} + \sum_{i=1}^L\left( \mathfrak D_{L-i,2i} + \mathfrak D_{L-i,2i+1} \right)\lesssim \mathcal{\bar D}_L.
		\end{aligned}
	\end{equation*}
	Meanwhile, from the estimate \eqref{DFrame04}, the commutator property \eqref{commt}, and the Poincar\'{e} inequality \eqref{qPoincare}, \eqref{vPoincare}, it can be shown without difficulty
	\begin{equation*}
		\begin{aligned}
			& \mathcal{\bar D}_L \lesssim \mathfrak D_L + \mathfrak D_{L,1} + \sum_{i=1}^L\left( \mathfrak D_{L-i,2i} + \mathfrak D_{L-i,2i+1} \right).
		\end{aligned}
	\end{equation*}
\eqref{energycmpare-e1}, \eqref{enerycmpare-e3} can be derived in a similar way. 	
\end{pf}

\subsection{Global Prior Estimate and Asymptotic Stability}
Let $ \mathfrak L(t) $ be defined as
\begin{equation*}
	\begin{aligned}
		& \mathfrak L(t) = \max_{0<s<t} \mathcal{\bar E}_L(t) + \int_{0}^{t} \mathcal{\bar D}_L(s)\,ds .
	\end{aligned}
\end{equation*}

\begin{lm}\label{lm:total-Energy}
	There is a $ \epsilon_0 > 0  $, such that $ \bar\omega \leq \epsilon_0 ,  \mathfrak L(t) \leq \epsilon_0 $ would imply
	\begin{equation}\label{est-decay1}
			 e^{\sigma t} \mathfrak{\bar E}_L(t) + \int_0^t e^{\sigma s} \mathfrak{\bar D}_L(s) \,ds \leq \mathfrak{\bar E}_L(0), \\
	\end{equation}
	for some $ \sigma >0 $.
\end{lm}

\begin{pf}
	From \eqref{est-total}, and \eqref{est-Nonlinear}, \eqref{est-coefficience} and the equivalence \eqref{energycmpare-d}, it holds
	\begin{equation*}
		\begin{aligned}
			& \dfrac{d}{dt}\mathfrak{\bar E}_L(t) + \mathfrak{\bar D}_L(t) \lesssim \mathfrak P (\bar\omega, \mathcal{\bar E}_L(t)) \left( \mathcal{\bar D}_L(t) + \mathfrak{\bar D}_L(t) \right) \lesssim \mathfrak P(\bar\omega,\mathfrak L(t) ) \mathfrak{\bar D}_L(t),
		\end{aligned}
	\end{equation*}
	where $ \mathfrak P = \mathfrak P(\cdot) $ is a polynomial with the property $ \mathfrak P(0) = 0$. Thus, $ \exists \epsilon_0 > 0 $ such that $ \mathfrak L \leq \epsilon_0, \bar\omega \leq \epsilon_0 $ would yield $\mathfrak P \leq \bar\epsilon $ with $\bar\epsilon > 0 $  small enough such that the following holds
	\begin{equation*}
		\dfrac{d}{dt} \mathfrak{\bar E}_L(t) + \mathfrak{\bar D}_L(t) \lesssim 0.
	\end{equation*}
	Multiply this inequality with $ e^{\sigma t} $, together with \eqref{energycmpare-e1},
	\begin{equation*}
		\dfrac{d}{dt} \left\lbrace e^{\sigma t} \mathfrak{\bar E}_L(t)\right\rbrace + e^{\sigma t} \mathfrak{\bar D}_L(t) \lesssim \sigma e^{\sigma t} \mathfrak{\bar E}_L(t) \lesssim \sigma e^{\sigma t} \mathfrak{\bar D}_L(t).
	\end{equation*}
	Thus \eqref{est-decay1} follows by integrating this inequality over temporal variable with an appropriately small $ \sigma > 0 $.
\end{pf}

Meanwhile, we shall demonstrate that the estimate \eqref{est-decay1} would in turn yield the boundedness of $ \mathfrak L(t) $. More precisely, 
\begin{lm}\label{lm:total-energy}
	There is a $ 0 < \epsilon_1 \leq \epsilon_0 $, such that under the assumption that $ \mathfrak{\bar E}_L(0) \leq \epsilon_1 $ and $ \mathcal{\bar E}_L(0) \leq \epsilon_1 $, \eqref{est-decay1} would imply
	\begin{equation}\label{est-decay2}
		e^{\sigma t}\mathcal{\bar E}_L(t) + \int_0^t e^{\sigma s} \mathcal{\bar D}_L(s) \,ds \leq \epsilon_0.
	\end{equation}
	In particular,
	\begin{equation}\label{est-decay3}
		\mathfrak L(t) \leq \epsilon_0.
	\end{equation}
\end{lm}

\begin{pf}
	From \eqref{energycmpare-e1} and \eqref{est-decay1}
	\begin{equation}\label{el1-est1}
		e^{\sigma t} \left( \norm{\dt^L v}{\Hnorm{1}}^2 + \sum_{i=0}^L \norm{\dt^i q}{\Hnorm{2L-2i+1}}^2 \right) \lesssim e^{\sigma t} \mathfrak{\bar E}(t) \leq \epsilon_1.
	\end{equation}
	Similarly, from \eqref{energycmpare-e2}, \eqref{enerycmpare-e3}, \eqref{energycmpare-d}.
	\begin{equation*}
		\begin{aligned}
			& \dfrac{d}{dt} \left( e^{\sigma t} \sum_{i=0}^{L-1}\norm{\dt^i v}{\Hnorm{2L-2i}}^2 \right) \lesssim \sigma  e^{\sigma t} \sum_{i=0}^{L-1}\norm{\dt^i v}{\Hnorm{2L-2i}}^2 + e^{\sigma t} \mathcal{\bar D}_L\\
			& ~~~~ \lesssim \sigma e^{\sigma t}\mathcal{\bar E}_L + e^{\sigma t} \mathcal{\bar D}_L \lesssim (\sigma+1) e^{\sigma t} \mathcal{\bar D}_L \lesssim (\sigma + 1 )e^{\sigma t} \mathfrak{\bar D}_L.
		\end{aligned}
	\end{equation*}
	Integration in the temporal variable then yields,
	\begin{equation}\label{el1-est2}
		e^{\sigma t} \sum_{i=0}^{L-1}\norm{\dt^i v}{\Hnorm{2L-2i}}^2  \lesssim \epsilon_1 + \mathcal{\bar E}_L(0) \lesssim \epsilon_1,
	\end{equation}
	as the consequence of \eqref{est-decay1}. Together with \eqref{el1-est1}, \eqref{energycmpare-d} and \eqref{est-decay1}, this then yields \eqref{est-decay2} with a sufficiently small $ \epsilon_1>0 $. 
	\eqref{est-decay3} is a direct consequence of \eqref{est-decay2}.
\end{pf}

\paragraph{Proof of Main Theorem}
Let $ \epsilon_0, \epsilon_1 > 0 $ be defined in Lemma \ref{lm:total-Energy} and Lemma \ref{lm:total-energy}. Then for $ \bar \omega < \epsilon_0 $, $ \mathfrak{\bar E}_L(0), \mathcal{\bar E}_L(0) < \epsilon_1 $ with $ L \geq 3 $, through continuous arguments, the following estimates hold
\begin{equation}
	\begin{aligned}
		&  e^{\sigma t}\mathcal{\bar E}_L(t) + \int_0^t e^{\sigma s} \mathcal{\bar D}_L(s) \,ds \leq \epsilon_0, \\
		& e^{\sigma t} \mathfrak{\bar E}_L(t) + \int_0^t e^{\sigma s} \mathfrak{\bar D}_L(s) \,ds \leq \epsilon_1.
	\end{aligned}
\end{equation}
In particular, the energy functional $ \mathfrak{\bar E}_L(t) + \mathcal{\bar E}_L(t) $ admits exponential decay as time grows up. This finishes the proof.


\paragraph{\bf Acknowledgements}
This work is part of the doctoral dissertation of the author under the supervision of Professor Zhouping Xin at the Institute of Mathematical Sciences of the Chinese University of Hong Kong, Hong Kong. The author would like to express great gratitude to Prof. Xin for his kindly support and professional guidance.



\end{document}